\documentclass[letterpaper, twoside, 11pt]{article}
\usepackage{newtxtext,newtxmath}
\DeclareMathAlphabet{\mathcal}{OMS}{cmsy}{m}{n}
\usepackage[textwidth=5.8in,textheight=9in,inner=1.35in]{geometry} 
\usepackage[T1]{fontenc}
\usepackage{amsmath,amssymb}
\usepackage{graphicx}
\usepackage{fancyhdr}
\usepackage{hyperref}

\newif\ifcourier

\ifcourier
\newgeometry{textwidth=6.5in,textheight=9.5in,inner=1in}
\usepackage{setspace} 

\newcommand{\mypagebreak}[1]  {}
\newcommand{\zzz}             {\\}
\newcommand{\ul}[1]           {\underline{#1}}
\newcommand{\bfx}[1]          {\underline{#1}}
\else
\newcommand{\mypagebreak}[1]  {\pagebreak#1}
\newcommand{\zzz}             {}
\newcommand{\ul}[1]           {\emph{#1}}
\newcommand{\bfx}[1]          {\textbf{#1}}
\fi

\makeatletter
\def\l@section{\@dottedtocline{1}{1em}{1.8em}}
\makeatother
\newcommand{\sectioncentred}[1]{\section[#1]{\centering#1}}
\newcommand{\sectioncentredunnumbered}[1]{%
  \section*{\hfil#1\hfil}
  \markboth{#1}{}
  \addcontentsline{toc}{section}{\hspace{1.8em}#1}
}

\newenvironment{mylist}{%
\renewcommand{\labelenumi}{\roman{enumi}.}
\vspace{-1.5ex}
\begin{enumerate}
\parskip2pt
}{%
\end{enumerate}
\vspace{-1ex}
}

\hypersetup{
    colorlinks=true,
    allcolors=[rgb]{0,0.2,0.3},
}

\makeatletter
\renewcommand\@makefntext[1]{\parindent .5em\makebox[0em][r]{\@makefnmark}~#1}
\makeatother
\newcounter{thisisdumb}

\hfuzz=\maxdimen
\makeatletter
\patchcmd\start@align{$$}{%
  $$%
  \displaywidth=\textwidth
  \displayindent=-\leftskip
}{}{\errmessage{Cannot patch \string\start@align}}
\patchcmd\mathdisplay{$$}{%
  $$%
  \displaywidth=\textwidth
  \displayindent=-\leftskip
}{}{\errmessage{Cannot patch \string\mathdisplay}}
\makeatother

\newcommand{\ppl}         {\ensuremath{\;{\scalebox{.8}{$\!\!\Big($}}}}
\newcommand{\ppr}         {\ensuremath{\;{\scalebox{.8}{$\!\!\Big)$}}}}
\newcommand{\PPL}         {\ensuremath{\;{\scalebox{.9}{$\!\!\Big($}}}}
\newcommand{\PPR}         {\ensuremath{\;{\scalebox{.9}{$\!\!\Big)$}}}}

\newcommand{\bb}[1]       {\ensuremath{\mathbb{#1}}}
\newcommand{\ssf}[1]      {\ensuremath{\textsf{#1}}} 
\newcommand{\sss}[1]      {\ensuremath{\mathcal{#1}}}
\newcommand{\rrm}[1]      {\ensuremath{\mathrm{#1}}}

\newcommand{\m}           {\ensuremath{\mkern 1mu}}
\newcommand{\mx}          {\ensuremath{\mkern -1mu}}
\newcommand{\n}           {\ensuremath{n}}
\newcommand{\w}           {\ensuremath{\omega}}
\newcommand{\bx}          {\ensuremath{\beta}}
\newcommand{\tx}          {\ensuremath{\theta}}

\newcommand{\Dx}          {\ensuremath{\Delta}}
\newcommand{\sx}          {\ensuremath{\sigma}}
\newcommand{\Sx}          {\ensuremath{\Sigma}}
\newcommand{\Xx}          {\ensuremath{\raisebox{1pt}{$\chi$}}}
\newcommand{\yx}          {\ensuremath{\psi}}
\newcommand{\px}          {\ensuremath{\varphi}}
\newcommand{\ee}          {\ensuremath{\varepsilon}}

\newcommand{\ax}          {\ensuremath{\land}}
\newcommand{\vx}          {\ensuremath{\m\lor}}
\newcommand{\nx}          {\ensuremath{\lnot\m}}
\newcommand{\Ax}          {\ensuremath{\forall}}
\newcommand{\Ex}          {\ensuremath{\exists\m}}
\newcommand{\imp}         {\ensuremath{\m\m\supset\m\m}}
\newcommand{\e}           {\ensuremath{\mathbin{\scalebox{.8}{$\in$}}}}
\newcommand{\notE}        {\ensuremath{\mathbin{\scalebox{.8}{$\notin$}}}}
\newcommand{\prece}       {\ensuremath{\preccurlyeq}}
\newcommand{\nprece}      {\ensuremath{\npreccurlyeq}}
\newcommand{\eqx}         {\ensuremath{\nolinebreak\mathbin{=}\nolinebreak}}

\newcommand{\gt}          {\mathrel{>}}
\newcommand{\lt}          {\mathrel{<}}
\newcommand{\gex}         {\ensuremath{\mathbin{\raisebox{.10ex}{\scalebox{.8}{$\ge$}}}}}
\newcommand{\lex}         {\ensuremath{\mathbin{\raisebox{.10ex}{\scalebox{.8}{$\le$}}}}}
\newcommand{\gtx}         {\ensuremath{\mathbin{\raisebox{.06ex}{\scalebox{.8}{$\gt$}}}}}
\newcommand{\ltx}         {\ensuremath{\mathbin{\raisebox{.06ex}{\scalebox{.8}{$\lt$}}}}}
\newcommand{\geX}         {\ensuremath{\mathbin{\raisebox{.05ex}{\scalebox{.9}{$\ge$}}}}}
\newcommand{\leX}         {\ensuremath{\mathbin{\raisebox{.05ex}{\scalebox{.9}{$\le$}}}}}
\newcommand{\gtX}         {\ensuremath{\mathbin{\raisebox{.02ex}{\scalebox{.9}{$\gt$}}}}}
\newcommand{\ltX}         {\ensuremath{\mathbin{\raisebox{.02ex}{\scalebox{.9}{$\lt$}}}}}

\newcommand{\provex}      {\ensuremath{\;{\scalebox{1.3}{$\vdash$}}\;}}
\newcommand{\nprovex}     {\ensuremath{\,\nvdash\,}}
\newcommand{\modelx}      {\ensuremath{\,\vDash\,}}

\newcommand{\A}           {\ensuremath{\sss A}}
\newcommand{\B}           {\ensuremath{\sss B}}
\newcommand{\M}           {\ensuremath{\resizebox{.875\width}{\height}{\sss M}}}
\newcommand{\I}           {\ensuremath{\resizebox{.750\width}{\height}{\sss I}}}
\newcommand{\Lw}          {\ensuremath{\sss L_{\infty\w}}}
\newcommand{\smallerP}    {\ensuremath{\mbox{{{\scalebox{.85}{$P$}}}}}}
\newcommand{\aaa}         {\ensuremath{\alpha}}
\newcommand{\bbb}         {\ensuremath{\beta}}

\newcommand{\AC}          {\ssf{AC}}
\newcommand{\wAC}         {\ssf{wAC}}
\newcommand{\BI}          {\ssf{BI}}
\newcommand{\cBI}         {\ssf{cBI}}
\newcommand{\CA}          {\ssf{CA}}
\newcommand{\DC}          {\ssf{DC}}
\newcommand{\KP}          {\ssf{KP}}
\newcommand{\GB}          {\ssf{GB}}
\newcommand{\HF}          {\rrm{HF}}
\newcommand{\PA}          {\ssf{PA}}
\newcommand{\PC}          {\ssf{PC}}
\newcommand{\PRA}         {\ssf{PRA}}
\newcommand{\s}           {\ssf{S}} 
\newcommand{\T}           {\ssf{T}}
\newcommand{\WKL}         {\ssf{WKL}}
\newcommand{\ZFC}         {\ssf{ZFC}}
\newcommand{\ZF}          {\ssf{ZF}}
\newcommand{\HY}          {\ssf{-}} 

\newcommand{\Pair}        {\ssf{Pair}}
\newcommand{\Union}       {\ssf{Union}}
\newcommand{\Induction}   {\ssf{Induction}}
\newcommand{\Collection}  {\ssf{Collection}}
\newcommand{\Bounding}    {\ssf{Bounding}}
\newcommand{\Separation}  {\ssf{Separation}}
\newcommand{\Replacement} {\ssf{Replacement}}
\newcommand{\Infinity}    {\ssf{Infinity}}
\newcommand{\Foundation}  {\ssf{Foundation}}
\newcommand{\Power}       {\ssf{Power}}
\newcommand{\Choice}      {\ssf{Choice}}

\newcommand{\Sat}         {\rrm{Sat}}
\newcommand{\Con}         {\rrm{Con}}
\newcommand{\rank}        {\rrm{rank}}
\newcommand{\rk}          {\rrm{rk}}
\newcommand{\sub}         {\rrm{sub}}
\newcommand{\wf}          {\rrm{wf}}
\newcommand{\ful}         {\ssf{ful}}
\newcommand{\Ful}         {\ssf{Ful}}
\newcommand{\FUL}         {\ssf{FUL}}
\newcommand{\rfn}         {\ssf{rfn}}
\newcommand{\Rfn}         {\ssf{Rfn}}
\newcommand{\RFN}         {\ssf{RFN}}
\newcommand{\Th}          {\rrm{Th}}
\newcommand{\True}        {\rrm{True}}
\newcommand{\good}        {\rrm{good}}
\newcommand{\domain}      {\rrm{domain}}

\newcommand{\Tr}          {\ensuremath{\mathrm{Tr}\,}}
\newcommand{\Prx}[1]      {\ensuremath{\mathrm{Pr}_{_{#1}}}} 
\newcommand{\isProof}[1]  {\ensuremath{\mathrm{Proof}_{_{#1}}}} 

\newcommand{\Myset}[2]    {\ensuremath{\{\:#1\,:\,#2\:\}}} 
\newcommand{\myset}[2]    {\ensuremath{\{\,#1  :  #2\,\}}} 
\newcommand{\mymap}[2]    {\ensuremath{\mx\mx:\mx#1\rightarrow#2}}
\newcommand{\godel}[1]    {\ensuremath{\ulcorner\mx\mx\mx#1\mx\urcorner}}
\newcommand{\lenm}[1]     {\ensuremath{\m|#1|^\text{--}}}
\newcommand{\lenx}[2]     {\ensuremath{\lenm{#1}\mx\mx\mx\eqx #2}}
\newcommand{\lenn}        {\ensuremath{\lenx\sx n}}
\newcommand{\seq}[1]      {\ensuremath{\langle#1\rangle}}
\newcommand*{\cupx}       {\mathbin{\raisebox{-.2ex}{\scalebox{1.4}{\ensuremath{\cup}}}}}
\newcommand*{\capx}       {\mathbin{\raisebox{-.2ex}{\scalebox{1.4}{\ensuremath{\cap}}}}}
\newcommand*{\aax}        {\mathbin{\scalebox{1.3}{\ensuremath{\bigwedge\nolimits}}}}
\newcommand*{\vvx}        {\mathbin{\scalebox{1.3}{\ensuremath{\bigvee  \nolimits}}}}
\newcommand{\up}          {\ensuremath{\raisebox{1.2pt}{$\upharpoonright$}}}
\newcommand{\phs}         {\ensuremath{\smash{\xrightarrow[{\raisebox{0ex}[0pt]{$\:^*\:$}}]{}}}}
\newcommand{\done}        {\ensuremath{\hfill\square}}
\newcommand{\circlex}[1]  {\ensuremath{\raisebox{-.45ex}{\scalebox{1.5}{\textcircled{\raisebox{.25ex}{$\scriptscriptstyle #1$}}}}}}
\newcommand{\ol}[1]       {\ensuremath{\overline{#1}}}
\newcommand{\thm}[1]      {\refstepcounter{thisisdumb}\noindent\bfx{#1}\;}

\newcommand{\myolscript}[1] {\vphantom{\ol -}\smash{\scriptscriptstyle{\ol{#1}}}}

\makeatletter
\newcommand*\MY@rightharpoonupfill@{%
    \arrowfill@\relbar\relbar\rightharpoonup
}
\newcommand*\overrightharpoon{%
    \mathpalette{\overarrow@\MY@rightharpoonupfill@}%
}
\newcommand*\@dblsty@mathpalette[2]{%
    \mathchoice
        {#1\displaystyle       \scriptstyle       {#2\,}}%
        {#1\textstyle          \scriptstyle       {#2\,}}%
        {#1\scriptstyle        \scriptscriptstyle {#2\,}}%
        {#1\scriptscriptstyle  \scriptscriptstyle {#2\,}}%
}
\newcommand*\@dblsty@overarrow@[4]{%
    \vbox{\ialign{##\crcr
        #1#3\crcr
        \noalign{\nointerlineskip}%
        $\m@th\hfil #2#4\hfil$\crcr
    }}%
}
\newcommand*\smalloverrightharpoon{%
    \@dblsty@mathpalette{\@dblsty@overarrow@\MY@rightharpoonupfill@}%
}
\makeatother
\newcommand{\vect}[1]     {\ensuremath{\smalloverrightharpoon{#1}}}
\newcommand{\vecth}[1]    {\ensuremath{\smalloverrightharpoon{\vphantom{\scalebox{.9}{d}}\smash{#1}}}}
\newcommand{\xvect}[1]    {\ensuremath{\m\m\vect{#1}}}
\newcommand{\xvecth}[1]   {\ensuremath{\m\m\vecth{#1}}}

\newcommand{\Av}          {\ensuremath{\Ax\m\m\vect}}
\newcommand{\Ev}          {\ensuremath{\Ex\m\vect}}
\newcommand{\ve}          {\ensuremath{\!\e}}
\newcommand{\subvect}[1]  {\vphantom{\xvect i}\smash{\m\scriptstyle{\xvect{#1}}}}

\begin{document}

\ifcourier
{\fontfamily{qcr}\selectfont
\begin{small}
\fi

\pagenumbering{roman}
\begin{titlepage}
\begin{center}
\ifcourier
\vspace*{.5in}
\else
\vspace*{1.75in}
\fi
\textbf{Some Problems in Logic:} \\
\bigskip
\textbf{APPLICATIONS OF KRIPKE'S NOTION OF FULFILMENT} \\
\vspace{.5in}
Joseph Emerson Quinsey \\
\vspace{1.75in}
Submitted in fulfilment of the \\
requirements for the degree of \\
Doctor of Philosophy \\
\vspace{1in}
St Catherine's College \\
\medskip
Oxford \\
\medskip
April 1980
\end{center}
\vfill
\noindent
Edited April 2019: The original thesis was type-written with hand-drawn symbols. This is a transcription of into LaTeX. No changes have been made, except for the correction of some typos, and the replacing of the Peano-Russell dot notation in formulas by parentheses.
\end{titlepage}

\clearpage
\null
\vfill
\centerline{Version date: \today}
\cleardoublepage

\ifcourier
\setstretch{0.9}
\else
\newgeometry{textwidth=6.2in,textheight=9.25in,inner=1.15in}
\fi
\thispagestyle{empty}
\begin{center}
SOME PROBLEMS IN LOGIC \\
Joseph Emerson Quinsey \\
St Catherine's College \\
Oxford \\
Hilary Term, 1980
\end{center}

\noindent
This work is a study of S. Kripke's notion of fulfilment. Motivated \zzz
by the result of [\ref{ref5}], Kripke was looking for a proof of G\"odel's Incompleteness \zzz
Theorem which was model-theoretic, natural (that is, without self-reference), \zzz
and easy. The resulting notion of fulfilment is very simple \zzz
and it provides a versatile tool for deriving a large number of results, \zzz
both new and old (indeed, the oldest), and in both Proof Theory and Model \zzz
Theory.

\noindent
In Chapter~I, we give short and elegant proofs to a number of \zzz
classical results; most of these are due to Kripke. There are two new \zzz
results. One, due to Kripke, is that there exists an easily definable \zzz
subring $R$ of the ring of primitive recursive functions such that for any \zzz
non-principal ultrafilter $D$ on \w, $R/D$ is a recursively saturated \zzz
model of Peano arithmetic. The other is that for any r.e.\ theory \T\ extending \zzz
\PRA\ and for any given r.e.\ set, we can feasibly find a $\Sx_1^0$ formula \zzz
which semi-represents in \T\ the given set; and if \T\ is not $\Sx_1^0$-sound,  \zzz
we can choose the formula to be $\Dx_1^0\,(\T)$.

\noindent
Chapter~II contains two distinct results. One answers a problem of \zzz
{[\ref{ref3}]} by showing that
\[
\Myset{\godel\phi\e\Pi^0_k}{\phi\text{ is }\Sx^0_k\text{-conservative over \PA}} \label{eq:*}\tag{$\ast$}
\]
is a complete $\Pi^0_2$ set. The second, when combined with the results of III, \zzz
gives a version of Herbrand's Theorem and the relationship between the \zzz
notions of \ul{proof} and \ul{fulfilment}.

\noindent
III gives an exposition and extension of the Hilbert-Ackermann method \zzz
of proving the consistency of \PA; our account is largely based on that \zzz
of [\ref{ref6}].

\noindent
IV is an exposition and extension of the method in [\ref{ref1}] for obtaining \zzz
conservation results of the form: $\Sx^1_2\HY\AC$ is $\Pi^1_3$-conservative over $(\Pi^1_1\HY\CA)_{<\ee_0}\up$. We \zzz
give a general version of this and from it derive a number of \zzz
results for arithmetic, analysis, and set theory. Using III, we may also \zzz
obtain \ul{uniform} versions: e.g.
\[
\Ax\aaa<\ee_0\;(\Pi^1_1\HY\CA)_{\aaa}\up\,\provex\RFN_{\Pi^1_3}({\Sx^1_2}\HY\AC)
\]
\noindent
In V and VI we give some model-theoretic applications of fulfilment. \zzz
V deals with non-\w-models and is based on [\ref{ref4}]. We also prove an extension \zzz
of the theorem of D. Scott and [\ref{ref2}] involving Weak K\"onig's Lemma, and we \zzz
describe the order types of some sets of elementary initial segments \zzz
of recursively saturated models of certain theories. VI is an extension of \zzz
{[\ref{ref2}]}'s theorem on minimal models of analysis. We develop the notion of indicator \zzz
for countable fragments of \Lw, and obtain a close parallel between \zzz
this and the first-order case. The chapter concludes with some \zzz
representability results in \w-logic and the analogue of \eqref{eq:*}.

\noindent
VII gives an exposition of the Paris-Harrington statement, deriving the \zzz
sharp negative results established by others, and we also give a strengthened \zzz
version of a key combinatorial lemma of [\ref{ref5}].
\vfill
{
\begin{enumerate}
\itemsep-.8ex\parsep0pt\topsep0pt\partopsep0pt
\ifcourier
\else
\itemindent1.6em
\fi
\renewcommand{\labelenumi}{[\arabic{enumi}]}
\item\label{ref1} H. Friedman, Iterated inductive definitions and $\Sx^1_2$-AC. \ul{Buffalo Conf.} 1970
\item\label{ref2} ---------------, Countable models of set theories, \ul{Cambrid\smash ge Summer School} 1973
\item\label{ref3} D. Guaspari, Partially conservative extensions of arithmetic, \ul{TAMS} 1979
\item\label{ref4} L. Kirby \& J. Paris, Initial segments of models of PA, \ul{LNM} 619 1977
\item\label{ref5} J. Paris \& L. Harrington, A mathematical incompleteness in PA, \ul{Handbook} 1977
\item\label{ref6} T. M. Scanlon, The consistency of number theory via Herbrand's Thm, \ul{JSL} 1973
\end{enumerate}
}

\ifcourier
\doublespacing
\else
\restoregeometry
\fi
\cleardoublepage

\newgeometry{textwidth=5in,textheight=9in,inner=1.75in}
\renewcommand*\contentsname{\hfil Table of Contents\hfil}
\addtocontents{toc}{\protect\vspace{2em}}
\tableofcontents
\restoregeometry
\cleardoublepage

\pagestyle{fancy}
\pagenumbering{arabic}
\thispagestyle{plain} \sectioncentredunnumbered{Introduction}

    The notion of fulfilment\footnote{The term ``fulfillability'' for
``the notion of fulfilment'' has the advantage of being one word
as opposed to four. But I will not use it,
and instead shall usually abbreviate ``the notion of fulfilment'' by
``fulfilment''.} is
due to S.~Kripke\footnote{I am greatly indebted to Professor Kripke
for informing me of some of his results. In Chapter~I,
results \hyperref[Chapter1:1.3]{1.3}, \hyperref[Chapter1:1.4]{1.4}, \hyperref[Chapter1:1.5]{1.5}, \hyperref[Chapter1:1.7]{1.7}, \hyperref[Chapter1:1.8]{1.8}, and \hyperref[Chapter1:1.13]{1.13} are due to
him, as is the \hyperref[Chapter2:FULRFN]{equivalence} $\RFN(\T)\equiv\FUL(\T)$ given
in Chapter~II. Professor Kripke
also informed me of the possibility of the exposition of the
\hyperref[Chapter3:HerbrandAckermannScanlon]{Herbrand-Ackermann-Scanlon} method given in Chapter~III.}. Motivated by the \zzz
Paris-Harrington result, he was looking for a proof of G\"odel's Incompleteness \zzz
Theorem which was model-theoretic, natural (that is, without self-reference) \zzz
and easy. The notion of fulfilment is very simple, and it \zzz
is implicit in much of mathematical logic. (Indeed, it is more or less \zzz
stated in some of the early works of Skolem, Herbrand and G\"odel.) In \zzz
making it explicit, we obtain a unifying notion with applications in \zzz
both Proof Theory and Model Theory.

    Let us first consider model-theoretic proofs of proof-theoretic results. \zzz
Fulfilment enables us to give enlightening proofs of some classical \zzz
results---proofs which do not use the fixed-point theorem. Included here are: \zzz
\hyperref[Chapter1:1.4]{1.4} Peano's arithmetic, \PA, is not finitely axiomatizable; \,\hyperref[Chapter1:1.5]{1.5}\, \PA\ \zzz
is not complete; and \hyperref[Chapter1:1.10]{1.10} the theory of \bb N is not arithmetical. To \zzz
dispense with the fixed-point theorem, however, we must occasionally pay the \zzz
price by considering only theories which include enough \Induction. This \zzz
is well illustrated by \hyperref[Chapter2:2.7]{2.7} where we give two proofs of the Essential \zzz
Unboundedness Theorem of G.~Kreisel and A.~Levy\,\cite{Krei68}. Let us consider our \zzz
proof \hyperref[Chapter1:1.5]{1.5} of the Incompleteness Theorem. When we apply it to theories \zzz
not in the language of arithmetic, say in the language of set theory or \zzz
analysis, the independent sentences we first obtain are quite complex as \zzz
measured in terms of the number and types of quantifiers, e.g.\ $\Pi_2$. But \zzz
by weakening the notion of fulfilment, we may obtain independent $\Pi^0_2$ \zzz
sentences. By weakening the notion still further \ul{and} using the fixed-point \zzz
theorem, we may obtain independent $\Pi^0_1$ sentences.

    Next let us consider some results of a more recursion-theoretic flavour. \zzz
Given any sufficiently strong r.e.\ theory \T\ and any $\Sx^0_1$ formula $\tx x$,\footnotemark\ \zzz
we show in \hyperref[Chapter1:1.12]{1.12} how to \ul{feasibly} (without the fixed-point theorem) \zzz
semi-represent the r.e.\ set $\myset{n\e\w}{\tx n}$ in \T\ by a $\Sx^0_1$ formula $\yx\mx x$; \zzz
indeed, if \T\ is not $\Sx^0_1$-sound, we may choose \yx\ to be $\Dx^0_1\,(\T)$. In \zzz
\hyperref[Chapter2:2.8]{2.8} we solve an open problem of D.~Guaspari\,\cite{Guas79} and R.~Solovay by showing
\[
\Myset{\tx\e\Pi^0_k}{\tx\text{ is }\Sx^0_k\text{-conservative over }\PA}
\]
to be a complete $\Pi^0_2$ set for each $k\ge 2$; we also give a complete (in \zzz
terms of quantifier complexity) classification of the analogous sets for \zzz
theories other than \PA.
\footnotetext{Footnote added 2019: The notation here means that \tx\ is a formula in the
language of \T\ which can be interpreted as being $\Sx^0_1$ in the language of arithmetic.
Also recall \ul{semi-represents} here means that for all \n~\e~\w,
$\bb N\modelx\tx n$ iff $\T\provex\yx\ol{n}$.}

    Thirdly, consider some proof-theoretic applications of fulfilment. \zzz
In Chapter~III we give an exposition and slight generalization of the \zzz
Hilbert-Ackermann method of proving the consistency of theories such as \zzz
\PA; our account is largely based on that of T.~Scanlon\,\cite{Scan73}. From our \zzz
result we may easily derive Herbrand's Theorem \hyperref[Chapter2:2.3]{2.3}, \hyperref[Chapter2:2.5]{2.5} and the \zzz
Reflexiveness Theorem \hyperref[Chapter2:2.4]{2.4}:
\[
\s+\Induction\provex\Ax n\:\big(\Prx\PC(\godel{\px\,\ol n\m\m})\imp\px\,n\big)\,,
\]
where \s\ is any weak theory which is sufficient to perform the required \zzz
coding, $\px\,n$ is any formula in the language of \s\ with a single free \zzz
variable of integer type, and the schema of \Induction\ ranges over all \zzz
formulae in the language of \s. Also from our proof we may obtain the \zzz
sharp bounds of G.~Minc\,\cite{Minc71} concerning subsystems of arithmetic, and \zzz
the ``No-Counter-Example Interpretation'' of G.~Kriesel. Our result also gives, \zzz
for example, the following. In \KP\ + \Infinity\ we may define by a \zzz
$\Dx_0$ predicate an ordering which is intuitively a well-ordering of order-type \zzz
$\ee_{On+1}$, the least $\ee$-number greater than the class of ordinals. Then \zzz
\KP\ plus $V\!=\!L$ plus the schema of \Foundation\ on this ordering implies \zzz
the schema of Uniform Reflection for \KP. We conclude Chapter~III with \zzz
some conservation results for the schema of \Induction\ over various \zzz
theories of analysis and set theory, and for the language of arithmetic \zzz
augmented with an extra constant $c$, conservation results for the schema \zzz
of \Induction\ \ul{up to $c$} over various theories of arithmetic.

    In Chapter~IV we give an exposition and extension of H.~Friedman's\,\cite{Frie70} \zzz
method of obtaining conservation results of the form, e.g.\ $\Sx^1_2\HY\AC\up$ is \zzz
$\Pi^1_3$-conservative over $\Pi^1_1\HY\CA\up$, and $\Sx^1_2\HY\AC$ is $\Pi^1_3$-conservative over \zzz
$(\Pi^1_1\HY\CA)_{<\ee_0}\up$, that is, axioms asserting that the relativized hyperjump may \zzz
be iterated \aaa\ times for all $\aaa<\ee_0$ but without the schema of \zzz
\Induction. The new results here are uniform versions of the above: for \zzz
example
\begin{align*}
\Ax n\,(\Pi^1_1\HY\CA)_n\up\,&\provex\RFN_{\Pi^1_3}(\Sx^1_2\HY\AC\up\,)\\
\Ax \aaa<\ee_0\,(\Pi^1_1\HY\CA)_{\aaa}\up\,&\provex\RFN_{\Pi^1_3}(\Sx^1_2\HY\AC)
\end{align*}
where e.g.\ $\Ax n\,(\Pi^1_1\HY\CA)_n\up$ is an axiom asserting that the relativized \zzz
hyperjump may be iterated \n\ times for each integer \n. We give a large \zzz
number of applications of our general theorems to set theory, analysis and \zzz
arithmetic; these results are for the most part known.

    Next let us consider some model-theoretic applications. Theorem \hyperref[Chapter1:1.13]{1.13}, \zzz
due to S.~Kripke, is interesting and very simple: there exists an easily \zzz
definable subring \sss F of the ring of rudimentary (or primitive recursive, or \zzz
recursive, etc.) functions such that for any non-principle ultrafilter \zzz
$D$ on \w, $\sss F/D$ is a recursively saturated model of \PA. Chapter~V is \zzz
concerned with non-\w-models of first-order theories. (This is a rapidly \zzz
expanding field, especially with regard to countable models of arithmetic, \zzz
and undoubtably many of the results contained herein are known to other \zzz
workers in the field. I have tried to give complete references, and, for \zzz
my own results, acknowledge any work independent of my own. But let me \zzz
say that most of the results of V were inspired by, or are due to, J.~Paris \zzz
and L.~Kirby.) We first note that many results concerning nonstandard \zzz
models of arithmetic have nothing whatsoever to do with arithmetic, so we \zzz
for the most part give our results in more generality; that is, they apply \zzz
to any model with sufficient coding abilities and satisfying the \zzz
appropriate \Collection\ axioms. Theorem \hyperref[Chapter5:5.3]{5.3} gives necessary and sufficient \zzz
conditions for the existence of certain initial segments which are models of \zzz
\ul{coded} theories, and \ul{indicators} for the same. We use this to give necessary \zzz
and sufficient conditions for the existence of initial segments which are \zzz
models of a given \ul{complete} theory. In \hyperref[Chapter5:5.11]{5.11} we give a description of the \zzz
order-types of certain sets of elementary initial segments of a recursively \zzz
saturated model of \PA. In \hyperref[Chapter5:5.9]{5.9} we give an extension of a result of L.~Kirby, \zzz
K.~McAloon and R.~Murawski concerning indicators in models of arithmetic for \zzz
models of analysis; the lemma \hyperref[Chapter5:5.9]{5.9} required is a common generalization of \zzz
a theorem of D.~Scott and H.~Friedman and the generalization of the \zzz
MacDowell-Specker Theorem (for countable models) by R.G.~Phillips and \zzz
H.~Gaifman. It concerns finding extensions of a countable model $M$ of arithmetic \zzz
which code a given countable class \sss X of subsets of $M$; the main requirement \zzz
is that $\seq{M,\sss X}$ satisfy \WKL: every infinite binary tree has an \zzz
infinite branch.

    In Chapter~VI we consider analogues for the results of V for \w-models, \zzz
or more generally, models of theories contained in countable fragments of \zzz
the infinitary language \Lw. We start by considering an extension of \zzz
a result of H.~Friedman\,\cite{Frie73}: there are no minimal models of analysis, where \zzz
analysis is assumed to contain the full schema of \AC. We present a \zzz
proof of this, one which is essentially that of Friedman\,\cite{Frie73}, but simpler \zzz
and more general. The main idea is that in a non-\bx-model we can construct \zzz
trees which are well-founded inside the model but not in the real world; \zzz
the desired substructure is then obtained from any infinite branch of such \zzz
a tree. Then we shall modify this construction using the notion of fulfilment \zzz
to obtain the following improvement: there are no minimal models of \zzz
$\Sx^1_1\HY\BI$. Indeed, in \hyperref[Chapter6:6.8]{6.8} we show that a sufficient condition for a model \zzz
\A\ of a given $\Pi^1_1$ theory \T\ \ul{not} to be a minimal model of \T\ is that the \zzz
notion of a well-founded linear order is not $\Sx^1_1/\A$, and this is also \zzz
a necessary condition for \A\ to be non-minimal in a certain strong \zzz
sense.

    We develop the ideas contained in the above proof into a theory of \zzz
indicators in, for example, non-\bx-models of $\Sx^1_1\HY\AC$ and non-well-founded \zzz
models of \KP\ for theories contained in countable fragments of \Lw. We \zzz
obtain a striking parallel between the \w- and non-\w-model cases. Typical \zzz
instances of our main theorem \hyperref[Chapter6:6.6]{6.6} for, say, the language of set theory \zzz
are as follows. If \T\ is an r.e.\ theory extending $\ZF^-$, then any non-\bbb-model \zzz
of \T\ has, for each $k\e\w$, a $k$-elementary transitive substructure \zzz
which is a model of \T. If \T\ is an r.e.\ theory extending \zzz
$\KP\up$ + $\Pi_1$-\Foundation, then any countable (or, more generally, locally countable) \zzz
nonstandard model of \T\ has a transitive substructure which is a \zzz
model of \T, and moreover, we can choose the substructure to have certain \zzz
saturation properties. Theorem \hyperref[Chapter6:6.6]{6.6} also provides model-theoretic proofs of \zzz
results in Feferman\,\cite{Fefe68}.

    Chapter~VI is concluded with an analogue of \hyperref[Chapter2:2.8]{2.8}. We show, for example, that \zzz
if \T\ is an r.e.\ theory in the language of analysis extending $\Sx^1_k\HY\AC$ with $k\ge 1$, then \zzz
\[
\myset{\tx\e\Sx^1_{k+1}}{\tx\text{ is }\Pi^1_{k+1}\text{-conservative over \T\ with the \w-rule}}
\]
and, if $k\ge 2$,
\[
\Myset{\tx\e\Pi^1_k}{\tx\text{ is }\Sx^1_k\text{-conservative over \T\ with the \w-rule}}
\]
are both $\Pi^0_1$ in Kleene's \sss O, and are complete for this class of sets. \zzz
We also prove a lemma concerning the semi-representation of $\Pi^1_1$ sets in r.e.\ \zzz
theories \T\ with the \w-rule: for example, if \T\ is not $\Pi^1_1$-sound, \zzz
then any $\Pi^1_1$ set may be semi-represented in \T\ with the \w-rule by \zzz
a formula in $\Dx^1_1\,(\T)$.

    The final chapter is an exposition of the Paris-Harrington result. \zzz
This chapter is presented as an introduction for the general reader, and \zzz
we give proofs of the truth and independence of the Paris-Harrington \zzz
statement, and we also show, by fairly simple model-theoretic means, the \zzz
sharp negative results established by the work of J.~Ketonen, L.~Kirby, \zzz
G.~Mills, J.~Paris, and R.~Solovay. We also prove an optimal result, \zzz
namely that for $c\eqx 2$,
\[
\PA\nprovex\Ax e\,\Ex n\:\big((e,n)\phs(e+1)^e_{\ol{c}}\,\big)
\]
(this has also been proved independently by others) and we give the equivalence \zzz
of the Paris-Harrington statement with the $\Sx_1$-Uniform Reflection \zzz
Principle of \PA. \\

\centerline{*\hskip 3em *\hskip 3em *}

    \refstepcounter{thisisdumb}\label{Introduction:BasicConventions}We shall
use the following notation and abuses of notation; areas of \zzz
logic where our abuses would lead to trouble are, while important, not of \zzz
interest to us here. We shall often not distinguish between formulae, etc., \zzz
and their codes, and it will often be implicitly assumed that some object \zzz
$x$ codes a finite sequence $\seq{(x)_0\m,(x)_1,\ldots,(x)_{\lenm x}}$.\footnote{Footnote added 2019: See
page \pageref{Chapter1:Definitionlenm} for the
notation \lenm{x}. In the 1980 thesis, this was denoted $\lambda x$.} $\vect x\ve y$ means that \zzz
$(x)_i\e y$ for all $i\le\lenm x$; $\vect x< y$, $\vect x\subseteq y$ are defined similarly.

    We shall often let \zzz
the same metavariable represent both formal variables and parameter \zzz
variables; for example, we might say ``if $\Ex x\m\tx x$ holds, choose a witness \zzz
$x$''. When we say, for example, that \Collection\ is the schema
\[
\Ax x\e a\,\Ex y\m\m\tx\imp\Ex b\,\Ax x\e a\,\Ex y\e b\m\m\tx\,,
\]
we mean that \Collection\ is the set of universal closure of sentences of this \zzz
form, where \tx\ ranges over all formulae of the language under consideration, \zzz
and where suitable precautions are taken to avoid collision of variables. \zzz
A \ul{theory} is a set of sentences.

    When there is no risk of confusion, the same symbol ``\e'' will have \zzz
three different uses: as a symbol of our formal language, as a binary relation \zzz
of some model under consideration, and with its usual informal meaning. \zzz
Likewise, we shall usually not need to distinguish between constant or function \zzz
symbols and their interpretation in some model. When we write $f\xvect x$ we \zzz
implicitly assume that the length of \xvect x is equal to the arity of $f$. If \zzz
$f$ is of arity \n, let $f''A = f[A^n]$.

    We shall be interested in models $\A=\seq{A,\e,\dots}$ of a weak set theory. \zzz
$\wf(\A)$ is the well-founded part of \A; we shall often implicitly suppose that \zzz
$\wf(\A)$ is equal to its transitive collapse. \A\ is \ul{nonstandard} if it \zzz
contains linear orderings which are internally well-founded but are not so \zzz
in the real world, and \A\ is an \ul{\w-model} if it does not contain non-standard \zzz
integers. We shall usually try to distinguish between ``contains'' \zzz
and ``includes''. A subset $B$ of $\wf(\A)$ is \ul{coded} in \A\ if there exists \zzz
$b\e A$ such that $B=\myset{a\e\wf(\A)}{\A\modelx a\e b}$. If there is no risk of confusion, \zzz
the structure in which sentences, especially quantifier-free ones, are to be \zzz
interpreted will only be implicitly mentioned; for example, we should usually \zzz
write the previous equation as $B=\myset{a\e\wf(\A)}{a\e b}$.

    \A\ is a model of \ul{overspill} ($\Sx_k$-\ul{overspill}) if for each ($\Sx_k$) formula \zzz
$\tx x$, possibly with parameters from \A, if $\tx a$ holds in \A\ for each \zzz
standard ordinal $a$ of \A, then $\tx a$ holds of some nonstandard ordinal of \zzz
\A, where here \ul{ordinal} refers to any element of any internally well-founded \zzz
linear ordering. A similar notion which is useful when dealing with weaker \zzz
theories is this: \A\ is \ul{recursively saturated} ($\Sx_{(k)}$-\ul{recursively saturated}) \zzz
if for any recursive set of ($\Sx_{(k)}$) formulae, with only a finite number \zzz
of free variables and perhaps a finite number of parameters from \A, if each \zzz
finite subset is realized in \A, then the whole set is realized in \A. \zzz
($\Sx_{(k)}$ is the closure of $\Pi_{k-1}$ under conjunction, disjunction, and existential \zzz
and bounded universal quantification; in the presence of $\Sx_k$-\Collection\ \zzz
and the pairing axiom, each $\Sx_{(k)}$ formula is equivalent to a $\Sx_k$ formula.) \zzz
A substructure \B\ of \A\ is an \ul{initial segment} (and \A\ is an \ul{end-extension} \zzz
of \B) if $x\e A$ and $y\e B$ with $x\e y$ implies $x\e B$; and a sub-structure \zzz
\B\ of \A\ is \ul{cofinal} in \A\ if for all $x\e A$ there exist \n\ \zzz
and $x_1,\dots,x_n\e A$, $y\e B$ such that $x\e x_1\e\dots\e x_n\e y$. \A\ is \zzz
\ul{locally countable} if $\myset{y\e A}{y\e x}$ is countable for all $x\e A$. \\

\centerline{*\hskip 3em *\hskip 3em *}

    I should like first of all to thank my supervisor, Professor Dana Scott, \zzz
and his substitute for the academic year 1978-9, Dr Robin Gandy;\footnote{Robin Gandy, 1919-1995.}\ and I am \zzz
grateful to the Association of Commonwealth Universities for enabling me \zzz
to study in Britain for three years under the Commonwealth Scholarship \zzz
scheme. I should also like to thank Professor Saul Kripke for informing \zzz
me of some of his results and for several conversations. I am most grateful \zzz
to Craig Smory\'nski and Alan Adamson\footnote{Alan Adamson, 1949-2012.}\ for many conversations and much \zzz
correspondence, and for trying to teach me a little logic. Alex Wilkie and Dan \zzz
Isaacson kindly lent me a great many preprints and were always ready to lend \zzz
an ear and offer suggestions; and I should like also to thank my fellow \zzz
students at the Mathematical Institute for patiently putting up with my \zzz
often foolish questions and comments. David Guaspari and Lawrence Kirby, \zzz
among others, generously sent me preprints of their work. I should also \zzz
like to thank all those who have kindly set aside an hour or so to talk \zzz
to me. I have left my greatest debt till last: because of the constant, \zzz
cheerful and uncomplaining support of my dearest wife Katherine, I have \zzz
dedicated this thesis to her.

\noindent
St Catherine's College, Oxford \\
1 April 1980

                                         \cleardoublepage
\thispagestyle{plain} \sectioncentred{The Main Definition}

    In this first chapter we define the notion of fulfilment, along with \zzz
several of its variants. After a lemma giving some simple properties of \zzz
this notion, we discuss its relationship to some of the early work in \zzz
Mathematical Logic by Skolem, Herbrand, and G\"odel. Then, using fulfilment, \zzz
we give results on the completeness and finite axiomatizability of \zzz
theories such as Peano's arithmetic---results whose proofs require only \zzz
\ul{very simple} model-theoretic arguments. Our independent statements, however, \zzz
suffer from the defect of being at best $\Pi^0_2$ whereas G\"odel's independent \zzz
sentences are $\Pi^0_1$. However, by making use of the \ul{fixed-point} theorem, \zzz
we shall generate interesting independent $\Pi^0_1$ (but self-referential) \zzz
statements of a model-theoretic nature. We shall close the chapter with \zzz
three miscellaneous results. The first is a model-theoretic proof (without \zzz
self-reference) of Tarski's Theorem: the theory of \bb N is not \zzz
arithmetical. The second shows how we may give a \ul{feasible} (without the \zzz
use of the fixed-point theorem or the use of such notions as pairs of \zzz
recursively inseparable sets) semi-representation of any r.e.\ set in \zzz
an r.e.\ theory which is not necessarily $\Sx^0_1$-sound. Finally, we give a \zzz
method of constructing models of various theories in the language of arithmetic \zzz
as quotients of certain ``constructively given'' function rings.

    Let \sss L be any finite first-order language. Let $(I,\prec)$ be any \zzz
linearly ordered set with a minimal element, which will be denoted by $0$, \zzz
and is such that each element $i$ of $I$ (except the maximal element, if \zzz
one exists) has an immediate successor, which will be denoted by $i+1$.\footnote{Footnote
added 2019: Here $0$ and $+$ are overloaded notations, and we should perhaps
write $0_I$ and $+_I$: for example, if $I$ is the set of positive odd integers,
then $0_I$ and $+_I$ are interpreted as stated, and the
first element of \sx\ would be $\sx_{0_{|\sx|}}$.  We will
also abuse the von Neumann ordinal notation: when we write
$\lenx{\sx}{3}$, we mean $\lenm\sx=\{0,1,2\}$ (as an ordered set), and so $|\sx|$ must
be $\{0,1,2,2\!+\!1\}$.}\, \zzz
\refstepcounter{thisisdumb}\label{Chapter1:Definitionlenm}Let $\sx = \seq{A_i}_{i\in I}$ be
a family of new unary predicate symbols, \zzz
let $\lenm\sx$ denote $I$ less (if it exists) the maximal element of $I$, and let $|\sx|\eqx I$. \zzz
\refstepcounter{thisisdumb}\label{Chapter1:DefinitionFulfil}For each
formula \px\ of \sss L and for each $i$ in $\lenm\sx$ define the (not \zzz
necessarily first-order) formula $\px^\sx_i$ inductively as follows. (Here $i$ and \zzz
$j$ are understood to range over $\lenm\sx$. When the context is clear, the \zzz
superscript \sx\ will be omitted.)
\begin{mylist}
\item If \tx\ is atomic, let $\tx_i\eqx \tx$ and $(\nx\tx)_i=\nx\tx$\,.
\item $(\Ex x\m\tx)_i=\Ex x\e A_{i+1}\m\tx_i$\,.
\item $(\Ax x\m\tx)_i=\Ax j\gex i\,\Ax x\e A_j\m\tx_j$, where $j$ is a new variable.
\item $(\nx\Ax x\m\tx)_i=(\Ex x\,\nx\tx)_i$\,; \\
      $(\nx\Ex x\m\tx)_i=(\Ax x\,\nx\tx)_i$\,.
\item Otherwise, $\tx_i$ is obtained from \tx\ by replacing each positive \zzz
instance of a Boolean atom \yx\ of \tx\ by $\yx_i$ and each negative \zzz
instance by $\nx(\nx\yx)_i$\,.
\end{mylist}
For example, if \tx\ is quantifier-free, then $(\Ax x\,\Ex y\,\Ax u\,\Ex v\m\tx)^\sx_0$ is just
\[
\Ax i\gex 0\,\Ax x\e A_i\,\Ex y\e A_{i+i}\,\Ax j\gex i\,\Ax u\e A_j \,\Ex v\e A_{j+1}\m\tx\,;
\]
and $(\px\imp\yx)^\sx_i$ is $\nx(\nx\px)^\sx_i\imp\yx^\sx_i$. \sx\ is \ul{increasing} if for all $i$ \zzz
in \lenm\sx, $A_i\subseteq A_{i+1}$, and \sx\ is \ul{closed} if for each $i$ in \lenm\sx\ and for each \zzz
function $f$ of \sss L, $f''A_i\subseteq A_{i+1}$, \ul{and} if each constant of \sss L is contained \zzz
in $A_0$. Let $\A=\seq{A,\dots}$ be a structure for \sss L and for each $i$ in \zzz
$I$, let $A_i$ be interpreted as a non-empty subset of $A$. Let \px\ be any \zzz
formula of \sss L and let \vect a\ be a valuation in $A$ of the free variables of \zzz
\px. Say \sx\ \ul{fulfils} \px\xvect a, and write $\A\modelx(\px\xvect a\m)^\sx$, if:
\[
\A\modelx(\px\xvect a\,)^\sx_0 \text{ \ul{and} } \sx \text{ is increasing and closed.}
\]
(The last two requirements ensure that \sx\ fulfils $\Ax x\,\Ex y\,(x=y)$ and \zzz
$\Av x\,\Ex y\,(f\m\m\vect x=y)$ for every function symbol $f$.) We shall occasionally use the phrase \zzz
\lenm\sx-fulfil when we wish to indicate \lenm\sx.  It will be convenient \zzz
to allow a closed sequence of length one $\seq{A_0}$ to vacuously fulfil any statement.

    The motivation behind this definition should be clear. For example, \zzz
suppose \px\ is a sentence which is true in a structure \A. Then there \zzz
exist many sequences of subsets of $A$ which fulfil \px: for example, the \zzz
trivial sequence $\seq{A,A,A,\dots}$. A more useful construction---one which is \zzz
due to Skolem---is as follows. Choose (using the axiom of choice) a \zzz
set $F$ of \ul{satisfaction functions} for \px: for example, if
\[
\px=\Ax x\,\Ex y\,\Ax u\,\Ex v\m\tx xyuv
\]
where \tx\ is quantifier-free, choose $f$ and $g$ such that
\[
\A\modelx\Ax x\,\Ax u\,\tx(x,f\mx x,u,gxu)
\]
and set $F=\{\,f, g\,\}$. Let $A_0\ne\emptyset$ contain the (interpretations of \zzz
the) constants of \sss L and let
\[
A_{i+1}=\cupx_{f\,\e\,\sss L\,\cup\,F}\, f''A_i
\]
for each $i$ in \w. Then $\seq{A_i}_{i<\w}$ fulfils \px, as does every finite \zzz
subsequence.

    \refstepcounter{thisisdumb}\label{Chapter1:DefinitionFulfilIII}In Chapter~III it
will be convenient to use a slight variant of the \zzz
above definition: let $i$ and $j$ range over $|\sx|$ rather than \lenm\sx, and change \zzz
clause (ii) above to read
\[
(\Ex x\m\tx)_i=
\begin{cases}
\Ex x\e A_{i+1}\m\tx_{i+1}\,,&\text{if }i\e\lenm\sx \\
\text{true,}&\text{otherwise.}
\end{cases}
\]
The motivation behind this definition is also simple. In the above \zzz
example, choose $F$ to be a set of so-called \ul{Skolem functions} for the \zzz
existential quantifiers of \px; that is, if \px\ is as before, choose $f$ \zzz
and $g$ such that
\vspace{-1ex}
\begin{flalign*}
           && &\A\modelx\Ax x\,\Ax u\,\Ex v\m\tx\m(x,f\mx x,u,v)&\\
\text{and} && &\A\modelx\Ax x\,\Ax y\,\Ax u\,\big(\Ex v\m\tx xyuv\imp\tx\m(x,y,u,gxyu)\big)\,.
\end{flalign*}
Now if we define $A_0\m,A_1,A_2,\dots$, as before, but using these Skolem \zzz
functions rather than the satisfaction functions, we obtain a sequence \zzz
which fulfils \px\ according to this second notion of fulfilment.

    It is clear that if \sx\ fulfils a sentence \px\ according to the \zzz
first definition, then it also fulfils \px\ according to the second. \zzz
Furthermore, a moment's thought shows that if $\sx=\seq{\sx_i}_{i\lex\lenm\sx}$ fulfils \zzz
\px\ according to the second definition, where $\lenm\sx\lex\w$, and if \px\ has $k$ \zzz
existential quantifiers, then the sequence
\[
\seq{\sx_{ik}}_{ik\lex\lenm\sx}
\]
fulfils \px\ according to the first. Thus for our purposes the two \zzz
notions are equivalent, and we shall always use the first except in \zzz
Chapters III and VII.

    There are three more variants of the notion of fulfilment which \zzz
naturally occur. The definitions are listed here below for ease of \zzz
reference, even though it would perhaps be less daunting to the reader, \zzz
and clearer, if they were deferred until they were actually needed.

    For the first variant, we do not start with a given fixed structure, \zzz
but rather we consider pairs $(\sx,\A)$ where $\seq{\sx_i}_{i\in |\sx|}$ is an increasing \zzz
sequence of sets of \ul{new constant symbols}, and \A\ is a structure for \sss L \zzz
with domain $\cupx_{i\in\m|\sx|}\,\sx_i$, and where \sx\ is closed under the functions of \sss L \zzz
as interpreted in \A. (We do not need the function symbols of \sss L to be \zzz
interpreted as \ul{total} functions; it suffices that they be defined wherever \zzz
needed.) In our notation we shall not mention the structure \A\ \zzz
explicitly. Say \sx\ \ul{$\ast$fulfils} a \ul{sentence} of \sss L, and write $\px^{*\sx}$, if
\[
\A\modelx\px^\sx\,.
\]
The \ul{raison d'\^etre} of the notion of $\ast$fulfilment is this. If \px\ is a \zzz
sentence of any finite language \sss L, then the sentence
\[
\text{for all $n\e\w$, there exists \sx\ with \lenn\ which $\ast$fulfils \px}
\]
may, by the usual coding techniques, be expressed by a $\Pi^0_1$ sentence of \zzz
arithmetic (for we can easily estimate an upper bound on the cardinality of $\sx_n$, \zzz
that is, of the structure \A; the bound is in fact an $\sss E_3$ function \zzz
of $n$ and the number of quantifiers in \px); whereas the sentence
\[
\text{for all $n\e\w$, there exists \sx\ with \lenn\ which fulfils \px}
\]
may not be expressed in arithmetic at all except when \px\ itself is in \zzz
the language of arithmetic, and even then it must in general be expressed \zzz
by a $\Pi^0_2$ sentence.

    For our next definition we consider a language \sss L which contains \zzz
(or, at least, some definitional extension of it contains) a type (or \zzz
a unary predicate symbol) $\ol{\w}$, constant symbols $\ol0$ and $\ol1$, and function \zzz
symbols $\ol{+\vphantom\times}$ and $\ol\times$ whose domain is $\ol\w$. A sequence \sx\ (or more \zzz
precisely, an ordered pair $(\sx,\A)$ \ul{half-$\ast$fulfils} a sentence \px, and \zzz
we write $\px^{\frac12\ast\sx}\!$, if: $\sx=\seq{\sx_i}_{i\in |\sx|}$ is an increasing sequence of \sss L sets \zzz
which contain \ul{both} constant symbols and integers; \A\ is a structure \zzz
with domain $\cupx_{i\in\m|\sx|}\,\sx_i$ in which the interpretation of $\ol{\w}$ consists of integers, \zzz
and in which the interpretation of the symbols
$\ol{0}$, $\ol{1}$, $\ol{+\vphantom\times}$, and $\ol{\times}$ is the \zzz
standard one; and
\[
\A\modelx\px^\sx\,.
\]
Thus if \px\ is a sentence in the language of, say, set theory or analysis, \zzz
we may express
\[
\text{for all $\n\e\w$, there exists \sx\ with \lenn\ which $\textstyle\frac12\ast$fulfils \px}
\]
by a $\Pi^0_2$ sentence of arithmetic. Note that we may also allow \px\ to \zzz
have free variables of integer type.

    \refstepcounter{thisisdumb}\label{Chapter1:Definitionifulfilment}In our
final variation of the notion of fulfilment, we consider \zzz
structures $\A=\seq{A,\dots}$ for the language of either set theory or \zzz
arithmetic. The notion of \ul{i-fulfilment} (for \ul{i}nitial) is like that of fulfilment, \zzz
except we require of our sequences $\sx=\seq{\sx_i}_{i\e|\sx|}$, in the case of \zzz
arithmetic, that each $\sx_i$ be an \ul{initial segment} of $A$ (i.e.\ $x\ltx y\e\sx_i$ \zzz
implies $x\e\sx_i$), and in the case of set theory, that for all $i$ in \lenm\sx, \zzz
$x\e y\e\sx_i$ implies $x\e\sx_{i+1}$. In the case of arithmetic, (the internalized \zzz
version of) i-fulfilment gives a particularly elegant formulation of fulfilment, \zzz
for rather than considering sequences of \ul{sets} of integers we can \zzz
just consider sequences of \ul{integers}, and say that $\sx=\seq{\sx_i}_{i\le n}$ i-fulfils \zzz
e.g.\ $\Ax x\,\Ex y\m\tx$ iff
\[
\Ax i\ltx n\,\Ax x\ltx\sx_i\,\Ex y\ltx\sx_{i+1}\m\tx\,.
\]
We shall use this observation to simplify the statements or proofs of \zzz
two or three results. One could combine the notions of i-fulfilment and \zzz
$\frac12\ast$fulfilment, but we shall have little need for this.

    \refstepcounter{thisisdumb}\label{Chapter1:Lemma11}Next we
shall consider a few simple properties of the notion of \zzz
fulfilment. The first two parts of the following lemma are due to Kripke.

\thm{1.1 Lemma} Let $\A=\seq{A,\dots}$ be any structure for \sss L, let \px\ be any \zzz
formula of \sss L (which may contain parameters from $A$) and let $\sx=\seq{A_i}_{i\in |\sx|}$ \zzz
be a sequence of subsets of $A$.
\begin{mylist}
\item If \lenm\sx\ is unbounded and if \sx\ fulfils \px, then \px\ is true in
\[
\B=\A\,\up\cupx_{i\in\m\lenm\sx}A_i
\]
\item Let $J\subseteq|\sx|$ have a minimal element and be such that each element of \zzz
$J$ (except the maximal element, if one exists) has an immediate successor \zzz
in $J$, and let $\sx\up J=\seq{A_i}_{i\in J}$. If \sx\ fulfils \px, then $\sx\up J$ \zzz
fulfils \px.
\item If $|\sx|\eqx \w$, if $\seq{A_i}_{i\le n}$ fulfils \px\ for each $n\e\w$, and if each \zzz
$A_i$ is \ul{finite}, then $\seq{A_i}_{i<\w}$ fulfils \px.
\end{mylist}
Furthermore, i, ii, and iii also hold for $\ast$, $\frac12\ast$, and i-fulfilment.

\noindent
Proof: (i)\, Let \yx\ be obtained from \px\ by first eliminating all \zzz
occurrences of ``\imp'' by use of ``\nx'' and ``\vx'', and then ``pushing'' all \zzz
negations inside. From the definition of fulfilment, $\yx^\sx$ iff $\px^\sx$. \zzz
Suppose $\px^\sx$ and so $\yx^\sx$. Then \Xx\ holds in \B, where \Xx\ is obtained \zzz
from $\yx^\sx$ by removing the bounds on the existential quantifiers. Working \zzz
from the innermost quantifiers of \Xx\ outwards, replace each quantifier \zzz
pair of the form ``$\m\Ax j\gex i\,\Ax x\e A_j$'' by ``$\m\Ax x$'': each successive \zzz
alteration maintains the truth in \B. Hence \yx, and \px, hold in \B.

\noindent
(ii)\, Again we may suppose that \px\ is in negation normal form. The \zzz
result is now obvious.

\noindent
(iii)\, The proof is by induction on the length of \px, which we can assume \zzz
is in negation normal form. There is nothing to prove if \px\ is quantifier-free, \zzz
and the conjunction and disjunction steps are trivial. Let $\sx\up\,(n+1)=\seq{A_i}_{i\lex n}$. For \zzz
the induction to work, we need to show a bit more: for all $k\geX 0$ and all \tx\ with \zzz
any evaluation of its free variables, $\tx^\sx_k$ iff $\Ax n\gtx k\,\tx^{\smash{\sx\upharpoonright (n+1)}}_k$.

\noindent
For the universal quantifier we have
\vspace{-.4ex}
\begin{flalign*}
&&                &\Ax n\gtx k\,(\Ax x\,\px)^{\smash{\sx\upharpoonright (n+1)}}_k &&&& \\
&&\text{iff\quad} &\Ax n\gtx k\,\Ax i\gex k\,\Ax x\e A_i\,\px^{\smash{\sx\upharpoonright (n+1)}}_i \text{, by definition, with $i$ new, $i\ltx n$}\\
&&\text{iff\quad} &\Ax i\gex k\,\Ax x\e A_i\,\Ax n\gtx i\,\px^{\smash{\sx\upharpoonright (n+1)}}_i    \\
&&\text{iff\quad} &\Ax i\gex k\,\Ax x\e A_i\,\px^\sx_i \text{, by our induction hypothesis}  \\
&&\text{iff\quad} &(\Ax x\,\px)^\sx_k, \text{ by definition.} \\
&\text{For the existential quantifier,}\hidewidth \\
&&                &\Ax n\gtx k\,(\Ex x\,\px)^{\smash{\sx\upharpoonright (n+1)}}_k \\
&&\text{iff\quad} &\Ax n\gtx k\,\,\Ex x\e A_{k+1}\,\px^{\smash{\sx\upharpoonright (n+1)}}_k \text{, by definition} \\
&\text{and because $A_{k+1}$ is finite,\footnotemark\ \footnotemark\ by the pigeon-hole principle this holds}\hidewidth  \\
&&\text{iff\quad} &\Ex x\e A_{k+1} \text{ such that for infinitely many $n\gtX k$, } \px^{\sx\upharpoonright(n+1)}_k \\
&&\text{iff\quad} &\Ex x\e A_{k+1}\,\Ax n\gtx k\,\px^{\smash{\sx\upharpoonright (n+1)}}_k \text{, by (ii)}\\
&&\text{iff\quad} &\Ex x\e A_{k+1}\,\px^\sx_k \text{, by the induction hypothesis}              \\
&&\text{iff\quad} &(\Ex x\,\px)^\sx_k, \text{ by definition.}\tag*{\done}
\end{flalign*}
\vspace{-.2ex}
\addtocounter{footnote}{-1}
\footnotetext{Footnote added 2019: The finiteness is
generally needed. For example, consider the sentence ``the odd numbers are bounded''
in the language \sss L = \{even,\, <\} given by $\Ex x\,\Ax y\,\big(x\ltx y\imp \rrm{even}(y)\big)$. Let $A_0$ be the
even natural numbers.  Add the first odd number to $A_1$, the second to $A_2$, etc. Each finite subsequence fulfils the
sentence, but the whole does not.}
\stepcounter{footnote}
\footnotetext{Footnote added 2019: But for many applications the finiteness restriction can be bypassed.
If \sss L has coding abilities, and \A\ is recursively saturated or is a model of $\Sx_1$-overspill,
and if \px\ is \n-fulfillable for all \n, then \px\ is \w-fulfillable.}
\indent
   Now that the above machinery has been set up, we can give a succinct \zzz
account of some of the earliest work in Mathematical Logic. As already \zzz
mentioned, Skolem in \cite{Skol20}, in his proof of L\"owenheim's Theorem, observed \zzz
that if a sentence \tx\ is true in a structure \A, then by the axiom of \zzz
choice we can find a sequence $\seq{\sx_i}_{i\in\m\w}$ of finite subsets of $A$ which fulfils \zzz
\tx; $\A\,\up\cupx_i\sx_i$ will then be a countable substructure of \A\ which is a model \zzz
of \tx.

    In \cite{Skol22}, \cite{Skol28}, and \cite{Skol29}, Skolem saw that this shows
\[
\text{if \tx\ is satisfiable, then } \Ax n\,\Ex\sx\,(\lenn\ax\tx^{\ast\sx})\,, \label{Chapter1:eq:1}\tag{1}
\]
and, moreover, the converse is also true, by a simple argument using \zzz
K\"onig's Infinity Lemma. Skolem pointed out (in \cite{Skol28}) that this gives a \zzz
proof procedure which, as we see, is cut-free and has the subformula \zzz
property. (Aside: Kripke has pointed out an elegant, nonstandard argument \zzz
for the converse of \eqref{Chapter1:eq:1}. Suppose $\Ax n\,\Ex\sx\,(\lenn\ax\tx^{\ast\sx})$. Then in any \zzz
proper elementary extension of \bb N we can find (the code of) a \zzz
sequence \sx\ of \ul{nonstandard} length which $\ast$fulfils \tx. Then by 1.1.i and \zzz
ii, \tx\ is true in the structure determined by \sx\ with domain $\cupx_{\,n\text{ standard} }\,\sx_n$.)

    In \cite{Herb30}, Herbrand gave an (incomplete) proof of
\[
\provex\nx\tx\quad\text{implies}\quad\nx\Ax n\,\Ex\sx\,(\lenn\ax\tx^{\ast\sx}) \label{Chapter1:eq:2}\tag{2}
\]
and its converse; moreover, his proof is \ul{effective} in that we may obtain \zzz
primitive recursively a witness $n$ from any proof of \nx\tx, and conversely, \zzz
a (bound on the) proof from any witness $n$.

    \refstepcounter{thisisdumb}\label{Chapter1:GodelsProof}From the converses of \eqref{Chapter1:eq:1} and \eqref{Chapter1:eq:2} we may immediately obtain, as many \zzz
people have remarked, the Completeness Theorem of G\"odel\,\cite{Godel30}. G\"odel's \zzz
proof was, very roughly, as follows. He proved the converse of \eqref{Chapter1:eq:1} as \zzz
Skolem did, so let us consider the converse of \eqref{Chapter1:eq:2}. First note that for \zzz
each formula \tx\ of \sss L and each integer $n$, we may express the notion
\[
\Ex\sx\,(\lenn\ax\tx^\sx) \label{Chapter1:eq:3}\tag{3}
\]
in the language of \sss L by writing out the elements of the terms in \zzz
\ul{explicitly}. For example, if $\tx=\Ax x\,\Ex y\m\m\yx\mx xy$, with \yx\ quantifier-free, and \zzz
if the language \sss L has no function symbols or constants, we may write \zzz
\eqref{Chapter1:eq:3} as
\[
\Ex x_0\m,x_1,\dots,x_n\,\aax_{\,i\ltx n}\yx\,x_i x_{i+1},
\]
With this understanding, we can show that, in any axiomatization of the \zzz
predicate calculus which we may wish to consider, for each $n\e\w$,
\[
\provex\tx\imp\Ex\sx\,(\lenx{\sx}{\ol n}\ax\tx^\sx).\label{Chapter1:eq:4}\tag{4}
\]
So now to show the converse of \eqref{Chapter1:eq:2}, suppose for some $n$, $\nx\Ex\sx\,(\lenn\ax\tx^{\ast\sx})$. \zzz
Then ${\modelx\nx\Ex\sx\,(\lenn\ax\tx^\sx)}$. But this is essentially a propositional formula, and \zzz
as the propositional calculus is complete, we have $\provex\nx\Ex\sx\,(\lenn\ax\tx^\sx)$. \zzz
Hence by \eqref{Chapter1:eq:4}, $\provex\nx\tx$. Note that this proof is also effective.

    This completes our discussion of early logic. We shall next consider \zzz
various incompleteness results.

   To obtain these results, we shall need to be able to formalize the \zzz
notion of fulfilment. We wish our results to be applicable to many different\refstepcounter{thisisdumb}\label{Chapter1:DefinitionMembership} \zzz
areas---arithmetic, analysis, and set theory---and so to state them in \zzz
the necessary generality, we shall usually consider some fixed but \zzz
arbitrary r.e.\ theory \s\ of \sss L in which a notion of \ul{membership} \e\ \zzz
may be defined so that \s\ proves
\begin{flalign*}
&\Pair\!:  &&\Ax x,y\,\Ex z\,\big(x\e z\ax y\e x\ax\Ax u\e z\,(u=x\vx u=y)\m\big)\,, \text{ and} &&&\\
&\Union\!: &&\Ax x\,\Ex y\,(\Ax z\e x\,\Ax u\e z\,u\e y\ax\Ax u\e y\,\Ex z\e x\,u\e z)\,.
\end{flalign*}
(For example, in arithmetic, say $x\e y$ if
\[
\Ex s,t\,(t\ltx 2^x\ax y=s2^{x+1}+2^x+t)\,,
\]
and in analysis, extend the usual notion of membership $\e\subset\w\times\sss P\w$ by the \zzz
class
\[
\Myset{(X,Y)}{\Ex n\,\Ax m\,(m\e X\equiv\seq{n,m}\e Y)}
\]
where $\lambda nm.\seq{m,n}$ is a standard pairing function on the integers.) \zzz
Furthermore, we shall assume that some notion of \ul{natural number} may be \zzz
defined in \s, and that \s\ proves $\Dx_1(\s)$-\Induction, where as usual we have
\vspace{-.5ex}
\begin{align*}
\Dx_0=\Pi_0=\Sx_0 &\,=\text{ closure of atomic formulae under } \ax,\vx,\nx,\,\Ax x\e y,\,\Ex x\e y\\
\Sx_{k+1} &\,=\text{ closure of $\Pi_k$ under $\Ex x$ and equivalences in the \PC}\\
\Pi_{k+1} &\,=\text{ dual of $\Sx_{k+1}$ }\,\\
\Dx_{k+1}(\T) &\,=\text{ formulae provably equivalent in \T\ to a $\Sx_{k+1}$ and a $\Pi_{k+1}$ formula}\,.
\end{align*}
Let $\ol0$ denote the first natural number as represented in \s, $\ol1$ the \zzz
second, etc. Finally, we shall suppose that \s\ proves that some $\Sx_1$ \zzz
formula is a satisfaction predicate for $\Dx_0$ formulae. This latter \zzz
assumption may often be omitted, but it makes our exposition somewhat \zzz
easier. In particular, it enables us to treat various restricted schemata \zzz
as \ul{single} sentences, which saves us the trouble of picking the appropriate \zzz
instances. Also, for certain proof-theoretic results we shall need that \s\ proves \zzz
the fourth Grzegorczyk function $\sss E_4$ is total.

    \refstepcounter{thisisdumb}\label{Chapter1:DefinitionInduction}Among the schemata which we will find useful to consider are:
\begin{flalign*}
&\Induction\!:  &&\tx\ol0\ax\Ax n\,(\tx n\imp\tx (n+1))\imp\Ax n\m\tx n\,,&&&&&\\
&\Foundation\!: &&\Ex x\m\tx x\imp\Ex x\,(\tx x\ax\Ax y\e x\,\nx\tx y)\,,\\
&\Separation\!: &&\Ex b\,\Ax x\e a\,(x\e b\,\equiv\m\tx)\,,\\
&\Collection\!: &&\Av x\ve a\,\Ev y\,\tx\imp\Ex b\,\Av x\ve a\,\Ev y\ve b\,\tx\,,\text{ and}\\
&\Bounding\!:   &&\Ex b\,\Ax x\e a\,(\Ex y\m\tx\imp\Ex y\e b\,\tx)
\end{flalign*}
where $b$ does not occur free in \tx. $\Lambda$-\Induction, etc., is the schema \zzz
of \Induction\ with \tx\ restricted to $\Lambda$. \Infinity\ is the axiom asserting \zzz
that the integers form a set. The variables $n$ and $m$ in will always be \zzz
restricted to the integers.

   \refstepcounter{thisisdumb}\label{Chapter1:DefinitionPAex}
   One base theory we shall occasionally consider is $\PA^-_\text{ex}$, that is, arithmetic with \zzz
axioms for addition, multiplication, exponentiation and membership, \zzz
and with induction limited to $\Dx_0$ formulae.

   The relation
\[
(\px\vect a)^\sx_i
\]
may be expressed as a $\Dx_1(\s)$ predicate of \sx, $i$, (the code of the sequence) \zzz
$\vect a$, and (the G\"odel number of) \px, and henceforth this notation will \zzz
implicitly refer to some such suitable representation. If \sx\ satisfies the extra \zzz
requirement for i-fulfilment, we write i(\sx). For a formula $\T\m(x)$ \zzz
(which, as the notation suggests, is to be considered as representing \zzz
a class of sentences) let $(\T)^\sx$ denote
{
\setlength{\abovedisplayskip}{1.5ex}
\setlength{\belowdisplayskip}{1ex}
\setlength{\abovedisplayshortskip}{.5ex}
\setlength{\belowdisplayshortskip}{1.0ex}
\[
\Ax\,\godel\tx\e\lenm\sx\,\big(\T\m(\godel\tx)\imp\tx^\sx\big)\,.
\]
\refstepcounter{thisisdumb}\label{Chapter1:DefinitionTr}For a formula $\tx\m\vect x$, let $(\Tr\tx\m\vect x)^\sx$ denote
\[
\Ax i\e\lenm\sx\,\Av x\ve\sx_i\,\big(\tx\m\vect x\imp(\tx\m\vect x)^\sx_i\big)\,,
\]
and let $(\Tr\Sx_k)^\sx$ denote
\[
\Ax\m i,\godel\tx\e\lenm\sx\,\Ax x\e\sx_i\,\big(\Sat_{\Sx_k}(\godel\tx,x)\imp(\tx\m x)^\sx_i\big)
\]
}where $\Sat_{\Sx_k}(e, x)$ is a standard satisfaction predicate for formulae \zzz
with one free variable. Parentheses will usually be omitted when there \zzz
is no risk of ambiguity.\footnote{Footnote added 2019: To clarify, for a formula $\T(x)$ in the language
of \s, we denote $\T=\{x:\T(x)\}$, and let
{
\setlength{\abovedisplayskip}{1ex}
\setlength{\belowdisplayskip}{.1ex}
\setlength{\abovedisplayshortskip}{.2ex}
\setlength{\belowdisplayshortskip}{1.0ex}
\[
\T\m\m^\sx\,\mathop:\!=\;\Ax\tx\e\T\cap\lenm\sx\;\;\tx\m^\sx
\]
\[
\Tr\T\m\m^\sx\,\mathop:\!=\;\Ax i\e\lenm\sx\;\;\Ax\tx\e\T\cap\sx_i\;\;(\tx\e\T)^\sx_i\,.
\]
The notation ``$\Tr$'' in the second phrase is perhaps unfortunate: it has little to do with truth,
and merely ensures that the notion of being an axiom of \T\ persists in the submodels we create.
}
}

    \refstepcounter{thisisdumb}\label{Chapter1:Lemma12}By formalizing Skolem's argument above, we obtain:

\thm{1.2 Lemma} For each formula $\px\m\vect x$ and each integer $k$,
\[
\Av x\,\big(\px\m\vect x\imp\Ax n\,\Ex\sx\,(\lenn\ax\px\m\vect x\m^\sx\ax\Tr\Sx_k\m^\sx)\m\big)
\]
is provable in $\s+\Induction$.\done

\noindent
We omit the proof and just remark that although \px\ may not have definable \zzz
satisfaction functions, using \Induction\ we can still define some \zzz
suitable finite approximations which will suffice to carry out the \zzz
construction. If the schema of \Collection\ is also assumed, we may \zzz
require that i\sx\ hold. In the presence of a strong theory such as \ZF, \zzz
we also have a simple proof Lemma 1.2 via the set-theoretical reflection \zzz
principle---although the usual proofs of this principle (see, for example, \zzz
Krivine\,\cite{Kriv71}) make implicit use of the notion of fulfilment. Thus for \zzz
the study of theories as strong as \ZF, the notion of fulfilment is \zzz
unnecessary; indeed, much of the value of fulfilment stems from the fact \zzz
that by using it, we can mimic arguments of \ZF\ in much weaker theories. \zzz
If the negation of \Infinity\ is assumed, the schema of \Collection\ is \zzz
implied in \s\ by that of \Induction\ (see remarks in Chapter~4 on page \pageref{Chapter4:ArithmeticCollectionComment}); and so \zzz
in particular we have
\[
\PA\provex\px\imp\Ax n\,\Ex\sx\,(\lenn\ax\text{i}\sx\ax\px^\sx\ax\Tr\Sx_k\m^\sx)\,.
\]
   \vspace{-.5ex}
   \refstepcounter{thisisdumb}\label{Chapter1:1.3}
   \indent
   The proof of the next theorem is the template for many later results. \zzz
Let $\ful_\px(n,\sx)$ be either
\begin{flalign*}
&          &&\lenn\ax n\e\sx_0\ax\text{i}\sx\ax\px^\sx &\\
&\text{or} &&\lenn\ax n\e\sx_0\ax\Tr\Dx_0\m^\sx\ax\px^\sx\,.
\end{flalign*}
\thm{1.3 Theorem}(Kripke) If \px\ is a sentence consistent with $\s+\Dx_0$-\Foundation,
\[
\px\imp\Ax n\,\Ex\sx\,\ful_\px(n,\sx) \label{Chapter1:eq:5}\tag{5}
\]
is not a theorem of the predicate calculus.

\noindent
Proof: Let $\A=\seq{A,\dots}$ be any non-\w-model of \s\ + $\Dx_0$-\Foundation, \px, and \eqref{Chapter1:eq:5}. Pick \zzz
any nonstandard \n\ in $A$ and by $\Dx_0$-\Foundation\ choose $\tau$ in $A$ so \zzz
that $\ful_\px(n,\tau)$ holds but for no $\sx\e\tau_n$, does $\ful_\px(n,\sx)$ hold. \zzz
By Lemma 1.1, $\B=\A\,\up\cupx_{i\in\m\w}\tau_i$ is a model for \px, and as the formula $\ful_\px$ \zzz
is absolute between \A\ and \B, for no element \sx\ of \B\ does $\ful_\px(n,\sx)$ \zzz
hold.\done

\refstepcounter{thisisdumb}\label{Chapter1:1.4}
\thm{1.4 Corollary}(Ryll-Nardzewski, Mostowski; this proof due to Kripke) \zzz
No consistent extension of \s\ + $\Dx_0$-\Foundation\ + \Induction\ is finitely axiomatizable.

\noindent
Proof: Suppose \px\ is an axiom for such a theory. By Lemma 1.2
\begin{flalign*}
&                      & \px\,&\provex\px\imp\Ax n\,\Ex\sx\,\ful_\px &&\\
&\text{and so}\hidewidth &    &\provex\px\imp\Ax n\,\Ex\sx\,\ful_\px\,.
\end{flalign*}
This contradicts the theorem.\done

\refstepcounter{thisisdumb}\label{Chapter1:1.5}
\thm{1.5 Corollary}(Kripke) Let \T\ be any consistent theory extending \s\ + $\Dx_0$-\Foundation\ \zzz
which is semi-representable in \T. Then
\[
\Ax n\,\Ex\sx\,(\lenn\ax n\e\sx_0\ax\Tr\Dx_0\m^\sx\ax\Tr\T\m^\sx\ax\T\m^\sx)
\]
is not provable in \T, where ``\T'' in this sentence is understood to \zzz
be any formula which semi-represents \T\ in \T. In particular, there is \zzz
a true $\Pi^0_2$ sentence not provable in \PA, and a true $\Sx_1$ sentence not \zzz
provable in \ZFC\ (assuming \ZFC\ has a standard model).\done

\noindent
This corollary may either be proved by the same method or, if ``\T'' is $\Pi_1$, be \zzz
derived from 1.3; we omit the proof. (See, however, 1.10 below.)

   Corollary 1.5 is of course weaker the G\"odel's First Incompleteness \zzz
Theorem\,\cite{Godel31} which gives (true) $\Pi^0_1$ sentences which are not provable \zzz
in \PA\ and \ZF; furthermore, 1.5 does not apply to theories without the \zzz
schema of $\Dx_0$-\Foundation. But we can improve our results as follows. \zzz
By use of the notion of $\frac12\ast$fulfilment, the proofs of 1.3 and 1.5 yield \zzz
$\Pi^0_2$ sentences which independent of, say $\ssf A_2$ and \ZF, respectively, \zzz
as stated in 1.6 below, and we also have that $\ssf A_2$ is not finitely \zzz
axiomatizable. By making use of the fixed-point theorem, however, we \zzz
may obtain independent $\Pi^0_1$ sentences as in 1.8 below.

\thm{1.6 Theorem} If \px\ is a sentence consistent with \s\ then
\[
\px\m\imp\m\Ax n\,\Ex\sx\,(\lenn\ax n\e\sx_0\ax\Tr\Dx^0_0\m^\sx\ax\px^{\frac12\ast\sx})
\]
is not a theorem of the predicate calculus, and so no consistent \zzz
extension of \s\ + \Induction\ is finitely axiomatizable. (Here it is \zzz
immaterial whether we write $\Tr\Dx^0_0\m^\sx$ or $\Tr\Dx^0_0\m^{\frac12\ast\sx}$.) And if \T\ is a \zzz
consistent theory extending \s\ which is semi-representable in (some \zzz
consistent extension of) \T, then
\[
\Ax n\,\Ex\sx\,(\lenn\ax n\e\sx_0\ax\Tr\Dx^0_0\m^\sx\ax\Tr\T\m^{\frac12\ast\sx}\ax\T^{\frac12\ast\sx}
\]
is not provable in \T, where ``\T'' in this sentence is understood to be \zzz
the representation of \T.\done

\noindent
The proof of 1.6 is almost word-for-word the same as those of 1.3 and 1.4 \zzz
and is left to the reader.

\refstepcounter{thisisdumb}\label{Chapter1:1.7}
    We have exhibited sentences which are not provable in various theories \zzz
\T, but without some extra assumptions on \T\ their negations may very \zzz
well be provable. Let us consider, for example, an r.e.\ theory \T\ in \zzz
the language of arithmetic. If \T\ is $\Sx^0_1$-sound, the sentence of 1.5 \zzz
is a true $\Pi^0_2$ not provable in \T, and so if \T\ is also $\Sx^0_2$-sound, its \zzz
negation is not provable. By use of the fixed-point theorem, however, we \zzz
may obtain independent sentences without this soundness hypothesis.

\thm{1.7 Theorem}(Rosser; this proof, Kripke) Suppose \T\ extends \s\ \zzz
and that \T\ is binumerated in (some consistent extension of) \T. \zzz
Then \T\ is not complete.

\noindent
Proof: Choose \tx\ so that \T\ proves:
\[
\tx\,\equiv\,\Ax n\,\big(\Ex\sx\,(\lenn\ax\T^\sx\ax\tx^\sx)\imp\Ex\sx\,(\lenn\ax\T^\sx\ax(\nx\tx)^\sx)\big)\,,
\]
where ``\T'' in this sentence is to be understood as above. We claim that \zzz
\tx\ is neither provable nor refutable in \T. Let \A\ be any model of \zzz
\T\ (plus the theory $\myset{\T\m(\godel\tx)}{\tx\e\T}$).

    Suppose \tx\ holds in \A. We can assume that overspill holds in \zzz
\A, for otherwise we just consider some recursively saturated elementary \zzz
extension of \A. For each \n\ in \w, $\Ex\sx\,(\lenn\ax\T^\sx\ax\tx^\sx)$ holds in \zzz
\A\ and so it must hold for some nonstandard \n\ in \A. Then \zzz
$\Ex\sx\,(\lenn\ax\T^\sx\ax(\nx\tx)^\sx)$ also holds for this \n. Let \sx\ be any witness \zzz
to this and let $\B=\A\,\up\cupx_{i\in\m\w}\sx_i$. Then \B\ is a model of \T\ and $\nx\tx$.

    Suppose $\nx\tx$ holds in \A. Let \n\ be any witness to $\nx\tx$; \n\ is \zzz
necessarily nonstandard. Then $\Ex\sx\,(\lenn\ax\T^\sx\ax\tx^\sx)$. Let \sx\ be any \zzz
witness to this and let $\B=\A\,\up\cupx_{i\in\m\w}\sx_i$. Then \B\ is a model of \T\ \zzz
and \tx.\done

   If \T\ is r.e., the above independent sentence is (equivalent to one \zzz
which is) $\Pi_2$. By using $\frac12\ast$fulfilment instead of fulfilment, we may \zzz
obtain an independent $\Pi^0_2$ sentence. However, with $\ast$fulfilment we can \zzz
obtain an independent $\Pi^0_1$ sentence, namely:

\refstepcounter{thisisdumb}\label{Chapter1:1.8}
\thm{1.8 Corollary}(Kripke)\footnote{This corollary, along with the notion of $\ast$fulfilment,
was also discovered independently by the author.} If \T\ is a consistent r.e.\ theory extending \s\ \zzz
and if
\[
\T\provex\tx\,\equiv\,\Ax n\,\big(\Ex\sx\,(\lenn\ax\T\m^{\ast\sx}\ax\tx\m^{\ast\sx})
 \imp\Ex\sx\,(\lenn\ax\T\m^{\ast\sx}\ax(\nx\tx)\m^{\ast\sx})\m\big)
\]
where ``\T'' here is understood to be any $\Dx^0_0$ formula which gives an \zzz
axiomatization of \T, then \tx\ is (equivalent to) a $\Pi^0_1$ sentence which \zzz
is independent of \T. If \tx\ is such that
\[
\T\provex\tx\,\equiv\,\Ax n\,\Ex\sx\,\big(\lenn\ax\T^{\ast\sx}\ax(\nx\tx)^{\ast\sx}\big)
\]
then \tx\ is (equivalent in \T\ to) a true $\Pi^0_1$ sentence not provable in \T.\done

\noindent
The proof is identical to that of 1.7, and is omitted. We remark that these \zzz
\tx's are not, however, very novel independent sentences: by corollary 2.5 \zzz
below, they are provably equivalent in \s\ to a Rosser sentence and to $\Con(\T)$, \zzz
respectively.

\thm{1.9 Corollary} If \T\ is, say, an r.e.\ extension of \s\ and if $k\geX 1$, \zzz
 there is a $\Pi_{k+1}$ sentence \yx\ which is $\Sx_{k+1}$-conservative over \T, and \zzz
whose negation is $\Pi_{k+1}$-conservative over \T; namely, choose any \tx\ for \zzz
which \T\ proves
\[
\tx\,\equiv\,\Ax n\,\big(\Ex\sx\,(\lenn\ax\T^\sx\ax\Tr\Sx_k\m^\sx\ax\tx^\sx)
 \imp\Ax x\,\Ex\sx\,(\lenn\ax x\e\sx_0\ax\T^\sx\ax\Tr\Sx_{k-1}\m^\sx\ax(\nx\tx)^\sx)\big)\,,
\]
and let \yx\ be the right-hand-side of this equivalence.\done

\noindent
For the definitions of the concepts involved here, and for an indication \zzz
of the proof, see Chapter~II, pages \pageref{Chapter2:DefinitionConservative}ff.

    For the remainder of the chapter we shall eschew the use of the \zzz
fixed-point theorem. The next theorem may also be obtained as a corollary of 1.5.

\refstepcounter{thisisdumb}\label{Chapter1:1.10}
\thm{1.10 Theorem}(Tarski)\footnote{Kripke has also obtained a similar proof.}
The theory of \bb N is not arithmetical.

\noindent
Proof: Suppose an arithmetical formula $\T\m(x)$ represents the theory of \bb N \zzz
in \bb N, and let $\tx(n,\sx)$ be
\begin{flalign*}
&                      & \lenn\,\ax\,& n\e\sx_0\ax\T^\sx\ax(\Tr\T)^\sx\ax\Tr\Dx_0\m^\sx\,. &&\\
&\text{Then}\hidewidth & &\bb N\modelx\Ax n\,\Ex\sx\,\tx(n,\sx)\,.\label{Chapter1:eq:6}\tag{6}
\end{flalign*}
Let \A\ be a proper elementary extension of \bb N. Choose a nonstandard \n\ \zzz
in $A$ and let $\sx'$ in $A$ be the least witness to $\Ex\sx\,\tx(n,\sx)$. Then \zzz
$\B=\A\,\up\cupx_{i\in\m\w}\sx'_i$ is also a model of the theory of \bb N which contains \n. The \zzz
formula \T\ is absolute between \A\ and \B, and hence so is \tx. But by \zzz
our choice of $\sx'$, there is no witness in \B\ for $\Ex\sx\,\tx(n,\sx)$, which \zzz
contradicts \eqref{Chapter1:eq:6}.\done

    By using the notion of $\frac12\ast$fulfilment, the above proof immediately \zzz
yields:

\thm{1.11 Corollary} If \A\ is any \w-model of \s, the theory of \A\ is not \zzz
definable over \A\ by any parameter-free formula of \sss L.\done

    D. Scott\,\cite{Scott62} was probably the first to use nonstandard models to \zzz
prove Tarski's Theorem. His proof, which is related to arguments of \zzz
Feferman and Tennenbaum, is as follows. By Henkin's proof of the \zzz
Completeness Theorem, if $\ssf{Th}(\bb N)$ is arithmetical there exists a \zzz
non-standard arithmetical model $\A=\seq{\w,0,1,\oplus,\otimes}$ of $\ssf{Th}(\bb N)$. Let
\[
\delta(n)=1\oplus1\oplus\dots\oplus1\quad \text{(\n\ times)};
\]
this is primitive recursive in the function $\oplus$. Choose a formula \yx\ \zzz
such that $\yx(n)$ holds in \bb N iff $\A\modelx\delta(n)\notE n$. Because \A\ and \zzz
\bb N are elementarily equivalent, $\bb N\modelx\yx(n)$ iff $\A\modelx\yx(\ol n)$. Let $m$ \zzz
be a (nonstandard) element of \A\ coding the set $\myset{n}{\A\modelx\yx(\ol n)}$. Then
\[
\A\modelx\ol n\e m\ \ \text{iff}\ \ \A\modelx\yx(\ol n)\ \ \text{iff}\ \ \bb N\modelx\yx(n)\ \ \text{iff}\ \ \A\modelx\ol n\notE n\,.
\]
The substitution of $m$ for \n\ gives a contradiction. This proof may \zzz
also be generalized to give 1.11 above.

\refstepcounter{thisisdumb}\label{Chapter1:1.12}
\thm{1.12 Theorem} Let $\T\supseteq\s$ be any consistent r.e.\ theory and
let $\myset{n}{\bb N\modelx\Ex m\Xx mn}$ \zzz
be any r.e.\ set, where \Xx\ is $\Dx^0_0$. Then we can feasibly semi-represent \zzz
this set in \T\ by the $\Sx^0_1$ formula
\[
\yx=\Ex m\,\Ex\sx\,\big(\lenx{\sx}{m}\ax\T^{\frac12\ast\sx}\ax\Tr\Dx^0_0\m^{\frac12\ast\sx}\ax\Xx mn\m\big)\,,
\]
where ``\T'' here is understood to be $\Dx^0_0$ formula which gives an \zzz
axiomatization of \T.\done

\noindent
The proof contains no new ideas, and is left to the reader. I am not \zzz
sure of the correct notion of \ul{feasible}; we have that
\[
\text{the length of }\yx=\text{the length of }\Xx+\text{ a constant independent of }\Xx\,;
\]
whereas a result of Parikh\,\cite{Parikh71} seems to suggest that any semi-representation \zzz
obtained in the usual manner via the fixed-point theorem would not have \zzz
this property. An extension of 1.12, giving a \yx\ in $\Dx^0_1\,(\T)$ for non-$\Sx^0_1$-sound \zzz
\T\ is contained in the proof of 2.8 in the next chapter.

    \refstepcounter{thisisdumb}\label{Chapter1:1.13}
    Our next result complements the discussion in Scott\,\cite{Scott61}.

\thm{1.13 Theorem}(Kripke, Kochen, Friedman) Let \sss G be any collection of functions \zzz
from \w\ to \w\ containing an unbounded non-decreasing function and \zzz
which is closed under composition and bounded recursion. Let \T\ be any \zzz
$\Sx^0_1$-sound set of sentences extending \PA\ whose characteristic function is \zzz
in \sss G. There exists a function $g(n,x)$ in \sss G such that if \sss F is \zzz
any collection of functions including \sss G and which is closed under \zzz
composition and bounded recursion and if $D$ is any non-principal ultrafilter \zzz
on \w, then
\[
\B=\myset{f\e\sss F}{f\text{ unary and }\Ex n\,\Ax x\,(fx\ltx gnx)}\,/\,D
\]
is a (recursively saturated) model of \T.

\noindent
Sketch of proof: Let $h(x)$ be an unbounded, non-decreasing function in \zzz
\sss G. Define:
\begin{align*}
m(x)&=
\begin{cases}
\max\,\myset{m\ltx hx}{\Ex\sx\ltx hx\,(\lenx{\sx}{m}\ax i\sx\ax\T^\sx)},\text{ if this exists,}\\
0,\text{ otherwise};
\end{cases}\\
\sx(x)&=
\begin{cases}
\min\,\myset{\sx\ltx hx}{\lenx{\sx}{mx}\ax i\sx\ax\T^\sx},\text{ if }
mx\ne0,\\
0,\text{ otherwise};
\end{cases}\\
g(n,x)&=
\begin{cases}
(\sx x)_n,\text{ if }n\lex mx,\\
0,\text{ otherwise},
\end{cases}
\end{align*}
where for this definition we are using the notion of i-fulfilment as \zzz
described on page \pageref{Chapter1:Definitionifulfilment}. For any non-principal ultrafilter $D$, \zzz
consider the structure $\A=\myset{f\e\sss F}{f\text{ unary}}/D$. Then $\sx/D$ is (the \zzz
code of) a sequence in \A\ of nonstandard length $m/D$ which i-fulfils \zzz
each standard sentence of \T. Hence $\bigcup_{n\e\,\w}\myset{x\e\A}{x\ltx(\sx/D)_n}$ is a \zzz
(recursively saturated) model of \T. But this is exactly the structure \B\ \zzz
in question.\done

    The above result is essentially due to Kripke. From the embedding \zzz
technique of Friedman\,\cite{Frie73}, one may obtain:

\thm{Theorem} If \A\ is a model of $\Dx^0_0$-overspill and \T\ is a $\Sx^0_1$-sound set of \zzz
sentences extending \PA\ coded in \A, then there exists a recursively \zzz
saturated initial segment \B\ of \A\ which is a model of \T.\done

\noindent
However, in a unpublished typescript (of which Kripke was not aware), \zzz
Friedman\,\cite{Frie7x} shows a result which to my mind is much more interesting: \zzz
one can \ul{canonically define} the initial segment. This, together with \zzz
the observation that we may take \A\ to be $\sss F/D$, will yield Theorem 1.13.
                           \cleardoublepage
\thispagestyle{plain} \sectioncentred{Herbrand's Theorem and Reflection Principles}

    In this second chapter we begin a more detailed study of the notion of \zzz
fulfilment. Using a combinatorial result from the next chapter, we use fulfilment \zzz
to derive a version of Herbrand's Theorem. As a corollary, we can \zzz
show that the notion of fulfilment is fairly stable: for example, if
\[
\provex\tx\equiv\px
\]
then in any structure \A, for all \n\ there exists a sequence of subsets \zzz
of $A$ which \n-fulfils \tx\ if and only if for all \n\ there exists such \zzz
a sequence which \n-fulfils \px; and moreover, there exists a primitive \zzz
recursive function $g(n)$ such that there is a constructive method of obtaining \zzz
a sequence \n-fulfilling \tx\ from one $g(n)$-fulfilling \px. As another \zzz
corollary, we obtain the Reflexiveness Theorem for theories such as \PA. We \zzz
go on to compare a number of minor variants of schemata involving fulfilment \zzz
with various proof-theoretic reflection principles, and then we give two \zzz
proofs of the main theorem of Kreisel and Levy\,\cite{Krei68}. We conclude the chapter by \zzz
solving a problem of Guaspari and Solovay concerning the complexity of
\[
\Myset{\tx\e\Lambda_1}{\tx\text{ is }\Lambda_2\text{-conservative over \T}}.
\]
     In this chapter $\Gamma$ will always denote an finite set of formulae and we shall \zzz
say that a sequence \sx\ fulfils $\Gamma\vx\nx\Gamma$, and
write $(\Gamma\vx\nx\Gamma)^\sx$, if \sx\ fulfils $\Av{x_{\tx}}\;(\tx\vx\nx\tx)$ \zzz
for each \tx\ in $\Gamma$, where $\vect{x_{\tx}}$ is a list of the free variables of \tx. \zzz
Let \s\ be as on page \pageref{Chapter1:DefinitionMembership}; the language of \T\ here \zzz
can be unrelated to that of \s.

    The following is the main lemma of this chapter.

\thm{2.1 Lemma} If a sentence \px\ is provable in a universal theory \T, then \zzz
there exists \n\ and a finite set of formulae $\Gamma$, depending on the proof, \zzz
such that in any model of \T\ there exists no sequence which \n-fulfils \nx\px\ \zzz
and $\Gamma\vx\nx\Gamma$.

\noindent
If we ignore the reference to $\Gamma$ (as we shall be able to do using 3.6 below) \zzz
this is just a version of Herbrand's Theorem.

    First we need a simple lemma. Define\label{Chapter2:DefinitionRank} the \ul{rank} of a formula as follows. \zzz
\tx\ is of rank $0$ if it is quantifier-free, \tx\ is of rank $k$ if \zzz
it is a Boolean combination of formulae of rank at most $k$, and $\ssf{Q}x\tx$ is of \zzz
rank $k+1$ if \tx\ is of rank $k$.

\thm{2.2 Lemma} If \tx\ is of rank $k$, then \tx\ax\nx\tx\ is not $k+1$-fulfillable \zzz
in any structure and for any valuation of its free variables.

\noindent
The proof of 2.2 is an easy induction on the complexity of \tx, but because \zzz
such a proof obscures the reason why 2.2 is true, we shall consider \zzz
instead a simple example and leave the complete proof to the reader. Let \zzz
$\tx=\Ax x\,\Ex y\,\Ax u\,\Ex v\m\m\yx$, where \yx\ is quantifier-free, and let $\sx=\seq{A_0,A_1,\dots,A_5}$. \zzz
Then $\tx^\sx$ implies
\[
\Ax x\e A_1\,\Ex y\e A_2\,\Ax u\e A_3\,\Ex v\e A_4\,\yx
\]
while $(\nx\tx)^\sx$ implies the negation. Hence not $(\tx\ax\nx\tx)^\sx$.\done

\noindent
Proof of 2.1: We use induction on the length of proof, and we shall \zzz
specify what \n\ and $\Gamma$ should be as we carry out our induction. First \zzz
we need to select some formal first-order system: for example, consider \zzz
the following Hilbert-style system.

\mypagebreak{}
\noindent
Axiom schemata:

\begingroup\leftskip=2em
\noindent
1. All tautologies. \\
2. All equality axioms. \\
3. All formulae of either the forms
\[
(\Ax x\m\tx x)\imp\tx t,\quad \tx t\imp(\Ex x\m\tx x)
\]
\quad where $t$ is any term free for $x$ in $\tx x$.
\par\endgroup

\vspace{.5em}
\noindent
Rules of Inference:

\begingroup\leftskip=2em
\noindent
1. Modus Ponens: $\quad\displaystyle\frac{\tx\quad\tx\imp\yx}{\yx}$

\noindent
2. Generalization: $\quad\displaystyle\frac{\tx\imp\yx x}{\tx\imp\Ax x\m\yx x}\quad\quad\displaystyle\frac{\yx x\imp\tx}{\Ex x\m\yx x\imp\tx}\quad$ where $x$ is not free in \tx.\\
\par\endgroup

\noindent
This system has the advantage that the only difficult step (i.e.\ the only \zzz
step involving $\Gamma$) in our induction corresponds to Modus Ponens. \zzz
(This is not to say that using, for example, a cut-free sequent calculus \zzz
would be any easier; in that case it would be the \Ex-introduction rule \zzz
which would be difficult.)

    Let \T\ be any set of universal sentences. By induction on the \zzz
length of proof, we shall show that for any formula \tx\ \zzz
there exists \n\ and $\Gamma$ such that $\Gamma\vx\nx\Gamma$ plus the negation of the \zzz
universal closure of \tx\ is not \n-fulfillable in any structure for \T.

    If \tx\ is an axiom of \T, then $\nx(\nx\tx)^\sx$ for any \sx\ because \tx\ \zzz
is universal; likewise for the equality axioms.

    Let \tx\ be a tautology. As fulfilment is preserved under truth \zzz
preserving Boolean transformations, we may consider the normal form of \zzz
\nx\tx\ consisting of a disjunction of conjuncts. Then each of these conjuncts \zzz
must contain both \yx\ and \nx\yx\ for some formula \yx. Let \n\ be greater \zzz
than the ranks of all such \yx. By 2.2, the existential closure of \tx\ \zzz
is not \n-fulfillable.

    Next consider the axiom $\tx t\imp\Ex x\m\tx$. Let the \ul{height} of a constant \zzz
or variable be $0$, and for each function symbol $f$, let the height of \zzz
$f \vect u$ be $1$ plus the maximum of the heights of its arguments. Let \vect v be \zzz
a list of the free variables of $\tx t\imp\Ex x\m\tx$ and let $\sx=\seq{A_0,\dots,A_n}$. \zzz
Unravelling the definition, $\nx\big(\nx\Av v\m(\tx t\imp\Ex x\m\tx)\big)^\sx$ iff $\Av v\ve A_1\,\big((\tx t)_0\imp\Ex i\ltx n\,\Ex x\e A_i\,\nx(\nx\tx x)_i\big)$. \zzz
Fix $\vect v\ve A_1$, and suppose $(\tx t)_0$. By Lemma 1.1.ii, $(\tx t)_k$ for all $k\ltX n$. \zzz
In particular, $(\tx t)_h$ for $h$ the height of $t$. \zzz
If \n\ is greater than $h$ plus the rank of \tx, then $\nx(\nx\tx t)_h$ \zzz
by Lemma 2.2. Thus $\Ex i\ltx n\,\Ex x\e A_i\,\nx(\nx\tx x)_i$, and we are done. Since fulfilment is \zzz
preserved by ``moving negations in or out'', the dual form follows.

    For Generalization, let $\sx=\seq{A_0,\dots,A_n}$ and note that \zzz
\begin{align*}
\nx(\nx\,\Av v,x\,(\tx\imp\yx x))^\sx
&\text{\quad iff\quad}\Av v,x\e A_1\,(\tx_0\imp\nx(\nx\yx x)_0)\\
&\text{\quad iff\quad}\Av v\ve A_1\,(\tx_0\imp\nx(\nx\Ax x\m\m\yx x)_0)\\
&\text{\quad iff\quad}\nx(\nx\,\Av v\,(\tx\imp\Ax x\m\m\yx x))^\sx.
\end{align*}

    For Modus Ponens, let \vect u, \vect v, and \vect w\ be listings of the free \zzz
variables of \yx, \tx, and \tx\m\imp\m\yx, respectively. Suppose that \tx\ is in \zzz
 $\Gamma$ and that we have a sequence $\tau$ which $n+1$-fulfils $\nx\Av u\m\yx\vect u$ and $\Gamma\vx\nx\Gamma$. \zzz
If \n\ is large enough then by 2.2, $\tau$ fulfils $\Av v\m\tx$. If we let \sx\ be \zzz
$\tau$ less its first term, say $\sx=\seq{A_0,\dots,A_n}$, then $\Ev u\ve A_0\,(\nx\yx)^\sx_0$ \zzz
and $\Av v\ve A_0\m\tx^\sx_0$. Hence $\Ev w\ve A_1\,(\tx^\sx_0\ax(\nx\yx)^\sx_0)$, that is, \sx\ \n-fulfils \zzz
$\nx\Av w\,(\tx\imp\yx)$ and $\Gamma\vx\nx\Gamma$, contradicting our induction hypothesis.

    This concludes the proof of 2.1.\done

    \refstepcounter{thisisdumb}\label{Chapter2:3.6}Next we shall quote a result, paraphrased, from the next chapter page \pageref{Chapter3:3.6}, which allows us to \zzz
eliminate the mention of $\Gamma$ in Lemma 2.1.
\begin{quotation}\thm{3.6 Theorem} For any structure \A\ for \sss L and any sentence \tx\ of \sss L, \zzz
for all \n\ there exists a sequence of subsets of $A$ which \n-fulfils \tx\ \zzz
if and only if for all \n\ and for all finite sets $\Gamma$ of formulae of \sss L \zzz
there exists such a sequence which \n-fulfils \tx\ and $\Gamma\vx\nx\Gamma$. Moreover, \zzz
the proof of this is effective in the following sense. Given \n, $\Gamma$, and \zzz
\tx, we may obtain an $m$ primitive recursively (in fact, by a function in \zzz
$\sss E_4$) from \n\ and the number of quantifiers in \tx\ and $\Gamma$ such that if \zzz
in some structure we have a sequence \sx\ which $m$-fulfils \tx, we may \zzz
constructively obtain a sequence $\sx'$ from \sx\ which \n-fulfils \tx\ and \zzz
$\Gamma\vx\nx\Gamma$. Thus we also have the formal versions:
\[
\PRA\provex\Ax\m\tx\,\ppl \,\Ax n\,\Ex\sx\,(\lenn\ax\tx^{*\sx})\,\equiv\,\Ax n,\Gamma\,\Ex\sx\,\big(\lenn\ax\tx^{*\sx}\ax(\Gamma\vx\nx\Gamma)^{*\sx}\big)\ppr
\]
where \tx\  ($\Gamma$) ranges over (the codes of) sentences (finite sets of formulae, \zzz
respectively) of \sss L, and assuming that \s\ proves the functions in $\sss E_4$ are total,
\[
\s\provex\Ax\m\tx\,\ppl \,\Ax n\,\Ex\sx\,(\lenn\ax\tx^{\sx})\,\equiv\,\Ax n,\Gamma\,\Ex\sx\,\big(\lenn\ax\tx^{\sx}\ax(\Gamma\vx\nx\Gamma)^{\sx}\big)\ppr
\]
where \tx\ and $\Gamma$ are as before but with \sss L the language of \s.\done\end{quotation}
    \refstepcounter{thisisdumb}\label{Chapter2:2.3}From 2.1 and 3.6 we may immediately obtain Herbrand's Theorem. A \zzz
formalized version of this is:

\thm{2.3 Corollary}(Herbrand)
\[
\s\provex\Ax\m\tx\,\ppl \Prx\PC(\tx)\imp\nx\Ax n\,\Ex\sx\,(\lenn\ax(\nx\tx)^\sx)\ppr
\]
where \Prx\PC\ is any natural proof predicate for the Predicate Calculus and \zzz
\tx\ ranges over (the codes of) sentences in the language of \s.\done

    \refstepcounter{thisisdumb}\label{Chapter2:2.4}
    Results 1.2 and 2.3 immediately give:

\thm{2.4 The Reflexiveness Theorem}\label{Chapter2:DefinitionDot} For each formula $\px x$ in the language of \s
\[
\s + \Induction \provex\Ax n\,\ppl \Prx\PC(\px\dot n)\imp\px n\ppr \,.
\]
where $\px\dot n$ is (the code of) the sentence obtained by substituting the \zzz
term $\ol n$ for $x$ in $\px x$.\done

    \refstepcounter{thisisdumb}\label{Chapter2:2.5}
    For any formula $\T\m(x)$ let $\Prx\T(\tx)=\Ex\text{ finite }X\subseteq\T\;\Prx\PC(\bigwedge X\imp\tx)$. \zzz
An alternative formalization of Herbrand's Theorem is:

\thm{2.5 Corollary} If \T\ is a formula in the language of \s,
\[
\s\provex\Ax\tx\,\ppl \Prx\T(\tx)\imp\nx\Ax n\,\Ex\sx\,\big(\lenn\ax\T^\sx\ax(\nx\tx)^\sx\big)\ppr \,,
\]
and if \T\ is a formula in the language of arithmetic,
\[
\PRA\provex\Ax\tx\,\ppl \Prx\T(\tx)\,\equiv\,\nx\Ax n\,\Ex\sx\,\big(\lenn\ax\T^{\ast\sx}\ax(\nx\tx)^{\ast\sx}\big)\ppr \,,
\]
where \tx\ has the appropriate ranges.

\noindent
Proof: That the converse direction also holds for $\ast$fulfilment is given \zzz
by G\"odel's proof of the Completeness Theorem as discussed on page \pageref{Chapter1:GodelsProof}.\done

    We may restate 2.5 perhaps more elegantly by using the fact that \zzz
\s\ proves:
\[
\Ax\tx\,\ppl \nx\Ax n\,\Ex\sx\,\bigl(\lenn\ax\T^\sx\ax(\nx\tx)^\sx\bigr)\,\equiv\,\Ax\px\,\bigl(\Ax n\,\Ex\sx\,(\lenn\ax\T^\sx\ax\px^\sx)\imp\Ax n\,\Ex\sx\,(\lenn\ax\T^\sx\ax\px^\sx\ax\tx^\sx)\bigr)\ppr
\]
and \PRA\ proves the $\ast$fulfilment version. The left to right direction is \zzz
immediate from 3.6 by taking $\Gamma=\{\nx\tx\}$, and the converse is obtained \zzz
from (the formalized version of) 2.2 by taking $\px\eqx \nx\tx$.

    If we wished to avoid the use of 3.6 we could, for example, restate \zzz
2.5 as
\[
\s\provex\Ax\tx\,\ppl \Prx\T(\tx)\,\equiv\,\nx\Ax n,\Gamma\,\Ex\sx\,\bigl(\lenn\ax\T^\sx\ax(\nx\tx)^\sx\ax(\Gamma\vx\nx\Gamma)^\sx\bigr)\m\ppr
\]
and similarly for $\ast$fulfilment. Likewise, in the next theorem we could \zzz
avoid the use of 3.6 by making analogous changes to all the following \zzz
schemata which involve fulfilment. Of course, if \s\ proves that the \zzz
theory \T\ contains all sentences of the form $\Av x\,(\tx\vx\nx\tx)$, then this \zzz
alteration has no effect, and it also has no effect for (a natural \zzz
axiomatization of) a theory such as \PA, where the induction schema \zzz
contains sentences which are essentially just $\Av x\,(\tx\imp\tx)$.

    \refstepcounter{thisisdumb}\label{Chapter2:DefinitionFULRFN}Next we shall use 2.5 to explore minor variants of different sentences \zzz
and schemata and to compare these with various reflection principles. To \zzz
concentrate on the case which is of most interest to us, let us suppose \zzz
that all formulae below are in the language of \s, that \s\ has a $\Dx^0_0$ \zzz
axiomatization, and that \s\ proves that \T\ is an extension of \s. First \zzz
consider the schemata:
\begin{flalign*}
&\Ful(\T)\!:           & \px\imp\Ax n\,\Ex\sx\,(&\lenn\ax\T^\sx\ax\px^\sx) && \\
&\FUL^1(\T)\!:         & \Ax m\,\ppl \px m\imp\Ax\n\,\Ex\sx\,\bigl(&\lenn\ax\T^\sx\ax(\px\dot m)^\sx\bigr)\ppr  \\
&\FUL^2(\T)\!:         & \Ax x\,\ppl \px x\imp\Ax\n\,\Ex\sx\,\bigl(&\lenn\ax\T^\sx\ax(\px     x)^\sx\bigr)\ppr  \\
&\text{and the sentences}\hidewidth \\
&\FUL^3_{\Sx_k}(\T)\!: & \Ax n  \,\Ex\sx\,(&\lenn            \ax\T^\sx\ax\Tr\Sx_k\m^\sx) \\
&\FUL^4_{\Sx_k}(\T)\!: & \Ax n,x\,\Ex\sx\,(&\lenn\ax\T^\sx\ax\Tr\Sx_k\m^\sx\ax x\e\sx_0)
\end{flalign*}
and similar sentences for $\Pi_k$. Let i-\Ful, i-$\FUL^1,\dots,\ $i-$\FUL^4_{\Sx_k}$ be the \zzz
analogous notions for i-fulfilment, and we can also consider the schemata \zzz
$\ast\Ful$ and $\ast\FUL^1$ for $\ast$fulfilment. Let $\ful$, $\ful^1,\dots,\ \ful^4_{\Sx_k}$ be the schemata
\begin{flalign*}
&\ful(\T)\!:   &              \px  \imp\Ex\sx\,\big(&\lenx\sx{\ol n}\ax\T^\sx\ax\px^\sx\big)\,,\;\;n\e\w\,, &&\\
&\ful^1(\T)\!: & \Ax m\,\ppl \px m\imp\Ex\sx\,\big(&\lenx\sx{\ol n}\ax\T^\sx\ax(\px m)^\sx\m\big)\ppr\,,\;\;n\e\w\,,\;\text{etc.}
\end{flalign*}
Finally consider \zzz
\begin{flalign*}
&\RFN(\T)\!:   & \Ax m\,(\Prx\T(\px\dot m)&\imp\px m) &&\\
&\Rfn(\T)\!:   &         \Prx\T(\px)&\imp\px\,,\;\text{and} \\
&\rfn(\T)\!:   &         \isProof\T(\ol p,\px)&\imp\px\,,\; p\e\w\,.
\end{flalign*}
\indent

\refstepcounter{thisisdumb}\label{Chapter2:FULRFN}
    The main feature of the following result, that $\FUL(\T)$ is equivalent \zzz
to $\RFN(\T)$, is due to Kripke. Suppose, for simplicity, that \T\ is $\Dx^0_0$.

\thm{2.6 Corollary} (a) In \s, for all $k$ in \w, with ``$(\T)$'' omitted for brevity:
\begin{mylist}
\item $\FUL^3_{\Sx_k    }\equiv
       \FUL^3_{\Pi_{k+1}}\equiv
       \FUL^1_{\Pi_{k+1}}\equiv
   \ast\FUL^1_{\Pi_{k+1}}\equiv
       \RFN_{\Sx_{k+1}  }$
\item $\FUL^4_{\Sx_k    }\equiv
       \FUL^4_{\Pi_{k+1}}\equiv
       \FUL^2_{\Pi_{k+1}}\equiv
       \FUL^2_{\Sx_{k+2}}\equiv
       \FUL^1_{\Sx_{k+2}}\equiv
   \ast\FUL^1_{\Sx_{k+2}}\equiv
       \RFN_{\Pi_{k+2}  }$
\item $\FUL^2_{\Dx_0    }\equiv
       \FUL^2_{\Sx_1    }\imp
       \FUL^1_{\Sx_1    }\imp
   \ast\FUL^1_{\Sx_1    }\equiv
       \RFN_{\Pi_1      }$
\item $\Ax n\,\Ex\sx\,(\lenn\ax\T^\sx)\,\imp\,
       \Ax n\,\Ex\sx\,(\lenn\ax\T^{\ast\sx})\,\equiv\,
       \Con(\T)$
\vspace{1ex}
\item $\Ful_{\Pi_k}\imp
   \ast\Ful_{\Pi_k}\equiv
       \Rfn_{\Sx_k}$
\item $\Ful_{\Sx_k}\imp
   \ast\Ful_{\Sx_k}\equiv
       \Rfn_{\Pi_k}$
\vspace{1ex}
\item $\ful^3_{\Sx_k    }\equiv
       \ful^3_{\Pi_{k+1}}\equiv
       \ful^1_{\Pi_{k+1}}\equiv
       \ful_{  \Pi_{k+1}}\equiv
       \rfn_{  \Sx_{k+1}}\equiv
   \ast\ful_{  \Pi_{k+1}}\equiv
   \ast\ful^1_{\Pi_{k+1}}$
\item $\ful^4_{\Sx_k    }\equiv
       \ful^4_{\Pi_{k+1}}\equiv
       \ful^2_{\Pi_{k+1}}\equiv
       \ful^2_{\Sx_{k+2}}\equiv
       \ful^1_{\Sx_{k+2}}\equiv
       \ful_{  \Sx_{k+2}}\equiv
       \rfn_{  \Pi_{k+2}}\equiv
   \ast\ful_{  \Sx_{k+2}}\equiv
   \ast\ful^1_{\Sx_{k+2}}$
\item $\ful^2_{\Dx_0}\equiv
       \ful^2_{\Sx_1}\imp
       \ful^1_{\Sx_1}\imp
       \ful_{  \Sx_1}\imp
       \rfn_{  \Pi_1}\equiv
   \ast\ful_{  \Sx_1}\equiv
   \ast\ful^1_{\Sx_1}$
\item $\myset{\Ex\sx\,(\lenx\sx{\ol n}\ax\T^\sx)}{n\e\w}\imp
       \myset{\Ex\sx\,(\lenx\sx{\ol n}\ax\T^{\ast\sx})}{n\e\w}\equiv
       \myset{\nx\isProof\T(\ol p,\godel{\,\text{0=1}}))}{p\e\w}$
\end{mylist}
(b) In \s\ + ``$\Collection \subseteq \T$'' (i.e.\ an axiom asserting that \T\ includes \zzz
the \Collection\ schema), the relations (i) to (x) with i-fulfilment \zzz
replacing fulfilment also hold. In fact, we may improve parts of (iii) \zzz
and (ix) by the equivalences:
\begin{flalign*}
&                     & \text{i-}\FUL^2_{\Sx_1}(\T)\,&\equiv\,\Ax n\,\Ax x\,\Ex\sx\,(\lenn                \ax x\e\sx_0\ax\text{i}\sx\ax\T^\sx) && \\
&\text{and}\hidewidth & \text{i-}\ful^2_{\Sx_1}(\T)\,&\equiv\,\myset{\Ax x\,\Ex\sx\,(\lenx\sx{\ol n}\ax x\e\sx_0\ax\text{i}\sx\ax\T^\sx)}{n\e\w}\,.
\end{flalign*}
Thus in the presence of \s\ + ``$\Collection \subseteq \T$'', fulfilment and i-fulfilment \zzz
are essentially the same notions, except in the $\Sx_1$ schemata.

\noindent
(c) The interest in (vii), (viii) and (ix) lies in the fact that \zzz
$\rfn_{\Pi_{k+1}}(\T)$ and $\rfn_{\Sx_{k+1}}(\T)$ axiomatize the $\Pi_{k+1}$ and  $\Sx_{k+1}$ consequences of \T, \zzz
respectively. The following are, in the presence of \s\ + ``$\Induction\subseteq\T$'', \zzz
further examples of such axiomatizations. First, for the $\Pi_{k+1}$ consequences, \zzz
we have:
\begin{mylist}
\setcounter{enumi}{10}
\item
$
\begin{aligned}[t]
         \myset{\Rfn_{    \Pi_{k+1}}(\T\cap\ol m)}{m\e\w}\,
&\equiv\,\myset{\ast\Ful_{\Sx_{k+1}}(\T\cap\ol m)}{m\e\w}\,,\;\text{for all $k$}\\
&\equiv\,\myset{  \;\Ful_{\Sx_{k+1}}(\T\cap\ol m)}{m\e\w}\,,\;\text{if $k\ge 1$}\\
&\equiv\,\myset{  \FUL^4_{\Sx_{k+1}}(\T\cap\ol m)}{m\e\w}\,,\;\text{if $k\ge 1$.}
\end{aligned}
$
\end{mylist}
For the $\Sx_{k+1}$ consequences, we have:
\begin{mylist}
\setcounter{enumi}{11}
\item
$
\begin{aligned}[t]
         \myset{\Rfn_{    \Sx_{k+1}}(\T\cap\ol m)}{m\e\w}\,
&\equiv\,\myset{\ast\Ful_{\Pi_{k+1}}(\T\cap\ol m)}{m\e\w}\,,\;\text{if $k\ge 1$}\\
&\equiv\,\myset{  \;\Ful_{\Pi_{k+1}}(\T\cap\ol m)}{m\e\w}\,,\;\text{if $k\ge 2$,}
\end{aligned}
$
\end{mylist}
and in the presence of ``$\Infinity\e\T$'' and ``$\Sx_{k+1}\HY\Collection\subseteq\T$'', we \zzz
also have the $\Sx_{k+1}$ form of
\[
\myset{\FUL^3_{\Sx_{k+1}}(\T\cap\ol m)}{m\e\w}\,.
\]
\noindent
(d) There are similar results for $\frac12\ast$fulfilment. In particular,
{
\setlength{\belowdisplayskip}{4pt}
\begin{flalign*}
&\tfrac12\!\ast\!\ful^3_{\Sx^0_k}(\T)\!:   & \Ex\sx\,\ppl &\lenx\sx{\ol n}\ax\T^{\frac12\ast\sx}\ax(\Tr\Sx^0_k)^{\frac12\ast\sx}\ppr ,\;\;n\e\w\,, &&\\
&\text{axiomatizes the $\Sx^0_{k+1}$ consequences of \T, and}\hidewidth \\
&\tfrac12\!\ast\!\ful^4_{\Sx^0_k}(\T)\!:   & \Ax m\,\Ex\sx\,\ppl &\lenx\sx{\ol n}\ax\T^{\frac12\ast\sx}\ax(\Tr\Sx^0_k)^{\frac12\ast\sx}\ax x\e\sx_0\ppr ,\;\;n\e\w\,,
\end{flalign*}
}
axiomatizes the $\Pi^0_{k+2}$ consequences of \T. Also, if \Induction\ is included \zzz
in \T, then $\myset{\Ax n\,\Ex\sx\,(\lenn\ax\T\cap\ol m\m\m^{\ast\sx})}{m\e\w}$ axiomatizes the $\Pi^0_1$ consequences \zzz
of \T.

\noindent
Proof: We shall only consider typical examples of nontrivial cases. \zzz
First let us consider $\FUL^2_{\Sx_1}(\T)\imp\FUL^1_{\Sx_1}(\T)$. Suppose $\Ex z\m\px z$ is true \zzz
in some structure, with \px\,\ $\Pi_1$. Choose any witness $a$ and let $b \eqx  \{a\}$. \zzz
Then any sequence which fulfils the true $\Pi_1$ formula
\[
\Ex x\e b\,(x=a\ax\px a)
\]
will fulfil $\Ex z\m\px z$.

    That $\ast\FUL(\T)$ implies $\RFN(\T)$ (and conversely) is an immediate \zzz
consequence of 2.5.

    Finally, suppose $\RFN_{\Pi_{k+1}}(\T)$ and that we wish to show $\FUL^4_{\Sx_{k-1}}(\T)$, \zzz
say. For each \n\ in \w\ we have
\[
\T\provex\Ax x\,\Ex\sx\,(\lenm\sx\mx\eqx\ol n\m\ax x\e\sx_0\ax\T^\sx\ax\Tr\Sx_{k-1}\m^\sx)\,.
\]
Moreover, this fact can be formalized, and so we have
\[
\s\provex\Ax n\,\Prx{\T}\ppl \Ax x\,\Ex\sx\,(\lenn\ax x\e\sx_0\ax\T^\sx\ax\Tr\Sx_{k-1}\m^\sx)\m\ppr \,.
\]
Hence by $\RFN_{\Pi_{k+1}}(\T)$ we obtain $\FUL^4_{\Sx_{k-1}}(\T)$.\done

    We cannot in general improve (iv), (v) or (vi). A result of \zzz
Feferman\,\cite{Fefe60} is that
\[
\s+\RFN_{\Sx^0_1}(\T)\provex\RFN_{\Sx^0_1}(\T+\Rfn(\T)).
\]
(Aside: this can be easily generalized from $\Sx^0_1$ to arbitrary $\Lambda\supseteq\Sx^0_1$. A \zzz
simple proof is given in Smory\'nski\,\cite{Smor77}; or by an equally simple proof, \zzz
using 2.5 we can show
\[
\s+\FUL^3_\Lambda(\T)\provex\FUL^3_\Lambda(\T+\Rfn(\T))\,.\,)
\]
Now by Matijacevic's Theorem and 2.5 we have that for any arithmetical \zzz
theory \T\ in the language of arithmetic,
\[
\s+``\Sx_1\HY\Induction\subseteq\T"\provex\ppl \Ax n\,\Ex\sx\,(\lenn\ax\T^\sx)\m\ppr \,\equiv\,\FUL^3_{\Pi^0_1}(\T)
\]
and so if \T\ is $\Dx^0_0$, then by either G\"odel's Second Incompleteness \zzz
Theorem or 1.5,
\[
\s+\Rfn(\T)\nprovex\Ax n\,\Ex\sx\,(\lenn\ax\T^\sx)\,.
\]
\indent
    \refstepcounter{thisisdumb}\label{Chapter2:2.7}
    Our next result is a version of the main lemma of Kreisel and Levy\,\cite{Krei68}.

\thm{2.7 Corollary} Let \ssf U extend \s\ and let \T\ be a subset of \ssf U which \zzz
is representable in (some consistent extension of) \ssf U. Suppose that \zzz
\ssf U implies $\RFN(\T)$, where ``\T'' in this schema is understood to be some \zzz
representation of \T. Then no consistent extension of \ssf U (plus \zzz
\myset{\T(\godel\tx)}{\tx\e\T}) may be obtained by adding a set of $\Sx_k$ sentences \zzz
to \T\ for any $k$.

\noindent
Proof: We shall sketch two proofs, the second uses the fixed-point \zzz
theorem while the first does not, but the first proof requires the \zzz
further assumption that \ssf U extends \Foundation. Fix $k$ in \w.

    By 2.6 we have that \ssf U implies
\[
\Ax n\,\Ex\sx\,(\lenn\ax n\e\sx_0\ax\T^\sx\ax\Tr\Sx_k\m^\sx)\,.\label{Chapter2:eq:1}\tag{1}
\]
We can suppose that $\T(x)$ is $\Pi_k$. Then by fixing a nonstandard \n\ and \zzz
choosing a \ul{minimal} witness \sx\ as in 1.3, we see by our usual method that \zzz
any nonstandard model of $\ssf U+\myset{\T(\godel\tx)}{\tx\e\T}$ has a $k$-elementary \zzz
substructure which is a model of \T\ and in which \eqref{Chapter2:eq:1} is false.

    For our second proof, let $\ssf U_0$ be a finite subset of  \ssf U which \zzz
suffices to prove that some formula $\Tr_{\Sx_k}(x)$ is a $\Sx_k$ truth predicate. \zzz
By 2.6, for all sentences \Xx
\[
\ssf U\provex\Xx\imp\Ax\px\,\big(\Tr_{\Sx_k}(\px)\imp\Ax n\,\Ex\sx\,(\lenn\ax(\Tr_{\Sx_k}(\px))^\sx\ax\T^\sx\ax\ssf U_0\m^\sx\ax\Xx^\sx)\bigr)\,. \label{Chapter2:eq:2}\tag{2}
\]
Choose \Xx\ such that
\[
\ssf U\provex\Xx\equiv\Ax\px\,\big(\Tr_{\Sx_k}(\px)\imp\Ax n\,\Ex\sx\,(\lenn\ax(\Tr_{\Sx_k}(\px))^\sx\ax\T^\sx\ax\ssf U_0\m^\sx\ax(\nx\Xx)^\sx)\bigr)\,.
\]
Then $\ssf U\provex\Xx$, for by \eqref{Chapter2:eq:2}, $\ssf U\provex\nx\Xx\imp\Xx$. But for any $\Sx_k$ sentence \px\ \zzz
consistent with $\ssf U+\myset{\T(\godel\tx)}{\tx\e\T}$, we may find a model of \zzz
$\Tr_{\Sx_k}(\px)$, \T, $\ssf U_0$, and $\nx\Xx$. Hence $\T+\px\,\nprovex\Xx$.\done

   \refstepcounter{thisisdumb}\label{Chapter2:2.8}
    We conclude the chapter by considering a problem posed by D.\ Guaspari\,\cite{Guas79}. \zzz
Given a class of sentences $\Lambda$ and a set of sentences \T, say \zzz
that a sentence \px\ is $\Lambda$-\ul{conservative}\label{Chapter2:DefinitionConservative} \ul{over} \T\ if \zzz
for all \tx\ in $\Lambda$, if $\T+\px\provex\tx$ then $\T\provex\tx$ already. For a \zzz
theory \T\ let
\[
(\Lambda_1,\Lambda_2)=\myset{\godel\px\e\Lambda_1}{\px\text{ is $\Lambda_2$-conservative over \T}}\,.
\]
R.\ Solovay and P.\ H\'ajek have shown that if \T\ is a consistent r.e.\ \zzz
reflexive theory then $(\Sx_k,\Pi_k)$ is a complete $\Pi^0_2$ set for all $k\ge 1$. \zzz
We shall prove this without the hypothesis of reflexivity and we shall \zzz
also consider the dual case. (In retrospect, I see that simple modifications \zzz
of H\'ajek's\,\cite{Hajek78} proof will give most of the results in (b) below: just use the \zzz
equivalence of the notions of reflection and fulfilment as given in 2.6.) \zzz
Part (a) below is due to Guaspari\,\cite{Guas79} who uses a construction due to Scott\,\cite{Scott62} \zzz
and Friedman\,\cite{Frie73}. He requires, however, that \T\ contain the power-set \zzz
axiom.

\thm{2.8 Corollary}\label{Chapter2:Corollary28} (a)\, Let \T\ be a consistent extension of \s. A sentence \px\ \zzz
is $\Sx_{k+1}$-conservative ($\Pi_{k+2}$-conservative) over \T\ iff all models of \zzz
\T\ + $\Sx_{k+1}$-overspill which code \T\ have (for any element $x$) a $k$-elementary \zzz
substructure (containing $x$) which is a model of \T\ + \px.

\noindent
(b)\, Let T be an r.e.\ extension of \s. Then: \zzz
\begin{mylist}
\item For all $k\ge 1$, $(\Sx_k,\Pi_k)$ is a complete $\Pi^0_2$ set.
\item For all $k\ge 2$, $(\Pi_k,\Sx_k)$ is a complete $\Pi^0_2$ set, as are $(\Pi_2,\Sx_1)$ \zzz
and $(\Sx_2,\Sx_1)$, and if \T\ proves \Infinity\, so is $(\Pi_1,\Sx_1)$.
\item $(\Dx^0_2(\T),\,\Sx^0_1)$ is a complete $\Pi^0_2$ set if T is not $\Sx^0_2$-sound, and is \zzz
$\Dx^0_2$ otherwise.
\item $(\Pi^0_1,\Sx^0_1)$ and $(\Pi^0_1,\Dx^0_1(\T))$ are complete $\Pi^0_2$ sets if \T\ is not \zzz
$\Sx^0_1$-sound, and are complete $\Pi^0_1$ sets otherwise.
\end{mylist}
Proof: Part (a) is immediate from 2.6, so let us consider (b). It is \zzz
clear by their definition that these sets are $\Pi^0_2$. It is also easy to \zzz
check that if \T\ is $\Sx^0_1$-sound, a sentence is $\Dx^0_1(\T)$-conservative over \zzz
\T\ iff it is consistent with \T, a $\Pi^0_1$ sentence is $\Sx^0_1$-conservative iff \zzz
it is true, and if \T\ is also $\Sx^0_2$-sound, a $\Pi^0_2$ sentence is $\Sx^0_1$-conservative \zzz
over \T\ iff it is true. This establishes the negative parts of \zzz
(iii) and (iv).

    We shall next show that ($\Pi_2,\Sx_1)$ and $(\Pi_{k+1},\Sx_{k+1})$ are complete \zzz
$\Pi^0_2$ sets for $k\ge 1$. We shall in fact give two proofs: the first avoids \zzz
the use of the fixed-point theorem but requires the further assumption \zzz
that \T\ extends \Foundation. Let $k\ge 1$, let $\yx n$ be $\Sx^0_1$, and let \zzz
$\Phi$ be the $\Pi_{k+1}$ sentence
\[
\Phi=\Ax n\,\ppl \Ex\sx\,(\lenn\ax\T^\sx\ax\Tr\Sx_k\m^\sx)\imp\yx n\ppr \,.
\]
Consider the following statements:
{\begin{mylist}
\parskip0pt
\itemindent.8em
\renewcommand{\labelenumi}{(\alph{enumi})}
\item for all $n\e\w$, $\T\provex\yx\ol n$;
\item $\Phi$ is $\Sx_{k+1}$-conservative over \T;
\item $\Phi$ is $\Pi_k    $-conservative over \T;
\item $\Phi$ is $\Sx^0_1  $-conservative over \T;
\item $\Phi$ is $\Dx^0_1(\T)$-conservative over \T.
\end{mylist}}
\noindent
We claim that (a), (b) and (c) are equivalent. \zzz
The result follows from this claim, for it \zzz
is easy to show by the reader's favourite method (e.g.\ 1.12) that
\[
\myset{\godel{\yx x}\e\Sx^0_1}{\Ax n\e\w,\,\T\provex\yx\ol n}
\]
is a complete $\Pi^0_2$ set.

    To prove the claim, note that (b)\imp(c) is trivial , and that (c)\imp(a) \zzz
is easy. So consider (a)\imp(b). Suppose (a). Let \A\ be any \zzz
non-\w-model of \T\ (+\Foundation). It suffices to show there exists a \zzz
$k$-elementary substructure \B\ of \A\ which is a model of \T\ and $\Phi$. \zzz
If $\Phi$ holds in \A, we are done. If not, choose the least witness $m$ \zzz
of $\nx\Phi$; $m$ is necessarily nonstandard by (a). Now choose a \ul{minimal} witness \zzz
\sx\ as in 1.3, and let $B=\cupx_{i\in\m\w}\sx_i$. Then $\B\!\prec_k\!\A$, and $\B\modelx\T$. It remains \zzz
to show that $\Phi$ holds in \B. Consider $\n\e B$. If $n<m$, then $\yx n$ holds \zzz
in \A\ by the choice of $m$, and so $\yx n$ holds in \B. If $n\ge m$, then \zzz
by the choice of \sx\ the antecedent in $\Phi$ must be false at \n\ in \B. \zzz
Hence $\Phi$ holds in \B.

    To prove that $(\Pi_{\max(k+1,2)},\Sx_{k+1})$ is a complete $\Pi^0_2$ set if \T\ does \zzz
not necessarily extend \Foundation, choose $\Phi$ such that \T\ proves:
\[
\Phi\equiv\Ax n\,\ppl \Ex\sx\,(\lenn\ax\T^\sx\ax\Tr\Sx_k\m^\sx\ax\Phi^\sx)\imp\yx n\ppr \,.\label{Chapter2:eq:3}\tag{3}
\]
Then, as before, (a), (b) and (c) are equivalent. Since the choice of $\Phi$ \zzz
can be made recursive in $\godel\yx$, we are done. Likewise, if \T\ proves \zzz
\Infinity, $(\Pi_1,\Sx_1)$ is a complete $\Pi^0_2$ set. (Note that the analogous \zzz
fixed-point and argument using the notion of $\ast$fulfilment does not give as \zzz
sharp a result: we obtain that $(``\m\Ax n\,\Ex m\ltx n\,\Pi_{k+1}\mx",\,\Sx_{k+1})$ is a complete $\Pi^0_2$ \zzz
set, where this first class is $\Pi_{k+1}$ prefixed by a universal and a bounded \zzz
existential numerical quantifier.)

    To show that $(\Sx_2,\Sx_1)$ is a complete $\Pi^0_2$ set, let $\yx m=\Ex l\,\yx_0\m ml$ be $\Sx^0_1$ \zzz
and consider the fixed point in \T:
\[
\Phi\equiv\nx\Ax n\,\ppl \Ex l\,\Ax m\ltx n\,\yx_0\m ml_m\,\imp\,\Ex\sx\,(\lenn\ax\T^\sx\ax\Tr\Dx^0_0\m^\sx\ax\Phi^\sx)\m\ppr \,.
\]
$\Phi$ is (equivalent to) a $\Sx_2$ sentence, and as above, (a), (b) (with $k\eqx 0$) \zzz
and (c) are equivalent.

    Next suppose \T\ is not $\Sx^0_2$-sound, and choose some true $\Pi^0_2$ sentence \zzz
$\Ax m\,\Ex l\,\tx ml$ which \T\ refutes. Let $\yx m=\Ex l\,\yx_0\m ml$ be $\Sx^0_1$, and let \zzz
$\Xx n=\Ex l\,\Ax m\ltx n\,\yx_0\m ml_m\ax\Ex l\,\Ax m\lex n\,\tx ml_m$, and suppose \T\ proves:
\[
\Phi\equiv\Ax n\,\ppl \Ex\sx\,(\lenn\ax\T^{\frac12\ast\sx}\ax\Tr\Dx^0_0\m^{\frac12\ast\sx}\ax\Phi^{\frac12\ast\sx})\,\imp\,\Xx n\ppr \,.
\]
Then (a) and (d) are equivalent. Since \T\ proves $\nx\Ax n\m\Xx n$, \T\ proves:
\[
\Phi\equiv\Ex n\,\ppl\nx\Ex\sx\,\big(\lenn\ax\T^{\frac12\ast\sx}\ax\Tr\Dx^0_0\m^{\frac12\ast\sx}\ax\Phi^{\frac12\ast\sx}\m\big)\ax\Xx(n-1)\m\ppr \,.
\]
Thus $\Phi$ is $\Dx^0_2(\T)$, and so $(\Dx^0_2(\T),\,\Sx^0_1)$ is a complete $\Pi^0_2$ set.

    Suppose now that \T\ is not $\Sx^0_1$-sound. D.\ Jensen and A.\ Ehrenfeuch\,\cite{Jens76} \zzz
and D.\ Guaspari\,\cite{Guas79} have shown by the fixed-point theorem that an r.e.\ \zzz
theory \T\ is $\Sx^0_1$-sound iff it decides every $\Dx^0_1(\T)$ sentence. Using a \zzz
similar calculation, we can show that if an r.e.\ theory \T\ is not $\Sx^0_1$-sound, \zzz
then any r.e.\ set may be semi-represented in \T\ by a formula which is $\Dx^0_1$ \zzz
in \T. We can, however, do better by giving a feasible representation as \zzz
in 1.12. Let $\Ex k\m\Xx\mx k$ be a false $\Sx^0_1$ sentence provable in \T, with $\Xx$ \zzz
$\Dx^0_0$, let $\tx mn$ be any $\Dx^0_0$ formula, and let
\[
\yx n=\Ex m,\sx\,(\lenx\sx m\ax\T^{\frac12\ast\sx}\ax\Tr\Dx^0_0\m^{\frac12\ast\sx}\ax\tx mn\ax\Ax k\ltx\godel{\sx\!}\,\nx\Xx k)\,.
\]
Then, since \s\ is assumed to include $\Dx^0_1(\s)$-\Induction, it is easy to \zzz
show that for all $n\e\w$, $\T\provex\yx\ol n$ iff $\Ex m\m\m\tx mn$ is true. Moreover, \zzz
because \T\ proves $\Ex k\m\Xx k$, \T\ also proves
\[
\Ax n\,\ppl \nx\yx n\equiv\Ex k\,\big(\Xx k\ax\nx\Ex m,\sx\lex k\,(\lenx\sx m\ax\T^{\frac12\ast\sx}\ax\Tr\Dx^0_0\m^{\frac12\ast\sx}\ax\tx mn)\m\big)\m\ppr \,.
\]
and so \yx\ is $\Dx^0_1(\T)$. Now if \T\ contains $\Dx^0_1(\T)$-\Induction, consider
\[
\Phi=\Ax n\,\ppl \Ex\sx\,(\lenn\ax\T^{\frac12\ast\sx}\ax\Tr\Dx^0_0\m^{\frac12\ast\sx})\imp\yx n\ppr \,;
\]
and if not, consider a fixed point similar to \eqref{Chapter2:eq:3}. Then we have that (a), (d) \zzz
and (e) are equivalent, and so $(\Pi^0_1,\Sx^0_1)$ and $(\Pi^0_1,\Dx^0_1(\T))$ are complete \zzz
$\Pi^0_2$ sets.

    Finally, for the sets $(\Sx_{k+1},\Pi_{k+1})$ consider either the fixed points:
\[
\Phi\,\equiv\,\nx\Ax n\,\ppl \Ax m\ltx n\,\yx m\imp \Ax x\,\Ex\sx\,(\lenn\ax x\e\sx_0\ax\T^\sx\ax\Tr\Sx_{k-1}\m^\sx\ax\Phi^\sx)\m\ppr \,;
\]
(which work for all $k\ge1$); or the fixed points:
\[
\Phi\,\equiv\,\nx\Ax n\,\ppl \Ax m\ltx n\,\yx m\imp \Ax\tx\,\big(\Tr_{\Sx_{k+1}}(\godel\tx)\imp\Ex\sx\,(\lenn\ax\T^{\ast\sx}\ax\tx^{\ast\sx}\ax\Phi^{\ast\sx})\m\big)\m\ppr \,,
\]
(which work for all $k$). The details are left to the reader.\done

    Theorem 6.12 on page \pageref{Chapter6:6.12} gives an analogue of this last \zzz
result for the \w-rule.
     \cleardoublepage
\thispagestyle{plain} \sectioncentred{The Substitution Method}

\refstepcounter{thisisdumb}\label{Chapter3:HerbrandAckermannScanlon}
Hilbert posed the following problem: can one prove the consistency \zzz
of Peano arithmetic in a constructive manner? He was interested in this \zzz
question because a positive answer implies that the (nonconstructive) \PA\ is \zzz
a conservative extension of the constructive techniques used with respect \zzz
to \ul{real} (that is, $\Pi^0_1$) sentences. (See the introduction of Smory\'nski\,\cite{Smor77} \zzz
for a further discussion.) As is well-known, solutions to this problem \zzz
were given by Gentzen, in \cite{Gent36} and \cite{Gent38}, who used as his main tool (in the \zzz
second paper at least) \ul{cut-elimination}, and by Ackermann in \cite{Acke40}, who \zzz
used the Hilbert substitution method. As Dreben and Denton\,\cite{Dreb70} point \zzz
out, the Gentzen approach, with its purely syntactical transformations \zzz
of formal proofs, has little use for the \ul{interpretations} of formulae; \zzz
the Ackermann approach, on the other hand, ``\kern .05em\lq finitistically\rq exploits the \zzz
oldest and most na\"ive idea in proof theory: a set of axioms is consistent \zzz
if it has a model''. An elegant extension of Ackermann's proof was given \zzz
in Tait\,\cite{Tait65a},\,\cite{Tait65b} using \ul{functionals}: Tait's proof applies to systems of \zzz
arithmetic with the schema of foundation on some arbitrary primitive \zzz
recursive linear ordering. Scanlon\,\cite{Scan73} building on the work of \zzz
Dreben and Denton\,\cite{Dreb70} obtains the same results for linear orderings \zzz
without the use of functionals; unfortunately, Scanlon's paper is very \zzz
long and intricate, and contains, as Kripke has pointed out, some unnecessary  \zzz
detours. This chapter is essentially an exposition of the above-mentioned \zzz
work of Ackermann, Dreben, Denton and Scanlon.

    The notion of fulfilment is of course not necessary in our treatment \zzz
(indeed one could regard it as superfluous since all the work is still \zzz
done with finite versions of Skolem functions) but I believe that it is \zzz
very helpful conceptually---especially as it provides a clear separation \zzz
between the logical and combinatorial parts of our argument.

    We shall assume that \s\ on page \pageref{Chapter1:DefinitionMembership} is finite. Let be \prece\ any formula of \zzz
\sss L, the language of \s, with two free variables; this should be thought \zzz
of as representing a pre-well-founded relation, i.e., the non-linear \zzz
analogue of a pre-well-ordering. Let $\prece^\text{strict}$ denote $x\prece y\ax y\nprece x$; \zzz
this will also be denoted by $\prec$. Let $\Foundation(\prece)$ be the schema
\begin{flalign*}
&\Foundation(\prece)\!:  &&\Ex x\m\tx x\imp\Ex x\,\ppl \tx x\ax\Ax y\,(\tx y\imp y\nprec x)\ppr\,.&&&&&
\end{flalign*}
Consider the consistency of the theory $\s+\Foundation(\prece)$. \zzz
By 2.6, it suffices to show that each finite subset of this theory plus any true sentence is \zzz
\n-fulfillable for all \n; and we would like to do this in as constructive \zzz
a manner as possible. Of course, by the Second Incompleteness Theorem \zzz
or by 1.5, we must assume a quite strong notion of constructiveness: \zzz
our proof may be formalized in an extension $\Foundation(\prece_\ee)$ defined \zzz
below. If the axioms of \s\ are $\Sx_3$ and the relation $\prece$ is $\Sx_2$, we \zzz
need only the schema $\Foundation(\prece_\ee)$ restricted to $\Sx_1$ formulae, and \zzz
if the axioms of \s\ have the form
\[
\Ev x\,\Av y\,\Ex!\xvect z\m\tx
\]
with \tx\ $\Dx_0$, and if $\prece$ has the form $\Av y\,\Ex!\xvect z\m\tx$, with \tx\ $\Dx_0$, then \zzz
we need only the schema ${\Foundation(\prece_\ee)}$ restricted to $\Dx_1(\s)$ predicates.

    We need some general definitions. Given a binary relation $\prece$, let \zzz
$\prec=\prece^\text{strict}$ be as above. Define a binary relation $\prece^\text{lex}$ on nonempty \zzz
finite sequences by $\seq{a_0\m,\dots,a_m}\prece^\text{lex} \seq{b_0\m,\dots,b_n}$ iff \ul{either} \zzz
there exists $k\le\min(m,n)$ such that $a_l\eqx b_l$ for all $l<k$ and \zzz
$a_k\prec b_k$, \ul{or} $m<n$ and $a_l\eqx b_l$ for all $l\le m$. Let $\prece^\w$ be the \zzz
disjoint union over \n\ of $\prece^\text{lex}$ restricted to sequences of length \n: \zzz
that is, $x\prece^\w y$  iff $|x|\eqx |y|$ and $x\prece^\text{lex} y$. A finite sequence $x$ \zzz
is \prece-\ul{descending} if $x_{i+1}\prec x_i$ for all $i\ltx\lenm x$; let $2^\prece$ denote $\prece^\text{lex}$ \zzz
restricted to \prece-descending sequences. If $\prec$ is well-founded and \zzz
nontrivial, $\prece^\text{lex}$ is \ul{not} well-founded, but the usual proofs (see, for \zzz
example Feferman\,\cite{Fefe77}) show that $\Foundation(\prece)$ implies $\Foundation(\prece^\w)$ and $\Foundation(2^\prece)$.

    Let us go back to our given \prece. Fix any element and denote it \zzz
by $\infty$. Define the relations:
\begin{align*}
x\prece_+ y   \text{ \ if \ }&y=\infty\text{ \ul{or} both }x\ne\infty\text{ and }x\prece y\,,\\
x\prece_0 y   \text{ \ if \ }&x\,(\prece_+)^\w\,y\,,\\
x\prece_{n+1}y\text{ \ if \ }&x\,(2^{\prece_n})\,y\,,
\end{align*}
and let $\prece_\ee$ be a disjoint union over \n\ of $\prece_n$: e.g.\ say
\[
(n,x)\prece_\ee (m,y)\text{ \ if \ }n=m\text{ and }x \prece_n y.
\]
We shall suppose that this relation may be represented in the language \zzz
of \s\ in some natural way. We note that if \prece\ is $\Dx_k(\s)$ then $\prece_\ee$ \zzz
is $\Dx_k(\s+\Sx_{k+1}\HY\Induction)$.

    Recall one of the twenty-odd schemata given on page \pageref{Chapter2:DefinitionFULRFN}, there labelled $\FUL^2$\,:
\begin{flalign*}
&\FUL(\T)\!:         & \Ax x\,\ppl \px x\imp\Ax\n\,\Ex\sx\,\big(&\lenn\ax\T\m^\sx\ax(\px     x)^\sx\big)\ppr\,. &&&&
\end{flalign*}
\thm{3.1 Theorem} $\s+\Foundation\m(\prece_\ee)\,\provex\,\FUL\m\big(\s+\Foundation\m(\prece)\m\big)$.

    We shall give an informal proof, which will be divided into four \zzz
lemmas.

    With each formula \tx\ with \n\ free variables associate an \n-ary \zzz
function symbol $f_\tx$. Define the \ul{Skolemization} $\sss S(\tx)$ of \tx\ by \zzz
induction on length:
\begin{mylist}
\item if \tx\ is an atom, let $\sss S(\tx)=\tx$;
\item let $\sss S(\tx\imp\yx)=\sss S(\tx)\imp\sss S(\yx)$, and similarly for the other Boolean connectives; and
\item let $\sss S(\ssf{Q}x\m\tx)=\sss S(\tx)(x/f_{\ssf{Q}x\m\tx}\vect v)$ where $\vect v$ is a listing of the free variables of $\ssf{Q}x\m\tx$.
\end{mylist}

    For the remainder of the proof fix an integer \n\ and a finite set \zzz
$H$ of formulae of the form $\Ex x\m\tx$. We shall first consider the problem \zzz
of finding a sequence which \n-fulfils the universal closure of $\Foundation(\prece)$
\[
\Ex x\m\tx\mx x\imp\Ex x\,\ppl \tx x\ax\Ax y\,(\tx y\imp y\nprec x)\ppr  \label{Chapter3:eq:1}\tag{1}
\]
for each $\Ex x\m\tx$ in $H$.

    Let $\ssf{sub}(H)$ be the collection of subformulae of formulae of $H$. If \zzz
$\ssf{Q}x\m\tx$ in $\ssf{sub}(H)\setminus H$ is of rank $r$ (as defined in Chapter~II page \pageref{Chapter2:DefinitionRank}), let the \zzz
\ul{rank} of $f_{\ssf Qx\m\tx}$ be $r$. Let $q-1$ be the maximum of the ranks of formulae \zzz
in $\ssf{sub}(H)\setminus H$, and for all $\Ex x\m\tx$ in $H$ let the rank of $f_{\Ex x\m\tx}$ be $q$. \zzz
Let the functions of \sss L have rank $1$. Let
$\sss F=\myset{f_{\ssf{Q}x\m\tx}}{\ssf{Q}x\m\tx\e\ssf{sub}(H)}$.

\thm{3.2 Lemma} Suppose we are given an interpretation of the functions of \sss F \zzz
and a set $A_0$ which contains (the interpretations of) the constants of \zzz
\sss L. Define by $\sx=\seq{A_i}_{i\lex n}$ by
\[
A_{i+1}=A_i\cup\cupx_{f\e\sss L\cup\sss F}f''A_i.
\]
Suppose:
\begin{mylist}
\item for each $\Ex\m\tx$ in $H$ with free variables $\vect u$,
\[
\Av u,x\,(\sss S(\tx)\imp x\nprec f_{\Ex x\m\tx}\,\vect u)
\]
where here, as in (ii) and (iii) below, \ul{the quantifiers} $\Av u$, $\Ex x$, \ul{and} \zzz
$\Ax x$ \ul{are restricted to} $A_{n-1}$;
\item for each $\Ex x\m\tx$ in $\ssf{sub}(H)$ with free variables $\vect u$,
\[
\Av u\,\ppl \Ex x\,\sss S(\tx)\imp\sss S(\Ex x\m\tx)\ppr ;
\]
\item for each $\Ax x\m\tx$ in $\ssf{sub}(H)$ with free variables $\vect u$,
\[
\Av u\,\ppl \sss S(\Ax x\m\tx)\imp\Ax x\,\sss S(\tx)\ppr ;\text{ and}
\]
\item for all $i\ltx n$, all $x,y$ in $A_i$, if $x\nprec y$ then $\seq{A_i,\dots,A_n}$ \zzz
fulfils $x\nprec y$ (in the sense of the second definition of fulfilment motivated by Skolem functions, as \zzz
given on page \pageref{Chapter1:DefinitionFulfilIII}).
\end{mylist}
Then for each $\Ex x\m\tx x$ in $H$ with free variables $\vect u$, \sx\ fulfils (in the \zzz
sense of this second definition) \eqref{Chapter3:eq:1} above.

\noindent
Proof: We shall first show by induction on length that if \tx\ in $\ssf{sub}(H)$ \zzz
has free variables $\vect u$, then for all $i\ltx n$, all $\vect u$ in $A_i$,
\[
\nx(\nx\tx)_i\imp\sss S(\tx)\text{ \ and \ }\sss S(\tx)\imp\tx_i\,.\label{Chapter3:eq:2}\tag{2}
\]
If \tx\ is quantifier-free this is trivial, and the inductive step is \zzz
easy for the Boolean connectives. So consider $\Ex x\m\tx$ in $\ssf{sub}(H)$, and \zzz
suppose the claim is true for \tx. Fix $i$ and the free variables $\vect u$ \zzz
of $\Ex x\m\tx$ in $A_i$.

    Suppose $\nx(\nx\Ex x\m\tx)_i$. By definition, there exists a $j$ with $i\lex j\ltx n$ \zzz
and an $x\e A_j$ such that $\nx(\nx\tx)_j$. By the inductive hypothesis, $\sss S(\tx)$. \zzz
Then $\Ex x\e A_{n-1}\,\sss S(\tx)$, and so by (ii), $\sss S(\Ex x\m\tx)$.

    Now suppose $\sss S(\Ex x\m\tx)$, that is, $\sss S(\tx)(x/f_{\Ex x\tx}\xvect u)$. Since $f_{\Ex x\tx}\xvect u\ve A_{i+1}$ \zzz
if $i+1\ltx n$ the inductive hypothesis gives $\tx_{i+1}(x/f_{\Ex x\tx}\xvect u)$, and so \zzz
$(\Ex x\m\tx)_i$. If $i+1\eqx n$ we have $(\Ex x\m\tx)_i$ anyway.

    The proof for formulae of the form $\Ax x\m\tx$ is dual to the above, and \zzz
is left to the reader. This establishes the claim.

    Let $\Ex x\m\tx$ in $H$ have free variables $\vect u$. We wish to show
\[
\Ax i\ltx n\,\Av u\ve A_i\,\PPL\nx(\nx\Ex x\m\tx)_i\imp\Ex x\e A_{i+1}
   \ppl(\tx x)_{i+1}\ax\Ax j\m.\,i\ltx j\ltx n\,\big(\nx(\nx\tx y)_j\imp(y\nprec x)_j\m\big)\ppr\PPR\,.
\]
Fix $i<n$ and $\vect u$ in $A_i$ and suppose $\nx(\nx\Ex x\m\tx)_i$. By $\eqref{Chapter3:eq:2}$, $\sss S(\Ex x\m\tx)$. \zzz
Let $x=f_{\Ex x\m\tx}(\xvect u)$. Then $\sss S(\tx)(x)$, and by $\eqref{Chapter3:eq:2}$ or by $i+1\eqx n$, $(\tx x)_{i+1}$. Fix \zzz
$j$ with $i<j<n$ and $y\e A_j$, and suppose $\nx(\nx\tx y)_j$. By $\eqref{Chapter3:eq:2}$, $\sss S(\tx)(x/y)$, \zzz
and so by (i), $y\nprec x$. By virtue of (iv), we are done.\done

    If $\xvect x$ is a $k$-tuple and $f$ is an $m$-ary function, let $\hat f''\xvect x$ denote (for this paragraph only) \zzz
the $k^m$-tuple $\seq{f\xvect y_i}_{i<k^m}$, where $\seq{\vect y_i}_{i<k^m}$ is some standard listing \zzz
of the $m$-tuples from $\vect x$. Given an interpretation of \sss F and a finite \zzz
sequence $\vect A_0$, define a sequence of sequences $\vect\sx=\seq{\vect A_0\m,\dots,\vect A_n}$ inductively \zzz
by
\[
\vect A_{i+1}=\vect A_i\m\smile\hat f''\xvect A_i\m\smile\hat g''\xvect A_i\m\smile\dots
\]
where $f,g,\dots$ is some fixed listing of the functions of $\sss L\cup\sss F$, \zzz
and where ``$_\smile$'' is sequence concatenation.

    A \ul{loop} is a finite sequence $a_0\m,\dots,a_{m+1}$ such that for all $i\le m$, \zzz
$a_i\prec a_{i+1}$ and $a_{m+1}\prec a_0$; $\Foundation(\prece)$ proves that no loops exist.

\thm{3.3 Lemma} There exists a constructive interpretation of \sss F which \zzz
satisfies the premises of Lemma 3.2.

\noindent
Proof: We shall first consider only the premises (i), (ii) and (iii). A \zzz
simple modification of our algorithm will then take care of (iv); this \zzz
is discussed following Lemma 3.5.

    We shall define a sequence $\seq{\sss F^i}_{i<\Omega}$ (where $\Omega\le\w$) of interpretations of \sss F, \zzz
where $\sss F^i=\myset{f^i}{f\e\sss F}$. Each $f^i$ will have constant value $\infty$ except at \zzz
a finite number of arguments. For each $\sss F^i$, define $\sx^i$ and $\vect\sx\m^i$ as \zzz
above, with $A_0$ any finite set containing the constants and $\infty$ and $\vect{A_0}$ a listing of $A_0$. Let \zzz
$i,j$ range over $\Omega$.

    For each $f$ in \sss F, let $f^0$ be the constant function of appropriate \zzz
arity with value $\infty$. Suppose we have defined $\sss F^i$. If $\sss F^i$ \zzz
satisfies the premises (i) to (iii), stop. Otherwise, choose some $\ssf Qx\m\tx$ \zzz
in $\ssf{sub}(H)$ whose corresponding condition is false, and choose $\vect u,x$ in \zzz
$A^i_{n-1}$ as witnesses. Define $\sss F^{i+1}$ by:
\begin{mylist}
\item if $f\e\sss F$, $f\ne f_{\ssf Qx\tx}$ and $\rank(f)\le r$, let $f^{i+1}=f^i$;
\item if $f\e\sss F$ and $\rank(f)>r$, let $f^{i+1}=f^0$; and
\item let $f^{i+1}_{\ssf Qx\tx}$ equal $f^i_{\ssf Qx\tx}$ everywhere except at $\vect u$, where we let $f^{i+1}_{\ssf Qx\tx}(\xvect u)=x$.
\end{mylist}
We shall call $(f_{\ssf Qx\tx},\xvect u)$ the $(i+1)^\text{st}$-\ul{critical pair}. We shall suppose that \zzz
the choice of the $(i+1)^\text{st}$-critical pair depends only on $\vect\sx\m^i$. The \ul{rank} \zzz
of $\sx^{i+1}$ will be $r$; let $\sx^0$ have rank $0$.

    We must show that our algorithm halts.

    First note that if $(f,\vect u)$ is the $i+1^\text{st}$-critical pair, then $f^{i+1}\xvect u\prec_+f^i\xvect u$. \zzz
If $f^i\xvect u=\infty$ this is trivial, so suppose $f^i\xvect u\ne\infty$. Then there exists a \zzz
greatest $j\le i$ such that $(f,\vect u)$ is the $j^\text{th}$ critical pair. Let $f = f_{\ssf Qx\tx}$, \zzz
and let $\pm\tx$ be $\tx$ or $\nx\tx$ according to whether \ssf Q is \Ex\ or \Ax. Then \zzz
$\sss S^i(\pm\tx)(f^j\xvect u,\vect u)$ is true, and for all $k$ with $j\le k\le i$, $\sss S^k(\pm\tx)(f^k\xvect u,\vect u)$ \zzz
remains true since nothing has altered. Thus, for $(f,\vect u)$ to be the $i+1^\text{st}$ \zzz
critical pair, the rank of $f$ must be $q$ and $f^{i+1}\xvect u\prec f^i\xvect u$. Finally, \zzz
$f^{i+1}\xvect u\ne\infty$, for otherwise there would exist a loop.

    Next note that if $i<j$, if the rank of $\sx^j$ is not less than that \zzz
of $\sx^i$, and if the rank of $\sx^k$ is greater than that of $\sx^j$ for all $k$ \zzz
between $i$ and $j$, then for all $f$ in \sss F, and all $\vect u$, $f^j\xvect u=f^i\xvect u$ except \zzz
when $(f,\vect u)$ is the $j^\text{th}$ critical pair, in which case $f^j\xvect u\prec_+f^i\xvect u$.

    Finally note that if for all $f$ in \sss F and all $\vect u$, $f^j\xvect u\!\prec_+\!f^i\xvect u$ or \zzz
$f^j\xvect u=f^i\xvect u$, then $\vect\sx\m^j=\vect\sx\m^i$ or $\vect\sx\m^j_n\ll\vect\sx\m\m^i_n$, where $\ll$ is the ordering $(\prec_+)\m^\w$.

    Let $\seq{i,\,\nu}$ denote $\seq{\vect\sx\m^j}_{i\le j<\nu}$, where $\nu\le\Omega$. $\seq{i,\,\nu}$ is \zzz
an $r$-\ul{subroutine} if the rank of $\sx^i$ is $\le r$, the rank of $\sx^j$ is greater \zzz
than $r$ for all $j$ between $i$ and $\nu$, and if $\nu<\Omega$, the rank of $\sx^\nu$ is
$\le r$.\footnote{I have reversed Scanlon's terminology here because my labelling
seems to be more natural if one hopes to extend this method to stronger
systems.}

\thm{3.4 Lemma} Let $\seq{i,\,i+s}$, $\seq{j,\,j+s}$ be two consecutive (i.e.\ $i+s\eqx j$) \zzz
$r$-subroutines and suppose that the rank of $\sx^j$ is not less than that of \zzz
$\sx^i$. Then there exists $p<\min(s,t)$ such that for all $l<p$, $\vect\sx\m^{i+l}=\vect\sx\m^{j+l}$, \zzz
and $\vect\sx\m^{j+p}_n\ll\vect\sx\m^{i+p}_n$.

\noindent
Proof: By the above remarks, $\vect\sx\m^j_n\ll\vect\sx\m^i_n$ or $\vect\sx\m^j=\vect\sx\m^i$. In the first case, \zzz
let $p\eqx 0$ and we are done. If $\vect\sx\m^j=\vect\sx\m^i$,
then $\vect\sx\m^{j+1}_n\ll\vect\sx\m^{i+1}_n$ or $\vect\sx\m^{i+1}=\vect\sx\m^{j+1}$.
In the first case we let $p\eqx 1$, and in the second we repeat our argument \zzz
again and again. We do in fact find a $p<\min(s,t)$, for suppose \zzz
otherwise. If $s<t\;\ (t<s)$, then $\vect\sx\m^{i+s-1}=\vect\sx\m^{j+s-1}\;\ (\vect\sx\m^{i+t-1}=\vect\sx\m^{j+t-1})$ implies \zzz
that the rank of $\sx^j\;\ (\sx^{j+t})$ is greater than $r$. And if $s \eqx  t$, then \zzz
$\vect\sx\m^{i+s-1}\eqx \vect\sx\m^{j+t-1}$ implies the existence of a common $j^\text{th}$ and $(j + t)^\text{th}$ \zzz
critical pair, contradicting the remarks above.\done

    For each $r\le q$ define an ordering $\ll_r$ on the \ul{finite} $r$-subroutines \zzz
inductively as follows. $q$\nobreakdash-Subroutines are always of length $1$, and we \zzz
say $\seq{\vect\sx\m^i}\ll_0\seq{\vect\sx\m^j}$ if $\vect\sx\m^i_n\ll\vect\sx\m^j_n$. For all $r<q$, each $r$-subroutine is \zzz
the concatenation of a unique sequence of $r+1$-subroutines: let $\ll_r$ \zzz
be the lexicographic ordering based on $\ll_{r+1}$. The next lemma is purely \zzz
combinatorial.

\thm{3.5 Lemma} Let $\seq{i,\,i+s}$, $\seq{j,\,j+t}$ be two $r$-subroutines such that \zzz
for some $p<\min(s,t)$, $\vect\sx\m^{j+p}_n\ll\vect\sx\m^{i+p}_n$ and $\vect\sx\m^{j+l}=\vect\sx\m^{i+l}$ for all $l<p$. \zzz
Then $\seq{j,\,j+t}\ll_r\seq{i,\,i+s}$.

\noindent
Proof: By induction on $r$. If $r\eqx q$ there is nothing to prove. Suppose \zzz
the lemma holds for $r+1$, and let $\seq{i,\,i+s}$, $\seq{j,\,j+t}$ satisfy the \zzz
hypotheses. Choose the greatest $k\le p$ such that $\sx^{i+k}$ is of rank $\le r+1$. \zzz
(Such a $k$ exists because $\sx^i$ has rank $r$.) Then $\sx^{j+k}$ has the same \zzz
rank. By the inductive hypothesis, the $r+1$-subroutine beginning with
$\sx^{j+k}$ is less than (in the order $\ll_{r+1}$) the one beginning with $\sx^{i+k}$. Since \zzz
for all $l<k$, $\vect\sx\m^{i+l}$ and $\vect\sx\m^{j+l}$ are equal and have equal ranks, we may \zzz
conclude that $\seq{j,\,j+t}\ll_r\seq{i,\,i+s}$.\done

    By induction on $r<q$ we see that each $r$-subroutine is a concatenation \zzz
of a $\ll_{r+1}$-descending sequence of finite $r+1$-subroutines, and so \zzz
each $r$-subroutine is finite. In particular, the (unique) $0$-subroutine \zzz
is finite, i.e.\ our algorithm terminates.

    For any formula \px\ (perhaps with parameters) the above algorithm \zzz
is easily modified to produce a \sx\ which also fulfils \px\ if \px\ is \zzz
true. If \px\ has definable satisfaction functions, we merely add these \zzz
to our language. If \px\ does not have (globally) definable satisfaction \zzz
functions, we may still define some finite approximations as we go along. \zzz
More precisely, to each partial computation, $\seq{\vect\sx\m^0,\vect\sx\m^1,\dots,\vect\sx\m^k}$ we assign \zzz
an \ul{index} as follows. For each $r\le q$ let $\infty_r$ be a new constant symbol \zzz
and extend the ordering $\ll_r$ to the ordering $\ll^+_r$ which has $\infty_r$ as its \zzz
maximal element. Let the index of $\seq{\vect\sx\m^0,\dots,\vect\sx\m^k}$ be $\seq{a_0\m,a_1,\dots,a_r}$, \zzz
where $r$ is the rank of $\sx^k$ for all $s\le r$,
\[
a_s=
\begin{cases}
\seq{j,\,l},\text{ where $l\le k$ and $j$ is the greatest $j<k$ for which this is an $s$-subroutine, if this exists,}\\
\infty_s, \text{ otherwise.}
\end{cases}
\]
Order these indices lexicographically, that is, $\seq{a_0\m,\dots,a_s}<\seq{b_0\m,\dots,b_t}$ \zzz
iff \ul{either} there exists $r<\min(s,t)$ such that $a_p\eqx b_p$ for all $p<r$ \zzz
and $a_r\ll^+_r b_r$ or $s<t$ and $a_r\eqx b_r$ for all $r\le s$.

    Now consider the set
\[
\myset{\aaa}{\text{\aaa\ is the index of some partial computation $\seq{\sx_0\m,\dots,\sx^k}$ which
is s.t.\ each $\sx^k$ fulfils $\Tr\Sx_l$}}\,.
\]
By $\Foundation(\prece_\ee)$ choose a minimal \aaa\ in this set and some partial \zzz
computation $\seq{\sx_0\m,\dots,\sx^k}$ which is a witness for this \aaa. Then $\sx^k$ \zzz
is in the halting state, for suppose not. As $\Foundation(\prece_\ee)$ implies \zzz
\Induction, we may choose a further finite approximation to the satisfaction \zzz
functions for $\Tr\Sx_l$ and then using these perform the next step of our algorithm. \zzz
This produces a $\sx_{k+1}$ which also fulfils $\Tr\Sx_l$. But $\seq{\sx_0\m,\dots,\sx^{k+1}}$ has index \zzz
less than \aaa, a contradiction.

    In particular, we can ensure that \sx\ fulfils the axioms of \s, and \zzz
similarly we can have $(\Tr x\nprec y)^\sx$. This completes the proof of 3.1.\done

    \refstepcounter{thisisdumb}\label{Chapter3:3.6}Next we fulfil a promise made in Chapter~II on page \pageref{Chapter2:3.6}.

\thm{3.6 Corollary} $\s\provex\Ax \px\,\ppl \Ax m\,\Ex\tau\,(\lenx\tau m\ax\px^\tau)\imp\Ax n,\Gamma\,\Ex\sx\,\big(\lenn\ax\px^\sx\ax(\Gamma\vx\nx\Gamma)^\sx\big)\m\ppr $, \zzz
and \PRA\ proves the $\ast$fulfilment version.

\noindent
Proof: We shall sketch an informal proof. Let \px, $\Gamma$, and \n\ be given, \zzz
and suppose that Skolem function symbols for the existential quantifiers \zzz
of \px\ are added to the language. Consider the proof of 3.1 when $\prece$ is \zzz
the empty relation and $\Gamma\eqx H$. Then it is easy to calculate an upper \zzz
bound on the length of any $\ll_i$-descending sequence. Hence we can \zzz
calculate an upper bound $b$ on the number of iterations of our algorithm \zzz
requires to find a sequence $n$-fulfilling $\Gamma\vx\nx\Gamma$ and which is closed \zzz
under the functions of the language---and the bound is independent of how \zzz
these functions are interpreted. Now choose $m=(b+1)(n+1)$, choose a \zzz
$\tau$ which $m$-fulfils \px, and interpret the Skolem functions for the existential \zzz
quantifiers of \px\ in accordance with the winning strategy for \zzz
the game associated with $\px^\tau$. We have chosen $m$ large enough so that \zzz
we may carry out our computation without encountering any undefined \zzz
values of these Skolem functions.\done

    If $\s+\Foundation(\prece)=\PA$, our proof gives the sharp bounds of Minc\,\cite{Minc71}.

\thm{3.7 Corollary} For any natural ordering $\prec$ of order type $2^{\w^2}_{k+1}$, \zzz
the consistency of \s\ + $\Sx_{k+1}$-\Induction\ is provable in $\Foundation(2^{\w^2}_{k+1})$ with \zzz
$\Foundation(2^{\w^2}_{k+1})$ restricted to primitive recursive predicates, where \zzz
$2^\aaa_0=\aaa$, $2^\aaa_{k+1}=2^{2^\aaa_k}$.

\noindent
Sketch of Proof: Note that $\PA^-$ implies
\[
\Sx_k\HY\Induction\equiv\Pi_k\HY\Induction\equiv\Sx_k\HY\Foundation(\w)\,.\label{Chapter3:eq:3}\tag{3}
\]
Also note that that $\w^2$ is essentially just the lexicographic ordering of \zzz
pairs of integers and so for any ``natural'' ordering of order type $\w^2$ \zzz
we have
\[
\Pi_k\HY\Foundation(\w^2)\imp\Sx_{k+1}\HY\Foundation(\w)\,.\label{Chapter3:eq:4}\tag{4}
\]
Inspection of the ordinals in the proof of 3.1 yields that
\[
\ssf{PR}\HY\Foundation\big(2^{(\w^2+1)^\w}_k\m\big)\,\provex\,\Ax n\,\Ex\sx\,\big(\lenn\ax(\Pi_k\HY\Foundation(\w^2))^\sx\ax(\Tr\Dx_0)^\sx\m\big)\,,
\]
where we obtain ``$\Tr\Dx_0\m^\sx$'' by considering $\Dx_0$ matrices rather than \zzz
quantifier-free ones, and adding the appropriate Skolem functions to \zzz
the language. This in turn implies, by 2.6 and the formalized versions \zzz
of \eqref{Chapter3:eq:3} and \eqref{Chapter3:eq:4},
\[
\ssf{PR}\HY\Foundation\big(2^{(\w^2+1)^\w}_k\m\big)\,\provex\,\RFN_{\Sx_1}(\PA^-+\Sx_{k+1}\HY\Induction)\,.
\]
It remains to note:
\[
(\w^2+1)^\w\le(\w^3)^\w=\w^\w=2^{\w^2}.\tag*{\done}
\]

    Our next corollary is the ``no-counter-example'' interpretation due \zzz
to Kreisel. Consider for example the sentence $\px=\Ex x\,\Ax y\,\Ex u\,\Ax v\,\yx xyuv$, \zzz
where \yx\ is $\Dx_0$. Third-order logic easily yields:
\begin{align*}
\px&\equiv\nx\Ax x\,\Ex y\,\Ax u\,\Ex v\m\nx\yx \\
   &\equiv\nx\Ex f,g\,\Ax x,u\,\nx\yx(x,f\mx\mx x,u,gxu) \\
   &\equiv\Ax f,g\,\Ex x,u\,\yx(x,f\mx\mx x,u,gxu) \\
   &\equiv\Ex H,K\,\Ax f,g\,\yx(Hfg,f(Hfg),Kfg,g(Hfg)(Kfg))\,.
\end{align*}
If \px\ is true in \bb N, we may choose $\seq{H,G}$ to be \ul{recursive}: let $\seq{Hfg,Kfg}$ \zzz
be the least pair $\seq{x,u}$ such that $\yx(x,fx,u,gxu)$.

\thm{3.8 Corollary} If $\PA\provex\px$, we may choose $\seq{H,K}$ above to be $<\!\ee_0$-recursive.

\noindent
Sketch of proof: We shall not define the notion of $<\!\ee_0$-recursiveness, but \zzz
merely indicate how to obtain $H$ and $K$. Suppose $\PA^-+\Sx_k\HY\Induction\provex\px$. \zzz
By the proofs of 3.1 and 2.5, by $2^{\w^2}_k$-recursion we may find, uniformly \zzz
in $f$ and $g$, a sequence \sx\ which $5$-fulfils \px\ and which is closed \zzz
under $f$ and $g$. Let $\seq{Hfg,Kfg}$ be the least pair $\seq{x,u}$ from $\sx_1\times\sx_3$ \zzz
such that $\px(x,fx,u,gxu)$.\done

    We shall end this chapter with an application which will be useful \zzz
in the next. First we give a general corollary. Let \s\ consist of a \zzz
finite number of $\Sx_3$ sentences.

\thm{3.9 Corollary} (i)\, Let $k\ge 1$ and suppose for some $\Sx_{k+2}$ formula $\Xx(x)$
\[
\s\provex\Ex!x\m\Xx(x)\ax\text{ ``$x$ is a binary relation''}\,.
\]
Denote this relation by $\prec$. Then the schema of $\Foundation(\prec)$ is a \zzz
$\Pi_{k+2}$-conservative extension of $\nabla_k\HY\Foundation(\prec)$ over \s, where
\[
\nabla_k = \text{closure of $\Sx_k$ under Boolean operations and bounded quantification.}
\]
(ii)\, More generally we have the following. Consider the schema of \zzz
so-called Bar-Induction:
\begin{flalign*}
&\Gamma\HY\BI_x\!:&               &\wf(x)\imp\yx x &&&& \\
&\Gamma\HY\BI\!:  & \quad \Ax x\,(&\wf(x)\imp\yx x)
\end{flalign*}
where $\yx\mx x$ ranges over $\breve\Gamma\HY\Foundation(x)$. \zzz
Let $\tx xm$ be $\Pi_{k+2}$ with $k\ge 1$. Then for any instance $\px x$ of $\BI_x$, \zzz
there is a instance $\px_0 x$ of $\nabla_k\HY\BI_x$ such that
\[
\s\provex\Ax m\ppl \Prx{\s}(\godel{\,\Ax x\,(\px x\imp\tx x\dot m)})\imp\Ax x\,(\px_0 x\imp\tx x\dot m)\m\ppr \,.
\]
Sketch of Proof: Note that (ii) implies (i), so let us consider (ii). \zzz
Choose an instance $\px_0 x$ of $\nabla_k\HY\BI_x$ so that
\[
\s\provex\Ax x\,\ppl \px_0 x\imp\px_1(2^{(x+1)^\w}_{\myolscript k})\m\ppr
\]
where $k\e\w$ and $\px_1 y\e\Pi_k\HY\BI_y$ are such that
\[
\s\provex\Ax x\,\ppl \px_1(2^{(x+1)^\w}_{\myolscript k})\imp\Ax y\,\big(\nx\tx y\imp\Ax n\,\Ex\sx\,(\lenn\ax x\e\sx_0\ax(\nx\tx y)^\sx\ax(\px x)^\sx)\m\big)\ppr \,.\tag*{\done}
\]

    \refstepcounter{thisisdumb}\label{Chapter3:FactDeferredProof}
In Chapter~VI page \pageref{Chapter6:6.10} we shall, however, prove a result (essentially due to \zzz
Friedman\,\cite{Frie76}) which in most applications is stronger, as follows. Let \zzz
\ul{countable} bar-induction be the schema:
\begin{flalign*}
&\cBI\!:&& \quad\quad\Ax x\,(\text{``$x$ countable''}\ax\wf(x)\imp\yx\mx x)\,,\text{ where $\yx\mx x$ ranges over \Foundation($x$).}&&
\end{flalign*}

\thm{Fact}: Let $k>1$ and let $\tx xm$ be $\Pi_{k+2}$. Then for any instance $\px x$ of $\cBI_x$, \zzz
there is an instance $\px_0 x$ of $\Pi_k$-\cBI\ such that
\[
\s+\Dx_0\HY\cBI\,\provex\,\Ax m\ppl \Prx{\s+\text{\w-rule}}(\godel{\,\Ax x\,(\px x\imp\tx x\dot m)})\imp\Ax x\,(\px_0 x\imp\tx x m)\m\ppr \,.\tag*{\done}
\]
We shall also mention, for the purposes of comparison, the following \zzz
well-known result.

\thm{Fact}: If $k>1$ and $\tx m$ is $\Pi_{k+2}$ then
\[
\s+\Sx_k\HY\DC\provex\Ax m\ppl \Prx{\s+\BI+\text{\w-rule}}(\godel{\tx\dot m})\imp\tx m\ppr \,.\tag*{\done}
\]

    The application we have in mind uses, however, only the weakest of \zzz
these conservation results: 3.9.i with $\prec$ the standard ordering on \w. \zzz
(Indeed, this is probably the most useful, for when one is formalizing \zzz
an informal proof, it is usually clear what instances of strong axioms \zzz
are used, but the instances of \Induction\ are often not so obvious.) It \zzz
is fairly straightforward to check that the usual proofs of the \zzz
Kondo-Addison Uniformization Theorem may be formalized to give:

for each $\Pi^1_1$ formula $\tx XY$ there exists a $\Pi^1_1$ formula $\yx XY$ such \zzz
that
\[
\Pi^1_1\HY\CA\provex\Ax X\,\ppl \Ex Y\m\tx XY\imp\Ex!Y\,(\tx XY\ax\yx XY)\m\ppr \,.
\]
We note that the right-hand-side is $\Pi^1_4$ and that the schema of $\Pi^1_1\HY\CA$ \zzz
is (equivalent in $\Pi^1_1\HY\CA\up$ to a schema which is) $\Pi^1_3$. Since $\Pi^1_2\HY\CA\up$ implies \zzz
$\nabla^1_2$-\Induction, we have that
\[
\Pi^1_2\HY\CA\up\,\provex``\Pi^1_1\ssf{-Uniformization}"
\]
and so by the usual arguments,
\[
\Pi^1_2\HY\CA\up\,\provex\Sx^1_2\HY\AC\up\,\,.
\]

    A similar application is as follows. If $k\ge3$, then for any $\Sx_k$ \zzz
formula $\tx XY$,
\[
\Pi^1_{k-1}\HY\CA\provex\Ax X,W\,\ppl V\eqx L(W)\imp\big(\Ex Y\m\tx XY\imp\Ex Y\,(\tx XY\ax\Ax Z\,(\tx XZ\imp Y\prece Z))\big)\m\ppr \,.
\]
where $\prece$ is a $\Pi^1_2$ formula with parameter $W$ representing Addison's \zzz
well-ordering of the reals constructible from $W$. (The $\Sx^1_k$-\BI\ needed for \zzz
this proof is derivable from $\Pi^1_k$-\CA.) Now the right-hand-side is $\Pi^1_{k+2}$ and \zzz
$\Pi^1_{k-1}\HY\CA\up$ is (essentially) $\Pi^1_{k+1}$, and so by 3.9 we have that this is provable \zzz
in $\Pi^1_k\HY\CA\up$ and so:
\[
\Pi^1_k\HY\CA\up\,+\Ex W.\,V\eqx L(W)\provex\Sx^1_2\HY\AC\up\,\,.
\]

    Our final application is a conservation result for arithmetic; let \zzz
$\s=\PA^-_\text{ex}$, as on page \pageref{Chapter1:DefinitionPAex}. Add a new constant $c$ to the language of arithmetic, \zzz
and consider the schema $\Induction(c)$ of induction up to $c$. From \zzz
the proof of 3.1 we see that for each $k,n\e\w$ and each finite subset $H$ \zzz
of \s\ + $\Induction(c)$ there exists $l\e\w$ such that
\[
\s+\Induction(2^c_{\myolscript l})\provex\Ax m\,\Ex\sx\,(\lenm\sx\mx\eqx\ol n\m\ax m\e\sx_0\ax H^\sx\ax\Tr\Sx_k\m^\sx)\,.
\]
Hence we may obtain results of the following form, where \ul{conservative extension} is as defined on
page \pageref{Chapter4:DefinitionConservative}.

\thm{3.10 Corollary} With respect to the language of arithmetic plus $c$, \zzz
$\Induction(c)$ is:
\begin{mylist}
\itemsep0ex
\item $\Pi_2$-conservative over \s\
\item $\Pi_3$-conservative over \s\ + $\Dx_1\HY\Induction(c)$
\item $\Pi_{k+2}$-conservative over \s\ + $\Sx_{k-1}$-\Collection\ + $\Sx_k\HY\Induction(2^c_{\myolscript l})_{\;l\,\e\,\w}$, and is
\item $\Pi_{k+2}$-conservative over \s\ + $\Sx_{k+2}\HY\Induction(2^c_{\myolscript l})_{\;l\,\e\,\w}$.\done
\end{mylist}
                       \cleardoublepage
\thispagestyle{plain} \sectioncentred{Some Conservation Results}

    This chapter consists, in part, of an exposition of Friedman's\,\cite{Frie70} \zzz
method of obtaining conservation results of the form e.g.\ $\Sx^1_{k+1}\HY\AC\up$ is \zzz
$\Pi^1_l$-conservative over $\Pi^1_k\HY\AC\up$, for various $k$ and $l$, and $\Sx^1_{k+1}\HY\AC$ is $\Pi^1_l$-conservative over \zzz
$\myset{(\Pi^1_l\HY\CA)_\aaa}{\aaa<\ee_0}$ (i.e., axioms asserting that the $\Pi^1_l$-jump may be iterated \zzz
\aaa\ times for all $\aaa<\ee_0$). The former type of result---the so-called \zzz
restricted case---is not, however, explicitly stated in Friedman\,\cite{Frie70}. This \zzz
is probably because at that time nobody was interested in this case, for it \zzz
seems to me improbable that one could discover the unrestricted versions \zzz
without being aware of the restricted. Let us briefly, then, indicate the \zzz
work done by others.

    Barwise and Schlipf\,\cite{BarwSchl75} gave a model theoretic proof for the $\Sx^1_1\HY\AC\up$ \zzz
case, which was considerably simplified by Feferman\,\cite{Fefe76} and, independently, \zzz
Stavi. Feferman\,\cite{Fefe76} also considered $\Sx^1_2\HY\AC\up$. These latter proofs are, \zzz
in essence, the same as those given below. Proof-theoretic proofs have been \zzz
given by Tait (\cite{Tait68} \& \cite{Tait70}) for $\Sx^1_1\HY\AC$ and $\Sx^1_2\HY\AC$, and by Feferman\,\cite{Fefe77} \zzz
and Feferman and Sieg\,\cite{Fefe80} for both $\Sx^1_k\HY\AC\up$ and $\Sx^1_k\HY\AC$ for all $k$.

    The other part of this chapter deals with \ul{uniform} versions of the above: \zzz
for example
\begin{flalign*}
&&                     \Ax n           \,(\Pi^1_1\HY\AC)_n   \up\,&\provex\RFN_{\Pi^1_3}(\Sx^1_2\HY\AC\up\,) &&\\
&\text{and}\hidewidth& \Ax\aaa\ltx\ee_0\,(\Pi^1_1\HY\AC)_\aaa\up\,&\provex\RFN_{\Pi^1_3}(\Sx^1_2\HY\AC)\,,
\end{flalign*}
where e.g.\ $\Ax\aaa\ltx\ee_0\,(\Pi^1_1\HY\AC)_\aaa\up$
consists of the \ul{axiom} asserting that for all $\aaa\ltx\ee_0$ and for all \zzz
$X$, there exists a hyperjump hierarchy of length \aaa\ relativized to $X$ (plus \zzz
other simple axioms).

    Let me briefly discuss the papers of Feferman\,\cite{Fefe77} and Feferman and \zzz
Sieg\,\cite{Fefe80}. The methods in these papers---namely, the normalization of  \zzz
infinite terms followed by a G\"odel-style functional interpretation, and  \zzz
cut-elimination arguments, respectively---are stronger than those used here \zzz
in that they give conservation results for certain extensions of the theory \zzz
of types, $\ssf Z^\infty+\ssf{QF-}\AC$, whereas mine do not. But if one is only interested \zzz
in, say, second-order conservation results, then these proof-theoretic \zzz
arguments are, to my mind, much more complicated than the corresponding \zzz
model-theoretic arguments. It seems likely, however that the uniform  \zzz
versions may also be obtained from Feferman and Sieg\,\cite{Fefe80}, and in this case \zzz
the two methods would be of about the same complexity.

    The plan of this chapter is as follows. After some preliminary \zzz
definitions, we shall consider the restricted case, then give some applications, \zzz
followed by the unrestricted case. It is hoped that the presentation here \zzz
of Friedman's construction in a more general setting (thus making it \zzz
necessary to state explicitly the principles used) helps rather than hinders \zzz
the reader.

    \refstepcounter{thisisdumb}\label{Chapter4:DefinitionConservative}Given a class of
sentences $\Lambda$ and sets of sentences $A$, $B$, and $C$, say \zzz
that $A$ is a $\Lambda$-\ul{conservative} extension of $B$ \ul{over} $C$ if $A+C$ extends \zzz
$B+C$ but for all \tx\ in $\Lambda$, if $A+C\provex\tx$ then $B+C\provex\tx$.

    Let \s\ be, as usual, a theory suitable for elementary set theory as on page \pageref{Chapter1:DefinitionMembership}. \zzz
For convenience, we shall assume that the axioms of \s\ are $\Sx_3$. In \zzz
Chapter~I we assumed, also mostly for convenience, that \s\ proves the \zzz
existence of a $\Sx_1$ satisfaction predicate for $\Dx_0$ formulae. There is, \zzz
however, an interesting application of Theorem 4.2 (Application (vii) below) in which \zzz
this is not the case, and so we note that this assumption is not necessary \zzz
for this theorem. Consider the schema \Bounding\ given on page \pageref{Chapter1:DefinitionInduction}. It is easy to show:

\thm{4.1 Fact} In the presence of $\s+\Dx_0$-\Separation, the following are equivalent:
\begin{flalign*}
&&\Pi_k\HY&\Bounding &&\\
&&\Sx_{k+1}\HY&\Bounding \\
&&\Sx_{k+1}\HY&\Collection + \Sx_{k+1}\HY\Separation.\tag*{\done}
\end{flalign*}

    Let $\s_1$ be the set of formulae of the form
\[
\Ev x\,\Ax a\,\Ex b\,\Av y\ve a\,\Ev z\ve b\,\tx\,,
\]
where $\Ev x\,\Av y\,\Ev z\m\tx$ is an axiom of \s. Obviously,
\[
\s + \Sx_1\HY\Collection\provex\s_1\,;
\]
indeed, in many natural examples, $\s\provex\s_1$.

\thm{4.2 Theorem} $\Sx_1$-\Collection\ is $\Pi_2$-conservative over $\s_1$, and for all $k\ge 1$, \zzz
$\Sx_{k+1}$-\Collection\ is a $\Pi_{k+2}$-conservative extension of $\Sx_k$-\Bounding\ over \s.

    We shall give two proofs which are essentially the same but the \zzz
second is more readily formalized. After a corollary of this second \zzz
proof, we shall give a long list of applications.

\noindent
Proof: Fix $k\ge 0$. Let $\A=\seq{A,\dots}$ be a model of $\s_1$ and, if $k\geX 1$, \zzz
$\Sx_k$-\Bounding, and suppose that \s\ is coded in \A\ and that \A\ is recursively \zzz
saturated. It suffices to show that for any $c\e A$ there exists a \zzz
$k$-elementary substructure $\B=\seq{B,\dots}$ of \A\ containing $c$ which is \zzz
a model of $\s+\Sx_{k+1}$-\Collection.

    Fix $c\e A$. For each $n\e\w$ there exists in \A\ a sequence \zzz
$\sx=\seq{\sx_0\m,\dots,\sx_{\lenm\sx}}$ such that:
\begin{mylist}
\item $\lenm\sx\ge n$;
\item $c\e\sx_0$, and the (interpretations of) the constants of the language are contained in $\sx_0$;
\item $\sx_i\e\sx_{i+1}$, $\sx_i\subseteq\sx_{i+1}$, and $x\e y\e\sx_i\!\imp\! x\e\sx_{i+1}$, for all $i<\lenm\sx$;
\item if $\Ev x\,\Av y\,\Ev z\m\tx$ is among the first \n\ axioms of \s, then there
      exists $\vect x\ve\sx_1$ such that for all $i<\lenm\sx$,
\[
\Av y\ve\sx_i\,\Ev z\ve\sx_{i+1}\m\tx\,\text{; and}
\]
\item if $k\ge 1$ and if $\Ev y\,\tx(\xvect x,\vect y)$ is among the first $n$ $\Sx_k$ formulae,
      then for all $i<\lenm\sx$,
\[
\Av x\ve\sx_i\,(\Ev y\m\tx\imp\Ev y\ve\sx_{i+1}\m\tx)\,.
\]
\end{mylist}
By recursive saturation, there exists \sx\ satisfying (i) to (v) for all $n\e\w$. \zzz
Let $B=\cupx_{i\in\m\w}\sx_i$. It is clear that $\B\modelx\s$ and that $\B\prec_k\A$. (If $k\eqx 0$, \zzz
then (iii) ensures that \B\ is a \ul{transitive} substructure, and so $\Dx_0$ \zzz
formulae are absolute.) A simple application of underspill will show that \zzz
$\Sx_{k+1}$-\Collection\ holds in \B, as follows.

    Suppose $\B\modelx\Av x\ve a\,\Ev y\m\tx$ where \tx\ is $\Pi_k$. Since \B\ is a transitive \zzz
substructure, the set $a$ has no more elements in \A\ than in \B\ and so \zzz
for any nonstandard $i\le\lenm\sx$,
\[
\A\modelx\Av x\ve a\,\Ev y\ve\sx_i\,\tx\,.
\]
By underspill, this must hold for some standard $i$. Since $\sx_i\e B$ for all \zzz
standard $i$, we may conclude
\[
\B\modelx\Ex b\,\Av x\ve a\,\Ev y\ve b\,\tx\,,
\]
and this completes the first proof of theorem 4.2.

    Our second proof is, for each fixed value of the parameter \n\ \zzz
purely internal. Fix $n\e\w$, and suppose \sx\ satisfies (i) to (v). We \zzz
claim that \sx\ fulfils the first \n\ axioms of \s, that \sx\ satisfies
\[
\Ax i\ltx\lenm\sx\,\,\Av x\ve\sx_i\,\big(\tx\xvect x\imp(\tx\xvect x)^\sx_i\big)
\]
for \tx\ among the first few $\Sx_k$ formulae (where here \ul{few} depends on \n, \zzz
on the G\"odel numbering, etc.), and that \sx\ fulfils the first few sentences \zzz
of the form
\[
\Ax a,\vect z\,\Ev x\ve a\,\Av y\,\Ex b\;(\tx\xvect x\xvect y\xvect z\imp\Av x'\!\e a\,\Ev y'\!\e b\,\tx\xvect x'\xvect y'\xvect z) \label{Chapter4:eq:1}\tag{1}
\]
where \tx\ is $\Pi_k$.

    We need only check the latter part of this claim. First observe that \zzz
for any $\Pi_k$ formula $\tx\xvect x\xvect y\xvect z$ whose G\"odel number is sufficiently small \zzz
with regard to \n, we have for all $i,j,l$, with $i\le j,\,l<\lenm\sx$, and for all $\vect x,\vect y,\vect z\ve\sx_i$,
\[
\tx^\sx_j\equiv\tx^\sx_l\equiv\nx(\nx\tx)^\sx_j \label{Chapter4:eq:2}\tag{2}
\]
because each of these holds if and only if \tx\ does. To show that \sx\ \zzz
fulfils \eqref{Chapter4:eq:1}, we use the terminology of game theory. Suppose Player I has \zzz
chosen $i<\lenm\sx$ and $a,\vect z\ve\sx_i$. Let
\[
j_0=
\begin{cases}
\min j<\lenm\sx\text{ s.t.\ }\Av x\ve a\,\Ev y\ve\sx_j\,\tx\xvect x\xvect y\xvect z \text{, if this exists,}\\
\lenm\sx \text{, otherwise.}
\end{cases}
\]
Let Player II choose $\vect x\ve a$ such that if $j_0\ne 0$, then
\[
\nx\Ev y\ve\sx_{j_0-1}\,\tx\xvect x\xvect y\xvect z\,.
\]
Now Player II has essentially won, for if Player I chooses $j<j_0$ and \zzz
$\vect y\ve\sx_i$, then
\[
\nx(\nx\tx\xvect x\xvect y\xvect z)^\sx_j
\]
is \ul{false}, and if Player I chooses $j\ge j_0$, then II picks $b$ to be $\sx_{j_0}$. This \zzz
proves the claim. The theorem now follows from 2.6.\done

    It is clear that for \ul{finite} \s\ the second proof may be formalized to \zzz
obtain by 2.6:

\thm{4.3 Corollary} $\s_1+\Sx_1\HY\Induction\provex\RFN_{\Pi_2}(\s+\Sx_1\HY\Collection)$, \zzz
and for $k\ge 1$, $\s+\Sx_{k+1}\HY\Induction+\Sx_k\HY\Bounding\provex\RFN_{\Pi_{k+2}}(\s+\Sx_{k+1}\HY\Collection)$. \done

    Further observing that we do not need to work with models of $\Sx_k$-\Bounding\ \zzz
but only with arbitrarily large ``finite approximations'' of such \zzz
models, we have, for example:
\begin{align*}
\PRA\provex\Ax k\,\Ax\tx\e\Sx_{k+2}\,\ppl
   &\Ax n\,\Ex\sx\,\big(\lenn\ax(\s+\Sx_k    \HY\Bounding  +\tx)^{\ast\sx}\big)\imp\\
   &\Ax n\,\Ex\sx\,\big(\lenn\ax(\s+\Sx_{k+1}\HY\Collection+\tx)^{\ast\sx}\big)\ppr\,,
\end{align*}
which immediately translates into:
\[
\PRA\provex\Ax k\,\Ax\tx\e\Pi_{k+2}\,\big(\Prx{\s+\Sx_{k+1}\HY\Collection}(\tx)\imp\Prx{\s+\Sx_k\HY\Bounding}(\tx)\big)\,.
\]
Similarly, we may obtain e.g.\ for $k\ge 1$,
\[
\s\provex\Ax l\gex k\,\big(\RFN_{\Pi_{k+2}}(\s+\Sx_l\HY\Bounding)\imp\RFN_{\Pi_{k+2}}(\s+\Sx_{l+1}\HY\Collection)\m\big)\,.
\]

    Next we shall consider some applications of 4.2, beginning first with \zzz
set theory. The following are immediate from our results.

\thm{Application (i)} Let $\KP\up$ be \KP\ with \Infinity\ but with the \Foundation\ \ul{schema} replaced \zzz
by the \Foundation\ \ul{axiom}. Then $\KP\up$ is a $\Pi_2$-conservative extension of \zzz
$\KP\up\setminus\Dx_0$-\Collection. Hence $\KP\up$ is a very weak set theory indeed: $L_{\w+\w}$ \zzz
is a model of its $\Pi_2$ consequences.\done

\thm{Application (ii)} For $k\ge 1$, $\Sx_{k+1}$-\Replacement\ is a $\Pi_{k+2}$-conservative extension of \zzz
$\Sx_k$-\Replacement\ + $\Sx_k$-\Separation\ over $\KP\up$. For any non-projectable $\gamma>\w$, \zzz
$L_\gamma$ is a model of the $\Pi_3$ consequences of $\Sx_2$-\Collection\ + $\KP\up$ (because $\gamma$ \zzz
is the limit of $\gamma$-stable ordinals).\done

\thm{Application (iii)} $\Sx_1\HY\Induction+\KP\up\setminus\Dx_0\HY\Collection\provex\RFN_{\Pi_2}(\KP\up)$, and for $k\ge 1$,
\begin{multline*}
\KP\up+\Sx_{k+1}\HY\Induction+\Sx_k\HY\Collection+\Sx_k\HY\Separation\\
\provex\RFN_{\Pi_{k+2}}(\KP\up+\Dx_{k+1}\HY\Induction+\Sx_{k+1}\HY\Collection+\Dx_{k+1}\HY\Separation)\,.\tag*{\done}
\end{multline*}

    Next let us briefly consider two-sorted theories with $x,y,\dots$ \zzz
intended to range over sets and $X,Y,\dots$ over classes. We obtain, for \zzz
example, that:

\thm{Application (iv)} $\GB\setminus\{\Power,\Choice\}$ plus the schema $\Sx^1_1$-\wAC\
\begin{flalign*}
&\ssf{weak-}\Sx^1_1\HY\AC\!: &&\Ax x\,\Ex Y\,\tx(x,Y)\imp\Ex Z\,\Ax x\,\Ex y\,\tx(x,(Z)_y)\,,\ \tx\text{ in }\Sx^1_1\,, &&&
\end{flalign*}
is conservative over $\ZF^-$: for $\Sx^1_1$-\wAC\ is $\Pi^1_2$-conservative over $\GB\setminus\{\Power,\Choice\}$, \zzz
and any model of $\ZF^-$ can obviously be extended to a model of $\GB\setminus\{\Power,\Choice\}$. \zzz
(A weaker but similar result was obtained by Moschovakis, c.1971.)\done

    Solovay has shown that any countable model \A\ of \ZFC\ may be extended to a model \zzz
$(\A,\prece)$ of \ZFC\ in the language with an extra binary predicate which \zzz
also satisfies
\[
\Ax x\,(\prece\cap\,x^2\text{ is a well-ordering})\,.
\]
Feferman\,\cite{Fefe76} observed that this yields:

\thm{Application (v)}(Schlipf\,\cite{Schl78}) \GB\ plus the schema
\begin{flalign*}
&\Sx^1_1\HY\AC\!: &&\Ax x\,\Ex Y\,\tx(x,Y)\imp\Ex Z\,\Ax x\,\tx(x,(Z)_x)\,,\ \tx\text{ in }\Sx^1_1\,, &&&
\end{flalign*}
is conservative over \ZFC.\done

    Now let us consider analysis; our language is that for $\seq{\w,\sss P\w;+,\times,0,1,\e}$. \zzz
Here \s\ consists of the \ul{axiom} of \Induction, $\Dx^0_0$-\CA, (where $\Dx^0_0$ is the \zzz
closure of the atomic formulae under $\ax,\vx,\nx,\Ex n\ltx m,\Ax n\ltx m$), and \zzz
simple axioms for addition and multiplication. As usual, \s\ plus a \zzz
schema such as
\begin{flalign*}
&\Gamma\HY\AC\!: &&\Ax n\,\Ex Y\,\tx(n,Y)\imp\Ex Y\,\Ax n\,\tx(n,(Y)_n)\,,\ \tx\e\Gamma\,, &&&
\end{flalign*}
will be denoted $\Gamma\HY\AC\up$, while $\Gamma\HY\AC$ will denote $\Gamma\HY\AC\up$ \ul{plus} the full \zzz
schema of \Induction. Then 4.2 gives:

\thm{Application (vi)}(Feferman and Sieg\,\cite{Fefe80}) $\Sx^1_1$-\wAC\up\ is $\Pi^1_2$-conservative over $\Dx^0_0\HY\CA\up$, over \zzz
$\ssf{PR}$-\AC\up, and also over $\Pi^0_1$-\AC\up. Thus:\done

\thm{Application (vii)}(ibid.) $\Sx^1_1$-\wAC\up\ is conservative over \PRA, and

\thm{Application (viii)}(Friedman\,\cite{Frie76}, Barwise and Schilpf\,\cite{BarwSchl75}) $\Sx^1_1$-\AC\up\ is conservative over \PA.\done

    Natural instances of \Bounding\ are just various basis theorems. For \zzz
example, the Kleene and the Kondo-Addison Basis Theorems immediately give \zzz
(upon checking that the proofs go through in the appropriate theory or by \zzz
using 4.1 and the discussion at the end of the last chapter) that:

\thm{Application (ix)}(essentially Friedman\,\cite{Frie70}, as are (x)-(xii) below) $\Sx^1_2$-\AC\up\ is a \zzz
$\Pi^1_3$-conservative extension of $\Pi^1_1$-\CA\up, and \done

\thm{Application (x)}(ibid.) $\Sx^1_3$-\AC\up\ is a $\Pi^1_4$-conservative extension of $\Pi^1_2$-\CA\up.\done

\noindent
For further basis theorems it seems that further assumptions are needed. \zzz
For example, using the relativized version of Addison's well-ordering \zzz
of the constructible reals, and that his proof may be formalized in the \zzz
appropriate theory, we have:

\thm{Application (xi)}(ibid.) For $k\ge 3$, $\Sx^1_{k+1}$-\AC\up\ is
a $\Pi^1_{k+2}$-conservative extension of $\Pi^1_k$-\CA\up\ \ul{over} \zzz
$\Ex X\m.\m V\eqx L(X)$.\done

\noindent
Since $\Ex X\m.\m V\eqx L(X)$ is $\Pi^1_4$-conservative over $\Pi^1_k$-\CA\up\ for $k\ge 3$ (and this \zzz
is easy to check), (xi) gives:

\thm{Application (xii)}(ibid.) For $k\ge 3$, $\Sx^1_{k+1}$-\AC\up\ is
a $\Pi^1_4$-conservative extension of $\Pi^1_4$-\CA\up.\done

\noindent
Uniform versions of (viii), (ix) and (x) are:

\thm{Application (xiii)} For $k=0,\,1\text{, or }2$
\[
\Pi^1_k\HY\CA\up+\Sx_{k+1}\HY\Induction\provex\RFN_{\Pi_{k+2}}(\Sx_{k+1}\HY\AC\up)\,,\ \text{(where here $\Pi^1_0=\Pi^0_1$)}
\]
and uniform versions may also be concocted for the general case.\done

    Next let us consider arithmetic. Here \s\ consists of $\PA^-_\text{ex}$,  \zzz
as given on page \pageref{Chapter1:DefinitionPAex}; all \zzz
theories below are assumed to contain \s. As is well-known, many schemata which \zzz
are very different in set theory become equivalent in arithmetic. \zzz
In particular,

\thm{Fact} In $\PA^-_\text{ex}$ for $k\ge1$,\label{Chapter4:ArithmeticCollectionComment}
\begin{gather*}
\Sx_k\HY\Induction\equiv\Pi_k\HY\Induction\equiv\Sx_k\HY\Separation\equiv\Sx_k\HY\Bounding \\
\Sx_k\HY\Induction\imp\Sx_k\HY\Collection \label{Chapter4:eq:3}\tag{3}
\end{gather*}
These results, due to Ch.\ Parsons (\cite{Pars70} and \cite{Pars72}) (and independently \zzz
rediscovered by many others), are all straightforward to prove, except for
\[
\Sx_k\HY\Induction\imp\Sx_k\HY\Separation
\]
for which we shall now give a proof using an elegant trick of H. Friedman's\,\cite{Frie79}. \zzz
Suppose we have already shown \eqref{Chapter4:eq:3}. It is easy to check that $\Sx_k$-\Collection\ \zzz
implies that any formula of the form $\Ax x\ltx y\m\tx$, with \tx\ in $\Sx_k$ is equivalent \zzz
to a $\Sx_k$ formula. Using this, it is easy to see that $\Sx_k$-\Induction\ \zzz
implies any class of the form
\[
\myset{x\ltx b}{\px x}\,,
\]
with \px\ in $\Sx_k$, has a maximum element. So given any element $a$ and any \zzz
$\Sx_k$ formula \tx, the class \myset{x\ltx a}{\tx x} is realized by the maximum \zzz
element of
\[
\myset{x\ltx 2^a}{\Ax y\e x\,\tx y}\,,\ Q.E.F.\tag*{\done}
\]

    Now these equivalences immediately give for all $k\ge 0$,

\thm{Application (xiv)}(Paris \& Kirby, Friedman) \label{Chapter4:Applicationxiv}
$\Sx_{k+1}$-\Collection\ is a $\Pi_{k+2}$-conservative extension of \zzz
$\Sx_k$-\Induction, and\done

\thm{Application (xv)}(cf.\ Parsons\,\cite{Pars71}, \cite{Pars72}) \label{Chapter4:Applicationxv}
$\Sx_{k+1}\HY\Induction\provex\RFN_{\Pi_{k+2}}(\Sx_{k+1}\HY\Collection)$.\done

    Finally let us note, as Friedman\,\cite{Frie70} remarks at the end of his paper, \zzz
that the first proof of 4.2 may be used to yield interesting results \zzz
concerning certain infinitary theories. We replace the notion of recursive \zzz
saturation by \bb A-saturation for the appropriate admissible set \bb A, and \zzz
obtain, for example:

\thm{Application (xvi)} In \w-logic, $\Sx^1_{k+1}$-\AC\ is
a $\Pi^1_l$-conservative extension of $(\Pi^1_k\HY\CA)_{<\smash{\w^{\text{CK}}_1}}$, \zzz
where $l=2,\ 3,\text{ or } 4$ according to whether $k=0,\ 1,\text{ or is } \ge 2$. (This notion \zzz
of iterated $\Pi^1_k$-\CA\ is described below; we mention this result here only to \zzz
compare it with Application (xxi).)\done

    Theorem 4.2 may be generalized to the so-called unrestricted case as \zzz
follows. We shall assume that \s\ contains $\Dx_0$-\Separation. Suppose that \zzz
for some sufficiently absolute formula $\yx(x,y)$, \s\ proves $\Ex!\,\Omega,\prec\mx(\yx(\Omega,\prec))$ and \zzz
that $\prec$ is a strict linear well-ordering of its domain $\Omega$. We shall \zzz
usually refer to $(\Omega,\prec)$ by $\Omega$, leaving the ordering implicit. As in \zzz
Chapter~III, we may define in \s\ and in a natural fashion from $\Omega$, the \zzz
linear orderings $(\Omega)^\w$, $2^{(\Omega+1)^\w}_{\myolscript k}\!\!$, and $\ee(\Omega)$. In dealing with these \zzz
orderings, it will be convenient to use the same terminology and notations \zzz
which one uses for von Neumann ordinals with $\aaa,\bbb,\dots$ ranging over $\ee(\Omega)$. \zzz
For $\sx=\seq{\sx_\bbb}_{\bbb<\aaa}$, let $\sx_{<\bbb}=\cupx_{\gamma<\bbb}\,\sx_\gamma$. Say \ul{\good}(\sx) if \zzz
\[
\Ax\bbb\e|\sx|\,\big(\cupx\sx_{<\bbb}\cup\sx_{<\bbb}\cup\{\sx_{<\bbb}\}\cup\{\sx\up\,\bbb\}\big)\subseteq\sx_\bbb\,;
\]
this is the analogue of clause (iii) in the proof of 4.2. For any \zzz
definable $\aaa\lex\ee(\Omega)$, consider the schemata:
\begin{flalign*}
&\Bounding_\aaa\!: & &\Ax b\,\Ex\sx\,\big(|\sx|=\aaa\ax b\e\sx_0\ax\good(\sx)\ax
\Ax\bbb\ltx\lenm\sx\,\Av x\ve\sx_\bbb\,(\Ev y\m\tx\imp\Ev y\ve\sx_{\bbb+1}\m\tx)\m\big)&&&\\
&\DC_\aaa\!: & \Av x\,\Ev y\,\tx\imp&\Ax b\,\Ex\sx\,\big(|\sx|=\aaa\ax b\e\sx_0\ax\good(\sx)\ax
\Ax\bbb\ltx\lenm\sx\,\Av x\ve\sx_\bbb\,\Ev y\ve\sx_{\bbb+1}\,\tx\big)\\
&\s_\aaa\!: & &\Ax b\,\Ex\sx\,\big(|\sx|=\aaa\ax b\e\sx_0\ax\good(\sx)\,\ax \\
&&&\quad\quad\aax_{\Ev x\Av y\Ev z\m\tx\;\e\;\s\cap\ol n}
\Ev x\ve\sx_1\,\Ax\bbb\lex\lenm\sx\,\Av y\ve\sx_\bbb\,\Ev z\ve\sx_{\bbb+1}\,\tx\,\big)
\end{flalign*}
Let $\Bounding_{<\ee(\Omega)}=\cupx_{n\,\e\,\w}\,\Bounding(2^{(\Omega+1)^\w}_{\myolscript n})$, and similarly for the other \zzz
schemata.

    Note that by taking the $\tx\vect x\vect y$ in $\DC_\aaa$ be
\[
\Ev z\m\tx\xvect x\xvect z\imp\tx\xvect x\xvect y\,,
\]
we have immediately:
\begin{align*}
\Sx_1    \HY\DC_\aaa&\provex\s_\aaa\,,\text{ and}\\
\Sx_{k+1}\HY\DC_\aaa&\provex\Sx_k\HY\Bounding_\aaa\,,\text{ for $k\ge 1$\,.}
\end{align*}

    Before stating our next theorem, we shall indicate, as briefly as \zzz
possible, some derivations of $\Sx_{k+1}\HY\DC_\aaa$ from more familiar principles. \zzz
For example, if $\ssf{WO}$ is the axiom asserting that every set can be well-ordered, \zzz
then
\[
\s+\ssf{WO}+(\ssf{card}(\aaa)\eqx\bbb)+\Sx_{k+1}\HY\DC_\bbb+\Sx_{k+2}\HY\Foundation(\aaa)\provex\Sx_{k+1}\HY\DC_\aaa\,.
\]
Alternatively, define a $\Sx_{k+1}\HY\ssf{Uniform\ Collection}$ schema to be one which \zzz
has the form:
\[
\Ax a,\vect z\,\big(\Av x\ve a\,\Ev y\,\tx\xvect x\xvect y\xvect z\imp\Ex!b\,(\yx_\tx\,ab\vect z\ax\Av x\ve a\,\Ev y\ve b\,\tx\xvect x\xvect y\xvect z)\m\big)
\]
where $\tx\mapsto\yx_\tx$ is a map from $\Pi_k$ formulae to $\Sx_{k+1}$ formula. (For \zzz
example, in set theory we may choose $b$ to be the least $V_\aaa$ (or perhaps \zzz
the least $L_\aaa$) such that $\Av x\ve a\,\Ev y\ve b\,\tx\xvect x\xvect y\xvect z$, and in analysis there are \zzz
the Uniformization results.) Then we have:
\[
\s+\Sx_{k+1}\HY\ssf{Uniform\ Collection}+\Sx_{k+1}\HY\Foundation(\aaa)\provex\Sx_{k+1}\HY\DC_\aaa\,.
\]

    Finally, note that in analysis the standard notion of $\Sx^1_k$-\DC\ \zzz
is equivalent to our notion $\Sx^1_k$-$\DC_\w$, which somewhat justifies our \zzz
notation. For suppose $\Sx^1_k$-\DC. This implies the schema:
\begin{flalign*}
&\Sx^1_k\ssf{-GDC}\!: &&\Ax n,X\,\Ex Y\,\tx nXY\imp\Ax U\,\Ex V\,\big(U=(V)_0\ax\Ax n\,\tx n(V)_n(V)_{n+1}\big)\,. &&&
\end{flalign*}
Suppose $\Ax X\,\Ex Y\m\tx$, where \tx\ is $\Pi^1_{k-1}$, and we wish to show
\[
\Ax A\,\Ex\sx\,(\lenm\sx=\w\ax A=\sx_0\ax\good(\sx)\ax\Ax n\,\Ax X\e\sx_n\,\Ex Y\e\sx_{n+1}\m\tx)\,. \label{Chapter4:eq:4}\tag{4}
\]
Let $\px nXY$ be:
\renewcommand{\labelenumi}{}
\vspace{-1.5ex}
\begin{enumerate}
\parskip-.8ex
\item if $n\eqx 0$, then $Y\eqx X$;
\item if $n\eqx 3(k+1)$ for some $k\e\w$, then $Y\eqx (X)_k$;
\item if $n=3(\seq{k,l}+1)+1$ for some $k,l\e\w$, then $Y=((X)_k)_l$;
\item if $n=3(k+1)+2$ for some $k\e\w$, then $\tx(X)_kY$,
\end{enumerate}
\vspace{-1.8ex}
and let $\yx UV$ be $(U=V_0)\ax\Ax n\m\px nV_nV_{n+1}$. Then $\Ax n\,\Ax X\,\Ex Y\m\yx nXY$, and so \zzz
by $\Sx^1_k\ssf{-GDC}$, $\Ax U\,\Ex V\m\yx$. Now by $\Sx^1_k$-\DC, $\Ax A\,\Ex\sx\,(\sx_0\eqx A\ax\Ax n\,\yx\sx_n\sx_{n+1})$. This \zzz
implies \eqref{Chapter4:eq:4}.

    These observations will be used when we give applications of:

\thm{4.4 Theorem} Suppose \s\ includes $\Dx_0$-\Separation. Then $\Sx_1\HY\DC_{<\ee(\Omega)}$ + \zzz
$\Foundation(\Omega)$ is a $\Pi_2$-conservative extension of $\s_{<\ee(\Omega)}$, and for $k\ge 1$, \zzz
$\Sx_{k+1}\HY\DC_{<\ee(\Omega)}$ + $\Foundation(\Omega)$ is a $\Pi_{k+2}$-conservative extension of \zzz
$\Sx_k\HY\Bounding_{<\ee(\Omega)}$ over \s.

\noindent
Proof: As before, we shall give two proofs. The first is mostly model \zzz
theoretic, but at one point it uses 3.1, namely that
\[
\s + \Foundation(\ee(\Omega))\provex\RFN(\Foundation(\Omega))\,;
\]
whereas the second proof is purely internal, but it involves details of \zzz
the Hilbert-Ackermann-Scanlon method for proving 3.1. The first proof is \zzz
essentially that of Friedman\,\cite{Frie70}, and it seems to require the further \zzz
assumption that \s\ proves $\Omega$ is infinite. The proof is in four stages.

    For the first, note that $\Foundation(\Omega)$ is $\Pi_{k+2}$-conservative over \zzz
\s\ + $\Sx_k\HY\Bounding_{<\ee(\Omega)}$ (or, if $k\eqx 0$, $\s_{<\ee(\Omega)}$). Indeed, because \s\ proves \zzz
$\Dx_0$-\Foundation, it is clear that any $\Sx_{k+2}$ sentence consistent with \s\ + \zzz
$\Sx_k\HY\Bounding_{\gamma\w}$ is consistent with \s\ + $\Sx_k\HY\Bounding_\gamma$ + $\Foundation(\Omega)$ (and \zzz
similarly for $\s_{\gamma\w}$), and that $2^{\smash{(\Omega+1)^\w}}_k\w<2^{\smash{(\Omega+1)^\w}}_{k+1}$ for large enough $k$.

    For the second stage, add new constant symbols $\aaa,\,l$ to the language, \zzz
and consider the theory:
\[
\T_\aaa=l\e\ol\w+(l>\ol n)_{n\e\w}+\aaa\eqx 2^{(\Omega+1)^\w}_l\!+\Sx_k\HY\Bounding_\aaa\text{\  (or, if $k\eqx 0$, $\s_\aaa$).}
\]
It is clear that $\T_\aaa$ is conservative over $\s=\Sx_k\HY\Bounding_{<\ee(\Omega)}+\Foundation(\Omega)$.

    For the third stage, by 3.1 we have
\[
\T_\aaa+\Foundation(\aaa)\provex\RFN(``\T_\aaa+\Sx_l\HY\Foundation(\Omega)")
\]
and so our proof of the Essential Unboundedness Theorem 2.7 gives that there \zzz
must be some instance \px\ of $\Foundation(\aaa)$ such that $\nx\px$ is $\Pi_{k+2}$-conservative over \zzz
the extensional theory represented by ``$\T_\aaa+\Sx_l\HY\Foundation(\Omega)$'', \zzz
i.e., $\T_\aaa$ + $\Foundation(\Omega)$. (Indeed, a more careful inspection yields that \zzz
we may take \px\ to be an instance of $\Sx_{k+1}\HY\Foundation(\aaa)$, but we shall not \zzz
need this fact.)

    For the last stage, the theorem now follows immediately from:

\thm{4.5 Lemma} For any instance \px\ of $\Foundation(\aaa)$, any $\Sx_{k+2}$ sentence of \zzz
\sss  L consistent with $\T_\aaa+\Foundation(\Omega)+\nx\px$ is consistent with \s\ + \zzz
$\Sx_{k+1}\HY\DC_{<\ee(\Omega)}+\Foundation(\Omega)$.

\noindent
Proof: Let \A\ be any model of $\T_\aaa+\Foundation(\Omega)+\nx\px$. Fix $m\e\w$ \zzz
sufficiently large so that $\Sx_{m+1}\HY\Foundation(\aaa)$ is false, and choose \zzz
$2^\gamma\le\aaa$ so that $\Sx_m\HY\Foundation(2^\gamma)$ holds but not $\Sx_{m+1}\HY\Foundation(2^\gamma)$. \zzz
This we can do by choosing the least (nonstandard) \n\ for which
\[
\Sx_{m+1}\HY\Foundation(2^{(\Omega+1)^\w}_n)
\]
is false, and setting $\gamma=2^{(\Omega+1)^\w}_{n-1}$; then $\Sx_m\HY\Foundation(2^\gamma)$ is implied by \zzz
$\Sx_{m+1}\HY\Foundation(\gamma)$. Now we may choose a $\Pi_{m+1}$ formula $Px$ (possibly \zzz
with parameters) such that the class
\[
I=\myset{\delta<\aaa}{P\delta}
\]
is a proper initial segment of \aaa\ for which the following overspill \zzz
property holds:
\[
\label{Chapter4:eq:5}\tag{5}
\parbox{3in}{if $\tx x$ is $\Sx_m$, possibly with parameters, and if $\tx\delta$
            holds for all $\delta\e I$ then $\tx\delta$ holds for some $\delta\e\aaa\setminus I$.}
\]

    Now let $a$ be an arbitrary element of $A$ and let \sx\ witness \zzz
$\Sx_k\HY\Bounding_\aaa$ (or, if $k\eqx 0$, $\s_\aaa$) with $a\e\sx_0$. Let $B=\cupx_{i\in\m\w}\sx_i$ and let \zzz
$\B=\A\up B$. It is clear that \B\ is a $k$-elementary transitive substructure \zzz
of \A\ and, as it is definable in \A, $\Foundation(\Omega)$ must hold in  \zzz
it. If $m\ge k+2$, it is easy to see that $\Sx_{k+1}$-\Collection\ also holds in \B\ \zzz
using the overspill property \eqref{Chapter4:eq:5} just as in the proof of 4.2. However, \zzz
we must show that $\Sx_{k+1}\HY\DC_{<\ee(\Omega)}$ holds in \B.

    Suppose $\B\modelx\Av x\,\Ev y\,\tx$, where \tx\ is $\Pi_k$, and let $b\e B$. It suffices \zzz
to show that
\[
\A\modelx\Ex\delta\e I\,\Ex!f\e\sx_\delta\,\circlex{H}\,(f,\bbb,\delta) \label{Chapter4:eq:6}\tag{6}
\]
for suitable \bbb, where
\begin{multline*}
\circlex{H}\,(f,\bbb,\delta)=\text{``$f:\bbb\mapsto\delta$ is increasing''}\ax b\e\sx_{f0}\,\ax\\
\Ax\gamma\ltx\bbb\,\big(\Av x\ve\sx_{f\gamma}\,\Ev y\ve\sx_{f(\gamma+1)}\,\tx\ax
\nx\Ex\xi\,(\gamma\ltx\xi\ltx f(\gamma+1)\ax\Av x\ve\sx_{f\gamma}\,\Ev y\ve\sx_\xi\,\tx)\big)\,.
\end{multline*}
This implies $\B\modelx\Sx_{k+1}\HY\DC_\bbb$, as follows. Choose $\delta\e I$ and $f\e\sx_\delta$ to \zzz
witness \eqref{Chapter4:eq:6}; then as \s\ includes $\Dx_0$-\Separation\ and as $\seq{\sx_\xi}_{\xi\ltx\delta}\e\sx_\delta$, \zzz
we have that $\seq{\sx_{f \xi}}_{\xi\ltx\bbb}\e\sx_{\delta+1}\subseteq B$. Now as $\circlex{H}$ is absolute between \zzz
\A\ and \B, this is our required witness to $\Sx_{k+1}\HY\DC_\bbb$.

    To prove \eqref{Chapter4:eq:6} we use induction on \bbb. We may suppose that \bbb\ is a \zzz
limit ordinal, for the successor case is the same as for $\Sx_{k+1}$-\Collection. \zzz
If \eqref{Chapter4:eq:6} is false, then by the induction hypothesis,
\[
\aaa\setminus I=\myset{\delta}{\A\modelx\Ax\gamma\ltx\bbb\,\Ex f\e\sx_\gamma\,\circlex{H}\,(f,\bbb,\delta)}\,.
\]
This gives us a $\Dx_0$ definition of $I$, for if $\lambda$ is a limit ordinal, for \zzz
all $\delta\ltx\lambda$ and $\vect x,\vect y\ve\sx_\delta$,
$\tx\xvect x\xvect y\equiv(\tx\xvect x\xvect y)^{\seq{\sx_{\delta+1},\dots,\sx_{\delta+k}}}$. But this \zzz
contradicts \eqref{Chapter4:eq:5}.

    This concludes our first proof of 4.4.\done

    Our second proof of 4.4 is purely internal, that is, it involves no \zzz
external model-theoretic considerations. We shall suppose that $k\ge 1$; the \zzz
$k\eqx 0$ case is similar and easier. Consider the following:

\thm{4.6 Lemma} Suppose \s\ contains the axiom of \Infinity. Let $l\e\w$ with \zzz
$l\ge 1$,and let $\bbb_0=\w^{2^{(\Omega+1)^\w}_l+1}$. Then
\[
\s+\Sx_k\HY\Bounding_{\bbb_0+1}\provex\Ax n,x\,\Ex\sx\,
\big(\lenn\ax x\e\sx_0\ax\Sx_l\HY\Foundation(\Omega)\m^\sx\ax\Sx_{k+1}\HY\Collection\m^\sx\ax\Tr\Sx_k\m^\sx\big)\,.
\]

    Comment before proof: If \s\ includes the negation of \Infinity, then the above ordinal is too \zzz
large for $\Sx_k\HY\Bounding_{\bbb_0+1}$ to make sense. We have, however the following \zzz
extension of 3.10. Add a new constant $c$ to the language of arithmetic.

\thm{4.7 Lemma} In arithmetic, for each $n,l\e\w$ there exists $m\e\w$, such that:
\[
\s+\Sx_k\HY\Bounding_{2^c_{\myolscript m}}\provex\Ax x\,\Ex\sx\,
\big(\lenn\ax x\e\sx_0\ax\Sx_l\HY\Induction(c)\m^\sx\ax\Sx_{k+1}\HY\Collection\m^\sx\ax\Tr\Sx_k\m^\sx\big)\,.\tag*{\done}
\]
The proof of 4.7 is similar to that of 4.6, and is omitted. After sketching an \zzz
informal proof of 4.6, we shall briefly describe how one could modify the \zzz
proof to show the existence of a sequence which fulfils $\Sx_{k+1}\HY\DC_\bbb$ for \zzz
suitable \bbb\ rather than $\Sx_{k+1}$-\Collection. From this amended version, theorem 4.4 \zzz
follows immediately. For the applications which we have in mind, however, \zzz
4.6 suffices.

\noindent
Proof of 4.6: Recall from Chapter~III that for any finite subset $\Gamma$ of \zzz
$\Foundation(\Omega)$, any $n\e\w$, and any finite sets $C,\sss G$ of constant and \zzz
function symbols, there is an algorithm $\bb A_{\Gamma,n}(C,\sss G)$ which, upon inputting any \zzz
(interpretation of $C$ as) elements $C$ and any (interpretation of the function \zzz
symbols \sss G as) functions \sss G (of the appropriate arity), produces a sequence \zzz
\sx\ which \n-fulfils $\Gamma$, has $C\subseteq\sx_0$, and is closed under \sss G. Further \zzz
recall that each state of \bb A (i.e.\ each iterate of the single loop of which \zzz
\bb A consists) depends upon only finite bits of \sss G. More specifically, the \zzz
$0^\text{th}$ (initial) state depends only upon the values of terms built up from \sss G \zzz
and $C$ (and the constants and functions of the language) whose height is \zzz
less than \n, and in general, to calculate the $i^\text{th}$ state requires, at \zzz
most, knowledge of the values of the terms whose height is less than $n(i+1)$.

    Also recall that with each state of \bb A is associated an ordinal such \zzz
that
\[
\aaa_0>\aaa_1>\aaa_2>\dots
\]
where $\aaa_i$ is the ordinal of the $i^\text{th}$ state. These $\aaa_i$ in general \zzz
depend upon \sss G, but we may obtain an upper bound for $\aaa_0$ depending only \zzz
on $\Gamma$ and \n:
\[
\aaa_0<2^{(\Omega+1)^m}_l\!,\text{ for some $m$}
\]
if $\Gamma$ is included in $\Sx_l\HY\Foundation(\Omega)$.

    Switching topics slightly, let \yx\ be any sentence and suppose some \zzz
sequence $\tau=\seq{\tau_0\m,\dots,\tau_{\lenm\tau}}$ fulfils \yx. Let \sss G be a set of satisfaction \zzz
function symbols for \yx. Then, as in 3.6, for any finite $C$ \zzz
there exists a finite interpretation of \sss G (also denoted by \sss G) by \zzz
partial functions such that:
\begin{mylist}
\item all terms of depth $\lex\lenm\tau$ built up from $C$ and \sss G have a value;
\item if $g\e\sss G$, $i\ltx\lenm\tau$, and $\vect x\ve\tau_i\cap domain(g)$, then $g\vect x\ve\tau_{i+1}$; and
\item \sss G is, in so far as it is defined, a set of satisfaction functions for \yx.
\end{mylist}
We just choose the interpretation for \sss G via a winning strategy for the \zzz
game associated with $\yx^\sx$. Also, if $\tau$ satisfies
\[
\Ax i\ltx\lenm\tau\;\Av x\ve\tau_i\,\big(\px\m\vect x\imp(\px\m\vect x)^\tau_i\big)\,,
\quad(\text{i.e. }(\Tr \px\m\vect x)^\tau)\,,
\]
we could make a similar remark for the satisfaction functions for \px.

    Now after the above preliminaries, let us proceed with the sketch of \zzz
the proof of the lemma. Fix any universal instance
\[
\Ax a,\vect z\,(\Ax x\e a\,\Ex y\m\tx\imp\Ex b\,\Ax x\e a\,\Ex y\e b\,\tx)
\]
of $\Sx_{k+1}$-Collection with \tx\ in $\Pi_k$, and let
\[
\yx=\Ax a,\vect z\,\Ex x\e a\,\Ax y\,\Ex b\,(\tx\imp\Ax x'\e a\,\Ex y'\e b\,\tx x'y')\,.
\]
Fix a finite subset $\Gamma$ of $\Sx_l\HY\Foundation(\Omega)$ and let $\px x$ be any complete \zzz
$\Pi_k$ formula. Fix \n\ and $x$. We shall show
\[
\Ex\sx\,(\lenn\ax x\e\sx_0\ax\yx^\sx\ax\Gamma^\sx\ax\Tr\px\m^\sx)\,.
\]
Let $\tau$ be a witness to $\Sx_k\HY\Bounding_{\bbb_0}$, with $x\e\tau_0$. Let $C\eqx \{x\}$ and \zzz
let \sss G consist of satisfaction function symbols for \yx\ and \px. Observe \zzz
that for each \ul{finite} subset $M$ of $\bbb_0$, the sequence $\seq{\tau_\gamma}_{\gamma\e M}$ fulfils \zzz
\yx. (And if the reader first proves this as an exercise, the remainder \zzz
of the proof should be clear.)\label{Chapter4:FiniteExercise} Furthermore, the remark \eqref{Chapter4:eq:2} on page \pageref{Chapter4:eq:2} \zzz
also holds, and so there exists some interpretation of \sss G so that (i), \zzz
(ii) and (iii) above hold for $\seq{\tau_\gamma}_{\gamma\e M}$ and \yx\ and $\px x$. The problem \zzz
is to choose $M$ and \sss G so that, with this choice, the computation \zzz
$\bb A=\bb A_{\Gamma,n}(C,\sss G)$ may be carried out, thus producing the required sequence. \zzz
But we cannot make an \ul{a priori} choice of $M$ and \sss G as we did in 3.6; \zzz
instead we shall choose them according to the algorithm \bb C incorporating \zzz
\bb A and another subroutine \bb B, which we shall now describe.

    Let \ssf{TBC} (\ul to \ul be \ul considered) be a subclass of the ordinals less than \zzz
$\bbb_0$, and suppose initially that all are \ssf{TBC}. Also initially set $M\eqx \{0\}$ and \zzz
let all the functions in \sss G have empty domains. Set $\aaa=2^{(\Omega+1)^m}_l$; \bbb\ \zzz
will always be the ordinal of the order type of \ssf{TBC}.

    Our algorithm \bb C is this: repeat subroutine \bb B \n\ times, and then \zzz
perform one step of \bb A. Set \aaa\ equal to the ordinal associated with the \zzz
present state of \bb A. If \bb A is not in its terminal state, repeat \bb B \zzz
\n\ times again, and so on.

    The subroutine \bb B is as follows. Suppose we are given \zzz
$M=\{\gamma_0\ltx\dotsb\ltx\gamma_j\}$, \sss G, \ssf{TBC} ($=\ssf{TBC}_\text{old}$), and \aaa. Now there exists a (unique) \zzz
sequence
\[
\mu_0\subseteq\dots\subseteq\mu_j
\]
of finite sets such that for all $i<j$, $\mu_i\subseteq\tau_{\gamma_i}$, and $\mu_i$ is included \zzz
in the domains of all the functions \sss G, and $\mu_{i+1}=\cupx_{g\e\sss G\cup\sss L}\;g''\mu_i$. \zzz
We extend the domains of the functions of \sss G to include $\mu_j$ and we reduce \zzz
$\ssf{TBC}_\text{old}$ to $\ssf{TBC}_\text{new}$, as follows. If $g\e\sss G$ is some satisfaction function \zzz
for \tx\ or \px, and $\vect y\ve\mu_j$, it is clear that we can choose the appropriate \zzz
$g\xvect y\ve\tau_{\delta_j+1}$, and \ssf{TBC} remains unaffected by this. The interesting cases \zzz
are the satisfaction functions for the quantifiers ``$\Ex x\e a$'' and ``$\Ex b$'' in \zzz
\yx. By the exercise just mentioned on page \pageref{Chapter4:FiniteExercise}, it should be clear that we can \zzz
define $g_{_{\Ex x\in a}}$ and $\ssf{TBC}_\text{new}$ so that:
\begin{mylist}
\item $\bbb_\text{new}$, the order type of $\ssf{TBC}_\text{new}$, is sufficiently large (we \zzz
shall explicitly state this requirement below), and if \zzz
\[
\gamma^+=\text{least ordinal in $\ssf{TBC}_\text{new}\cupx\{\bbb_0\}$ greater than $\gamma$}\,,
\]
we have:
\item $\Ax a,\vect z\ve\mu_j\,\Ax\gamma\,\ppl\gamma_j\lex\gamma\ltx\bbb_0\imp\Ax y\e\tau_\gamma\,
\big(\tx(g_{_{\Ex x\in a}}(a,\vect z),y,\vect z)\imp\Ax x'\e a\,\Ex y'\e\tau_{<\gamma^+}\,\tx\m x'y'\xvect z\big)\ppr$.
\end{mylist}
Now we can choose $g_{_{\Ex b}}$ so that for all $a,\vect z$ in $\mu_j$, if
\[
\gamma=\text{ least }\gamma\text{ s.t. }\gamma_j\lex\gamma\ltx\bbb_0\ax\Ex y\e\tau_\gamma\,\tx(g_{_{\Ex x\in a}}(a,\vect z),y,\vect z)
\]
then $g_{_{\Ex b}}(a,\vect z,y)=\tau_{<\gamma^+}$ for all $y\e\tau_{\bbb_0}$. Add $\gamma^+_j$ to $M$.

    In requirement (i), we wish to have \bbb\ sufficiently large so that \zzz
immediately after each \n-fold iteration of \bb B, we have
\[
\bbb\gex\w^\aaa\,.
\]
That it is possible to satisfy this requirement follows from the easy \zzz
combinatorial Lemma 4.8 given below. In satisfying it, we ensure that we \zzz
will always have enough \ssf{TBC} ordinals to continue extending the domains of \zzz
the partial functions interpreting \bb G until such time as the algorithm \zzz
\bb A halts.

    This completes the description of the subroutine \bb B.

\thm{4.8 Lemma} Let $\delta\ge\w^{\aaa+1}$ and let $N\e\w$. Consider the following game. \zzz
Set $\delta_0=\delta$. On the $(i+1)^\text{st}$ move, Player I chooses a partition
\[
\delta_i=\eta+\nu
\]
and Player II chooses $\mu<\eta$. Set $\delta_{i+1}=\mu+\nu$. Player II \ul{wins} if \zzz
$\delta_N\ge\w^\aaa$. We claim that Player II has a winning strategy.

\noindent
Proof: Let Player II play as follows: if $\w^\aaa j<\eta<\w^\aaa(j+1)$ for \zzz
$j\e\w$, let $\mu=\w^\aaa j$, and if $\eta\ge\w^{\aaa+1}$, set $\mu=\w^\aaa N$. It is easy \zzz
to check that this works.\done

    This concludes our sketch of the proof of Lemma 4.6. We hope our \zzz
description has been clear enough for the interested reader to fill in \zzz
the missing details.

    Theorem 4.4 for $k\ge 1$ follows from:

\thm{4.9 Lemma} For $l\ge 1$ and for suitable \bbb, we have
\[
\s+\Sx_k\HY\Bounding(\bbb^{2^{(\Omega+1)^\w}_{\myolscript l}}\!)\provex
\Ax n,x\,\Ex\sx\,(\lenn\ax x\e\sx_0\ax\Sx_l\HY\Foundation(\Omega)^\sx\ax\Sx_{k+1}\HY\DC_\bbb\m^\sx\ax\Tr\Sx_k\m^\sx).\done
\]
The proof of 4.9 is similar to that of 4.6, but as we shall not make use \zzz
of this result, we shall omit all of it except for the statement of the \zzz
analogous (equally easy) combinatorial lemma.

\thm{4.10 Lemma} Let $\delta\ge\bbb^{(\aaa+1)\w}$ and let $N\e\w$. Let $\delta_0=\delta$ and consider the \zzz
game where on the $(i+1)^\text{st}$ move Player I presents a partition
\[
\delta_i=\textstyle\sum_{\gamma<\bbb}\,\eta_\gamma+\nu
\]
and Player II chooses $\delta_{i+1}=\eta_\gamma+\nu$ for some $\gamma<\bbb$. Player II \ul{wins} \zzz
if $\delta_N\ge\bbb^{\aaa\w}$. We claim that II has a winning strategy.\done

    As before, we may use the effective proofs to obtain uniform versions:

\thm{4.11 Corollary} If \s\ includes $\Dx_0$-\Separation,
\[
\Ax\aaa\ltx\ee(\Omega)\,\s_\aaa\provex\RFN_{\Pi_2}\big(\s+\Sx_1\HY\DC_{<\ee(\Omega)}+\Foundation(\Omega)\big)\,,
\vspace{-1.5ex}
\]
and for $k\ge 1$,
\[
\s+\Ax\aaa\ltx\ee(\Omega)\,\Sx_k\HY\Bounding_\aaa\provex\RFN_{\Pi_{k+2}}\big(\s+\Sx_{k+1}\HY\DC_{<\ee(\Omega)}+\Foundation(\Omega)\big)\,.\tag*{\done}
\]
Note that in this terminology, we can restate the content of Corollary 4.3 \zzz
as:
\[
\Ax n\,\s_\n\provex\RFN_{\Pi_2}(\s+\Sx_1\HY\Collection)\,,
\vspace{-1.5ex}
\]
and for $k\ge 1$,
\[
\s+\Ax n\,\Sx_k\HY\Bounding_n\provex\RFN_{\Pi_{k+2}}(\s+\Sx_{k+1}\HY\Collection)\,.
\]

    Let us now briefly consider some applications of 4.4. Our notation \zzz
is as before.

\thm{Application (xvii)}(J\"ager\,\cite{Jager78}) In set theory we have, for example, that \KP\up\ + \zzz
\Induction\ is a $\Pi_2$-conservative extension of $\KP\up\setminus\Dx_0\HY\Collection$ + \zzz
$\myset{\Ax x\,\Ex y\,y=L_\aaa(x)}{\aaa\ltx\ee_0}$; hence $L_{\ee_0}$ is a model of the $\Pi_2$ consequences of \zzz
\KP\up\ + \Induction.\done

\thm{Application (xviii)} Also if \aaa\ is the $\aaa^\text{th}$ nonprojectible, and if $\bbb<\aaa$ is $\Sx_1$-definable, \zzz
then $L_\aaa$ is a model of the $\Pi_3$ consequences of $\KP\up+\Foundation(\bbb)+\Sx_2\HY\Collection$.\done

    In analysis we have, for example:

\thm{Application (xix)}(Friedman\,\cite{Frie70}) For $k = 0, 1, \text{or } 2$, $\Sx_{k+1}$-\AC\ is a $\Pi^1_{k+2}$-conservative extension of $(\Pi^1_k\HY\CA)_{<\ee_0}\up$;\done

\thm{Application (xx)} For $k = 0, 1, \text{or } 2$, $\Ax\aaa\ltx\ee_0\,(\Pi^1_k\HY\CA)_\aaa\up\provex\RFN_{\Pi_{k+2}}(\Sx_{k+1}\HY\AC)\,$; and\done

\thm{Application (xxi)} For $k = 0, 1, \text{or } 2$, $\Sx_{k+1}$-\AC\ plus the schema of
\[
\BI_\prec
\]
where $\prec$ ranges over primitive recursive orderings or well-orderings is \zzz
$\Pi_{k+2}$-conservative over the schema
\[
\wf(\prec)\imp(\Pi^1_k\HY\CA)_\prec\up\,,
\]
where $\prec$ has the same range.\done

\thm{Application (xxii)} In arithmetic, extend the usual language \sss  L to $\sss  L(c)$ by adding \zzz
a new constant $c$. Then from 4.7 we have that for all $k\ge 0$, $\Sx_{k+1}$-\Collection\ \zzz
+ $\Induction(c)$ is a $\Pi_{k+2}$-conservative (w.r.t.\ the language $\sss  L(c)$) extension \zzz
of $\Sx_k$-\Induction.\done

    The proofs of 4.4 also give the following result (which, however, does \zzz
not follow from our present statement of 4.4):

\thm{Application (xxiii)} $\Sx^1_1$-\wAC\ is $\Pi^1_2$-conservative over $\ssf{PR}\HY\CA\up$ + the induction \zzz
schema restricted to first order formulae with second order parameters + the \zzz
schema
\[
\Ax X\,\Ex H\,\,\text{``$H$ is the Wainer hierarchy relativized to $X$ up to \aaa''}\,,\ \aaa<\ee_0\,.\tag*{\done}
\]
                     \cleardoublepage
\thispagestyle{plain} 
\sectioncentred{Some Model-theoretic Applications: Non-\texorpdfstring{\w}{omega}-models}

    This and the following chapter contain some model-theoretic applications \zzz
of the notion of fulfilment. In this chapter we shall consider only non-\w-models. \zzz
Let us say at the outset that much of this chapter was inspired \zzz
by the work of J. Paris and L. Kirby.

    Our first theorem gives necessary and sufficient conditions for the \zzz
existence of certain initial segments which model a given \ul{coded theory}, \zzz
and \ul{indicators} for the same. By an iteration of this, we are able to give \zzz
necessary and sufficient conditions for the existence of an initial segment \zzz
which models a given \ul{complete} theory in Corollary \hyperref[Chapter5:5.5]{5.5}. We next consider \zzz
a result of Kirby, McAloon, and Murawski\,\cite{KirbyMcAlMura79}: they noticed that for any \zzz
countable model \M\ of \PRA\ and for any theory \T\ in the language of \zzz
analysis extending $\Pi^0_1$-\CA\ which is coded in \M\ there exists an indicator \zzz
for the initial segments \I\ of \M\ which are such that $\seq{I,\sss R_I(M)}$ is \zzz
a model of \T. (These terms are defined below.) We shall extend this \zzz
result by weakening the condition on \T: namely, in \hyperref[Chapter5:5.8]{5.8} we only require that \T\ \zzz
extend $\Dx^0_0\HY\CA\up$, \WKL\ (\ul{w}eak \ul{K}\"onig's \ul{l}emma) plus (something weaker than) \zzz
$\Sx^0_1$-\Induction. Our last result \hyperref[Chapter5:5.11]{5.11} gives a description of the order type of \zzz
the set of elementary initial segments of a recursively saturated (r.s.) \zzz
model of \PA, and it extends results which were known to hold for \ul{countable} \zzz
r.s.\ models (or, more generally, \ul{resplendent} models).

    By use of the arithmetized Completeness Theorem of Hilbert and Bernays, \zzz
we can obtain the following well-known result.

\thm{5.1 Result} Let \M\ be a nonstandard model of $\PA^-$ + $\Dx_2$-\Induction, and \zzz
let \T\ be any theory coded in \M\ extending \PA. There exists a r.s.\ end-extension \zzz
of \M\ which is a model of \T\ iff the $\Sx_1$ theory of \M\ is consistent \zzz
with \T.\done

    By the use of fulfilment we may easily obtain the dual of 5.1, that \zzz
is, with ``end-extension'' replaced by ``initial segment'' and ``$\Sx_1$'' by ``$\Pi_1$''. \zzz
Our first theorem is an elaboration upon this idea.

    But first we must consider the question: what \ul{is} an initial segment? \zzz
For models of arithmetic the answer is clear, and for set theory we have \zzz
the notion of transitive substructure. For analysis, however, there is no \zzz
obvious choice. If we use the notion determined by the membership \zzz
relation (as defined on page \pageref{Chapter1:DefinitionMembership}) and have $\ol\w$ as a constant of our language, \zzz
then the notion of \ul{an initial segment which is a model of $\ssf{QF}\HY\AC\up$} is the \zzz
same as that of \ul{an \w-absolute substructure which is a model of $\ssf{QF}\HY\AC\up$}, which \zzz
is indeed a very natural one. Or we could base the notion of initial \zzz
segment on some quite different ordering, such as \ul{$x$ is constructible} \zzz
\ul{before $y$}, or \ul{$x$ has hyper}-(or \ul{Turing} or \ul{Wadge}) \ul{-degree less than that of $y$}. \zzz
We shall, however, always assume that the notion of initial segment is the \zzz
one based on the membership relation, that is, $B$ is an \ul{initial segment} \zzz
of a structure \A, denoted by $B\subseteq^\text{end}\A$, if for all $x,y\e A$, $x\e y\e B$ \zzz
implies $x\e B$ (where here, of course, the relation $x\e y$ is to be interpreted \zzz
in \A). While this assumption is partially for simplicity, it is \zzz
also because of doubts about the model-theoretic importance of these other \zzz
notions. At the same time, we would like to emphasize that our constructions \zzz
are for the most part perfectly general, and work for any appropriate \zzz
interpretation of ``bounded quantifier''.

    Let \A\ be a model of \s\ as on page \pageref{Chapter1:DefinitionMembership}, let $B\subseteq A$ and let $b\e A$. Say $B\subseteq b$ if \zzz
for all $a\e B$, $\A\modelx a\e b$. A formula $\tx\m xyz$, possibly with parameters from \zzz
\A, is an \ul{indicator} for a collection \sss S of initial segments of \A\ if for \zzz
all $a,b\e A$, there exists $I\e\sss S$ such that $a\e I\subseteq b$ iff for all $n\e\w$, \zzz
$\A\modelx\tx\m nab$. Two sets, \sss S and \sss S', are \ul{symbiotic} if for all $a,b\e A$, there \zzz
exists $I\e\sss S$ with $a\e I\subseteq b$ iff there exists $I'\e\sss S'$ with $a\e I'\subseteq b$. \zzz
Say \ul{\w\ codes} an initial segment $I$ if there exists a sequence \sx\ of \zzz
nonstandard length such that $I=\cupx_{i\in\m\w}\sx_i$. \sss S is \ul{cofinal} in \A\ iff for all \zzz
$a\e A$, there exists $I\e\sss S$ with $a\e I$.

    The next lemma is well-known (perhaps due to Gaifman or Puritz), and \zzz
we omit the proof. Given structures \A, \sss C, with $A\subseteq C$, there exists a \zzz
unique set $B$ that $A\subseteq^\text{cof}B\subseteq^\text{end}C$. If $\A\modelx\Sx_1$-\Collection\ and $\A\prece_0\sss C$ \zzz
we have that $\sss B=\sss C\up B$ is a structure. Moreover, if $\A\modelx\Dx_0$-\Separation:

\thm{5.2 Lemma} i.\ If $\A\prece_k\sss C$ and $\Sx_{k+1}$-Collection\ holds in \A, then $\B\prece_k\sss C$.

\noindent
\quad ii.\ If $\A\prece_0\sss B$ and $\Sx_{k+1}$-\Collection\ holds in \A, then $\A\prece_{k+1}\sss B$.\done

\refstepcounter{thisisdumb}\label{Chapter5:5.3}
\thm{5.3 Theorem} \label{Chapter5:Theorem53}Let $k\e\w$, let \A\ be a non-\w-model of \s, let \T\ be any theory coded in \A, and let $\sss S=\myset{I\prece^\text{end}_k\A}{I\modelx\T}$.
\begin{mylist}
\item Suppose \T\ extends \s\ and \Collection\ and that $\Sx_{k+1}$-overspill holds in \A. Then
  \begin{mylist}
  \renewcommand{\labelenumii}{\textbullet}
  \vspace{.8ex}
  \itemsep0ex
  \item \sss S has an indicator, and
  \item \sss S is symbiotic with $\myset{I\e\sss S}{I\text{ r.s.\ and \w\ codes }I}$.
  \item $\sss S\setminus\{A\}$ is nonempty iff the $\Pi_{k+1}$ theory of \A\ is consistent with \T, and
  \item \sss S is cofinal in \A\ iff the $\Sx_{k+2}$ theory of \A\ is consistent with \T.
  \end{mylist}
\item Suppose \A\ is $\Sx_{(k+1)}$-recursively saturated and is locally countable. Then
  \begin{mylist}
  \renewcommand{\labelenumii}{\textbullet}
  \vspace{.8ex}
  \itemsep0ex
  \item \sss S is symbiotic with $\myset{I\e\sss S}{I\text{ r.s.}}$.
  \item \sss S is nonempty iff the $\Pi_{(k+1)}$ theory of \A\ is consistent with \T, and
  \item \sss S equals $\{A\}$ or is cofinal in \A\ iff the $\Sx_{(k+2)}$ theory of \A\ is consistent with \T.
  \item If $\Sx_{k+1}$-\Collection\ either holds in \A\ or is included in \T, \sss S has an indicator, and if we also have either \ssf{QF}-\Foundation\ or \ssf{QF}-\Separation\ in \A\ or \T, then $\sss S\ne\{A\}$.
  \end{mylist}
\end{mylist}

    Before the proof, let us first make several remarks.

    The notion of indicator in arithmetic is due to J. Paris and L. Kirby; \zzz
see Kirby\,\cite{KirbySolo77}. The existence of the indicators for arithmetic indicated \zzz
above is also due to them: they use a game-theoretic argument.

    It will be clear from our proof that if the power-set axiom holds in \zzz
\A, we only need $\Pi_k$-overspill and $\Pi_{(k)}$-recursive saturation in order to \zzz
show the existence of the indicators and the symbiosis.

    For information concerning the order types of these sets of initial \zzz
segments, see theorem \hyperref[Chapter5:5.11]{5.11} below.

    Part (i) (without the condition of \w\ coding $I$ and with \T\ also \zzz
assumed to contain $\Dx_0$-\Separation) actually follows from (ii) by taking \zzz
countable elementary submodels of \A\ and using 5.2 and a result of \zzz
Smory\'nski and Stavi\,\cite{Smor79}, namely that cofinal elementary extensions of r.s.\ \zzz
models of \s\ + \Separation\ are recursively saturated, However, as there are \zzz
short elegant proofs using i-fulfilment, we shall prove (i) directly.

    It is easy to check that if $\sx\e A$, if i\sx\ holds and if \sx\ fulfils \zzz
the universal closure of \tx\vx\nx\tx, then it also fulfils the universal \zzz
closure of
\[
\Ex x\e z\,\tx\imp\Ex\e z\,\big(\tx\ax\Ax y\e x\,\nx\tx(y/x)\m\big)\,.
\]
This places a restriction (at least in the arithmetical case) on the ease \zzz
with which we may construct r.s.\ initial segments in which \Collection\ \zzz
does not hold---for example, our method requires the countability condition.

    Finally, let us remark that in a pre-print (but not in the published) \zzz
version of Paris and Kirby\,\cite{ParisKirby78}, a result is (essentially) proved which is \zzz
almost the dual of ours; let \M\ be a countable model of  $\PA^-_\text{ex}$ as given on page \pageref{Chapter1:DefinitionPAex}, \zzz
let \T\ be any r.e.\ theory extending \PA; then there exists a $k$-elementary \zzz
end-extension of \M\ which is a model of \T\ iff \M\ is a model of \zzz
$\Sx_{k+1}$-\Collection, the $\Pi_{k+1}$ theory of \T, and \ul{\T-provable-$\Pi_{k+1}$-overspill}, \zzz
that is, if $\px x$ is $\Pi_{k+1}$ and for all $n\e\w$, \T\ proves $\px\ol n$, then $\px m$ \zzz
holds in \M\ for some nonstandard $m$.

\noindent
Proof of 5.3: Let us first quickly check the assertions of (i) in the \zzz
order given.

    Let $a,b\e A$ and suppose $a\e I\subseteq b$ for some $I\e\sss S$. Then for each \zzz
$n\e\w$,
\[
\Ex\sx\,\big(\lenn\ax a\e\sx_0\ax\vect\sx\subseteq b\ax\text i\sx\ax\T^\sx\ax\Tr\Sx_k\m^\sx\big) \label{Chapter5:eq:1}\tag{1}
\]
is true in \A, where ``\T'' here is represented by any code for \T. \zzz
Conversely, if \eqref{Chapter5:eq:1} holds for each $n\e\w$ then, using overspill, there \zzz
is an $I$ in \sss S with $a\e I\subseteq b$: just let $I=\cupx_{i\in\m\w}\m\sx_i$ where \sx\ witnesses \zzz
\eqref{Chapter5:eq:1} for some nonstandard \n. As \sx\ also fulfils the universal closure \zzz
of \tx\vx\nx\tx\ for all formulae \tx, it provides a global satisfaction predicate \zzz
for $I$ in \A, and so $I$ is recursively saturated.

    If \sss S is nonempty (respectively, cofinal), then it is clear that \A\ is a model \zzz
of the $\Sx_{k+1}$ ($\Pi_{k+2}$) theory of \T. Conversely, if \A\ is a \zzz
model of the $\Sx_{k+1}$ theory of \T, then \eqref{Chapter5:eq:1}, with references to $a$ and $b$ \zzz
deleted, holds for each $n\e\w$. So, using overspill, \sss S is not empty. \zzz
If \A\ is a model of the $\Pi_{k+1}$ theory of \T, then \eqref{Chapter5:eq:1}, prefixed by a \zzz
universal quantification of $a$ and with $b$ deleted, holds for each $n\e\w$. \zzz
Hence \sss S is cofinal in \A. This establishes (i).

    \refstepcounter{thisisdumb}\label{Chapter5:5.3.ii}Now consider 5.3.ii. Given $n\e\w$ and two closed, increasing sequences \zzz
$\sx,\tau$, say $\sx\prece_n\tau$ if $|\sx|\eqx |\tau|$, for all $i<|\sx|$, $\sx_i\subseteq\tau_i$, and for all \zzz
$i<\lenm\sx$, all $\godel{\px\m\vect v}<n$, and all $\vect v\ve\sx_i$, $(\px\m\vect v)^\sx_i\imp(\px\m\vect v)^\tau_i$. Fix $a,b\e A$. \zzz
Choose some code $t$ for \T: we shall suppose that \T\ contains the \zzz
universal closure of \tx\vx\nx\tx\ for all formulae \tx. For each $n\e\w$, \zzz
define a $\Sx_{(k+1)}$ predicate \ul{$\rrm{extendible}_n(\sx)$}, with parameters $a,b$, and $t$, as follows. \zzz
Let $\rrm{extendible}_0(\sx)$ be $0\eqx 0$, and let $\rrm{extendible}_{n+1}(\sx)$ hold iff
\[
\lenm\sx\gtx\ol\n\ax a\e\sx_0\ax\vect\sx\subseteq b\ax t^\sx\ax\Tr\Sx_k\m^\sx\ax
\Ax i\ltx\lenm\sx\,\Ax x\e\sx_i\,\Ax y\e x\,\Ex \tau\,\big(y\e\tau_0\ax\sx\prece_n\tau\ax\rrm{extendible}_n(\tau)\m\big)\,.
\]
Let \ul{$n$-extendible} be the corresponding informal property: with this notion we are \zzz
trying to capture some of the properties of a sequence defined (externally) \zzz
via satisfaction functions.

    If there exists $I$ in \sss S such that $a\e I\subseteq b$, then it is clear that \zzz
for each $n\e\w$ there exists an \n-extendible sequence in $A$. By \zzz
$\Sx_{(k+1)}$-recursive saturation, we may choose a sequence \sx\ in $A$ which is \zzz
\n-extendible for all $n\e\w$. Let $\sx^0\eqx \sx$. Let $y\e A$ be such that there \zzz
exist $i\e\w$ and $x\e\sx_i$ with $y\e x$. Then for each $n\e\w$ there exists $\tau$ \zzz
such that
\[
y\e\tau_0\ax\sx\prece_n\tau\ax\rrm{extendible}_n(\tau)\,.
\]
By $\Sx_{(k+1)}$-recursive saturation, choose $\tau$ such that this holds for all \zzz
$n\e\w$; let $\sx^1\eqx \tau$. Continuing in some such manner, define a sequence \zzz
$\seq{\sx^i}_{i\in\w}$. Let $B_i=\cupx_{j\in\m\w}\m(\sx^i)_j$. The $B_i$'s form an ascending elementary \zzz
chain of $k$-elementary substructures of \A\ which are recursively saturated \zzz
models of \T, with $B_i\subseteq b$ for all $i\e\w$. Let $I=\cupx_{i\in\m\w}\m\B_i$. Then $I$ is \zzz
also a $k$-elementary substructure of \A\ which is a r.s.\ model of \T, with \zzz
$I\subseteq b$. If we make the successive choices of $x$ and $y$ in a sufficiently \zzz
orderly manner and use the local countability of \A, we may ensure that $I$ \zzz
is also an initial segment. Thus \sss S is symbiotic with its r.s.\ members.

    Suppose that the $\Sx_{(k+1)}$ theory of \T\ holds in \A. Redefine $\rrm{extendible}_n$ \zzz
by deleting the references to $a$ and $b$. In any model of \T, for each \zzz
$n\e\w$, \underline{\smash{there exists a sequence which is \n-extendible}}. Just as in our proof \zzz
of the Completeness Theorem on page \pageref{Chapter1:GodelsProof}, this underlined statement may be \zzz
expressed as a $\Sx_{(k+1)}$ sentence \ul{without} the coding of either sequences or \zzz
formulae and without the use of the $\Sx_k$ satisfaction predicate, and, so \zzz
expressed, it is provable in \T. Hence it must hold in \A, and so for \zzz
each $n\e\w$, $\Ex\sx\,\rrm{extendible}_n(\sx)$ holds in \A. By $\Sx_{(k+1)}$-recursive saturation \zzz
and the above construction, we may conclude that \sss S is nonempty. The \zzz
cofinal case is argued similarly.

    If $\Sx_{k+1}$-\Collection\ holds in either \A\ or \T, we may express the notion \zzz
of $n$-extendible uniformly in $n$ to obtain an indicator. Also in this case, we \zzz
can find $b\e A$ such that $B\subseteq b$ for the \B\ constructed above. If we have \zzz
QF-\Foundation, then $b\notE b$, and if we have QF-\Separation, consider $c=\myset{x\e b}{x\notE x}$. \zzz
Obviously $c\notE b$. The next result shows that, in general, $\Sx_{k+1}$-\Collection\ is \zzz
both a necessary and sufficient condition to ensure $\sss S\ne\{A\}$.\done

\thm{5.4 Lemma} Let $k\e\w$ and let \A\ be a countable model of \s\ + $\Sx_{k+1}$-overspill. \zzz
The following are equivalent.
\begin{mylist}
\item For all formulae $\px\xvect a$ with parameters $\vect a\ve A$, there exists $b\e A$ \zzz
and $B\subseteq b$ such that $\vect a\ve B$, $\B\prec^\text{end}_k\A$, and $\px\xvect a$ is absolute between \A\ and \B.
\item $\Sx_{k+1}$-\Collection\ holds in \A.
\end{mylist}
Proof: Suppose (ii). Replace ``\T'' with ``$\px\xvect a$'' in the above definition \zzz
of \n-extendible. By $\Sx_{k+1}$-\Collection\ we may bound the quantifiers ``$\Ex\tau$'', and \zzz
then (i) follows by the above construction.

    Suppose (i), and suppose $\Av x\ve a\,\Ev y\m\tx$ holds in \A, where $\tx\e\Sx_{k+1}$. \zzz
Choose $b\e A$, $B\subseteq b$ such that $\B\prec^\text{end}_k\A$ and $\B\modelx\Av x\ve a\,\Ev y\m\tx$. Then \zzz
$\A\modelx\Av x\ve a\,\Ev y\ve b\,\tx$. Thus (ii).\done

\thm{5.5 Corollary} \label{Chapter5:5.5}Let \A\ be a model of \s\ and $\Sx_{k+1}$-overspill, and let \T\ \zzz
be a \ul{complete} theory extending \s\ plus \Collection\ plus \Separation. The following are \zzz
equivalent.
\begin{mylist}
\item There exists a $k$-elementary initial segment of \A\ which is a \zzz
nonstandard model of \T.
\item \A\ is a model of $\T\cap\Sx_{k+1}$, and for all $n\e\w$, $\T\cap\Sx_n$ is \zzz
coded in \A.
\end{mylist}

    Before the proof, let us make a number of remarks.

    This result is (in the arithmetical case) the dual of a result of \zzz
Wilkie\,\cite{Wilk77}. (The proof below is modelled on Lessan's\,\cite{Less78} proof of this dual.) \zzz
Together these results answer a question in the introduction of Friedman\,\cite{Frie73}: \zzz
are the standard systems of the nonstandard models of \PA\ which contain all \zzz
arithmetic sets the same as those of models of the theory of \bb N?

    Let us consider more closely the arithmetical case when \T\ extends \zzz
\PA: we shall quickly sketch two alternative proofs of 5.5. First, by 5.3.i, \zzz
we can construct a chain of nonstandard initial segments
\[
\A\succ_k\B_0\succ_{k+1}\B_1\succ_{k+2}\B_2\succ_{k+3}\dotsb \label{Chapter5:eq:2}\tag{2}
\]
where $\B_k$ is a model of $\T\cap\Sx_{n+k+2}\,\cup\,\PA$. Let $I=\capx_{i\in\m\w}\m\B_i$. As \T\ has \zzz
definable satisfaction functions, it is easy to check that $\I\prec_k\A$ and \zzz
that \I\ is a model of \T. It is also easy to ensure that \I\ is nonstandard.

    For the second proof, let \B\ be the minimal model for \T, or, if \zzz
$\T=\Th(\bb N)$, let \B\ be a conservative extension of \bb N (that is, its \zzz
standard system consists of exactly the arithmetical sets; see Phillips\,\cite{Phil74} \zzz
or 5.9 below). A subset of \w\ is coded in \B\ iff it is recursive \zzz
in $\T\cap\Sx_n$ for some $n\e\w$. Thus the set of reals of \B\ is included in \zzz
that of \A, and so we may, by Friedman\,\cite{Frie73}, construct a $k$-elementary \zzz
embedding of \B\ into \A. Then the initial segment of \A\ determined \zzz
by the image of \B\ is by 5.2 a $k$-elementary substructure of \A\ and is a \zzz
model of \T.

\noindent
Third proof of 5.5: The left to right implication is easy, for \s\ + \Collection\provex  \zzz
\Induction, which provides the overspill required to code the true $\Sx_n$ sentences. \zzz
For the converse, by 5.2 we may suppose $A$ is countable. Let $\seq{c_i}_{i\in\m\w}$ be \zzz
an infinite sequence of new constants. Let $\seq{\px_i}_{i\in\m\w}$ be a listing of \zzz
all sentences of $\sss L\cup\seq{c_i}_{i\in\m\w}$ of the form $\Ex x\m\tx$, and suppose is $\px_i$ is in $\Sx_{i+k}$ and \zzz
only contains the constants  $\seq{c_j}_{j\le i}$. We shall define a descending chain of \zzz
initial segments as in \eqref{Chapter5:eq:2} above, and an interpretation of $\seq{c_i}_{i\in\m\w}$ as follows. \zzz
(We shall denote the interpretation of $c_i$ also by $c_i$.) Choose $\B_0\prec_k\A$ to \zzz
be any r.s.\ initial segment which is a model of $\T\cap\Pi_{k+3}$, and let $c_0$ be any \zzz
nonstandard integer of $\B_0$. Suppose we have defined $\B_i$, a r.s.\ model of \zzz
$\T\cap\Pi_{i+k+3}$, and have interpreted $\seq{c_j}_{j\le i}$ in $\B_i$. If $\px_i=\Ex x\m\tx$ holds in $\B_i$, \zzz
let $c_{i+1}$ be interpreted by a witness for $\px_i$; otherwise, let $c_{i+1}$ be the first \zzz
element of some fixed listing of $A$ which is a member of $\B_i\setminus\{c_0\m,\dots,c_i\}$. \zzz
Now choose $\B_{i+1}$ to be any $(k+i)$-elementary r.s.\ initial segment of $\B_i$ \zzz
which includes $\{c_0\m,\dots,c_{i+1}\}$ and is a model of $\T\cap\Pi_{i+k+4}$.

    Let $I=\capx_{i\in\m\w}\m\B_i=\myset{c_i}{i\e\w}$. By construction, $I$ is a nonstandard \zzz
model of \T.\done

\thm{5.5.ii Corollary} Let \A\ be a countable model of \s\ plus $\Sx_{(k+1)}$-recursive \zzz
saturation, and let \T\ be a complete theory extending \s. Then the following \zzz
are equivalent.
\begin{mylist}
\item There exists a $k$-elementary initial segment \B\ of \A\ which is \zzz
a model of \T\ and which is $\Sx_{(n)}$-recursively saturated for all $n\e\w$.
\item \A\ is a model of $\T\cap\Sx_{(k+1)}$, and for all $n\e\w$, $\T\cap\Sx_{(n)}$ is \zzz
coded in \A.
\end{mylist}
If \A\ is also a model of $\Sx_{k+1}$-\Collection, we may ensure that $\B\ne\A$.

\noindent
Proof: The proof is as above but, as the intersection of recursively \zzz
saturated initial segments may not be recursively saturated, we must in \zzz
our construction also add witnesses to ensure  $\Sx_{(n)}$-recursive \zzz
saturation.\done

    We can extend 5.3 as follows. Let \T\ be a theory in any language \zzz
in which the language of arithmetic may be interpreted. Let \M\ be a \zzz
nonstandard model of \PRA\ and suppose \T\ is coded in \M. Since (an \zzz
axiomatization of) the arithmetical consequences of \T, $\rrm{arith}(\T)$, will \zzz
also be coded in \M, we may use 5.3 to consider, say,
\[
\sss S=\myset{\I\prece^\text{end}_k\A}{\I\modelx\rrm{arith}(\T)}\,.
\]
There is however another, perhaps more elegant approach using $\frac12\ast$fulfilment, \zzz
and this also allows us to consider
\[
\myset{\I\e\sss S}{\text{\I\ is expandable to a model of \T}}\,.
\]
In the countable case and for many theories \T\ (e.g.\ any theory extending \zzz
$\Dx^1_1\HY\CA\up$), this coincides with \myset{\I\e\sss S}{\text{\I\ r.s}}, so we do not have any \zzz
advantage over 5.3.ii. But we may improve 4.3.i to the following, also \zzz
obtained independently by Kirby, McAloon, and Murawski\,\cite{KirbyMcAlMura79} for theories \zzz
extending $\Pi^0_1\HY\CA\up$.

\thm{5.6 Theorem} Let $k\e\w$ and let \M\ be a model of \PRA\ plus \zzz
$\Sx^0_{k+1}$-overspill. Let \T\ be any theory coded in \M\ extending \s\ and \zzz
also extending the schema:
\begin{flalign*}
&\ssf{Arith-Collection}\!: &&\Ax m\ltx n\,\Ex p\m\tx\imp\Ex q\,\Ax m\ltx n\,\Ex p\ltx q\,\tx\,. &&&
\end{flalign*}
Let $\sss S=\myset{\I\prece^\text{end}_k\A}{\I\modelx\rrm{arith}(\T)}$. Then \sss S has an indicator, and is \zzz
symbiotic with
\[
\myset{\I\e\sss S}{\text{\I\ r.s., \w\ codes $I$, and \I\ is expandable to a model of \T}}\,.
\]
$\sss S\setminus\{\M\}$ is nonempty iff the $\Pi^0_{k+1}$ theory of \M\ is consistent \zzz
with \T, and \sss S is cofinal in \M\ iff the $\Sx^0_{k+2}$ theory of \M\ is \zzz
consistent with \T. \done

\noindent
The proof of 5.6 is clear using the notion of $\frac12\ast$fulfilment, and we omit \zzz
it, only pausing to remark that the schema of \Induction\ implies that of \zzz
\ssf{Arith-Collection}.

    We shall consider two more model-theoretic ``tricks''. The first is \zzz
closely related to $\frac12\ast$fulfilment and 5.6. Let $\PA_\in$ be the version of \zzz
Peano arithmetic formulated in the language with a single binary relation \zzz
$\in$. Let \T\ be any consistent, recursive theory in the language $\{\in\}$ which \zzz
proves that there exists no $\in$-loops, For example, \T\ might be \ZFC\ or \zzz
some theory of analysis over \HF. A well-known folklore result, easily \zzz
proved, using the compactness theorem, is that any model of \T\ may be embedded \zzz
in a model of $\PA_\in$. Using fulfilment, we easily obtain the dual result:

\thm{5.7 Theorem} Any nonstandard model \A\ of $\PA_\in$ has a substructure \B\ \zzz
which is a recursively saturated model of \T. Moreover, if we assume that \zzz
\T\ is strong enough to perform some coding, that \T\ includes \ssf{Arith-Collection}, \zzz
and that \A\ is a model of the theory of \T, then we may \zzz
choose \B\ so that $\HF^\B$ (i.e.\ $\myset{x\e B}{x\e^\B\HF\m^\B}$) is an initial segment \zzz
of \A. And the usual hierarchical and cofinal variants also hold.

\noindent
Proof: We simply note that any finite poset with no loops is isomorphic \zzz
to a subset of $\seq{\HF,\in}$. \done

    For our our second ``trick'', let \T\ be a theory in any finite language \zzz
in which we may interpret arithmetic, and let be a nonstandard \zzz
model of \PRA\ which codes \T. In \T, add a constant \ol{\smallerP\w} for the \zzz
definable type consisting of subsets of \ol\w. We have been considering \zzz
models \A\ of \T\ ``coded'' in \M\ in which $\ol\w\m^\A$ (i.e.\ \myset{x\e A}{x\e^\A\ol\w\m^\A}) is \zzz
an initial segment of \M. In such a case we have
\[
\Myset{\myset{x\e\ol\w\m^\A}{x\e^\A X}}{X\e^\A\ol{\smallerP\w}\m^\A}\,\subseteq\;\sss R(\m\ol\w\m^\A\mx,\m\M\m)\,,\label{Chapter5:eq:3}\tag{3}
\]
where this latter set is the collection of those subsets of $\ol\w\m^\A$ which are \zzz
coded in \M. Kirby, McAloon, and Murawski\,\cite{KirbyMcAlMura79} noticed if \T\ is a \zzz
second-order theory extending $\Pi^0_1\HY\CA\up$ and if \M\ is countable, then the \zzz
collection of initial segments \I\ of \M\ which are expandable to a \zzz
model \A\ of \T\ is symbiotic with its subset consisting of those $\I=\ol\w\m^\A$ \zzz
for which the relation \eqref{Chapter5:eq:3} is in fact equality. We may extend this result \zzz
as follows.\footnote{Footnote added 2019: This next result and its two lemmas have nothing to do with fulfillability.}

\thm{5.8 Theorem}\label{Chapter5:5.8} Let $k\geX 1$, and suppose \A\ is a model of
\begin{mylist}
\renewcommand{\labelenumi}{\textbullet}
\itemsep0ex
\item $\Dx^0_0\HY\CA\up$, $\Dx^0_0$\ssf{-Arith-Collection} (where these two schemata are allowed to have second-order parameters),
\item an axiom asserting that the $k^\text{th}$ Turing jump of the empty set exists, and
\item \WKL\ (an axiom asserting that every infinite binary tree has an infinite branch).
\item Further suppose that $\Pi^0_1$-overspill holds in \A, and that
\item $\ol{\smallerP\w}\m^\A$ is countable.
\end{mylist}
Then for each $c\e\ol\w\m^\A$ there exists a structure \B\
\begin{mylist}
\renewcommand{\labelenumi}{\textbullet}
\itemsep0ex
\item isomorphic to \A\ such that
\item $\ol\w\m^\B$ is a $(k-1)$-elementary initial segment of $\ol\w\m^\A$,
\item $\ol{\smallerP\w}\m^\B=\sss R(\ol\w\m^\A,\,\ol\w\m^\B)$, and
\item the isomorphism fixes $c$.
\end{mylist}
If in addition we have that $\Sx^1_1$-overspill holds in \A\ and that \zzz
$\Ex X\m.\m X=\emptyset^{_{(\ol k)}}$ holds in \A\ for all $k\e\w$, then \zzz
we may require that $\ol\w\m^\B\prec\ol\w\m^\A$.

    Before the proof, we need a definition and two lemmas. \zzz
First note that without loss of generality we can assume that \A\ is \zzz
a structure for the language of analysis, say $\A=\seq{\M,\sss X}$.  A pseudo-standard formula of \M\ is \zzz
one of the form $\px(\ol a,\ol b,\dots,u,v,\dots)$ where $\px(x,y,\dots,u,v,\dots)$ is a \zzz
formula and $a,b,\dots$ are elements of \M. If we add new constants $\myset{d_a}{a\e\M}$ \zzz
to our metalanguage, with each pseudo-standard formula $\px(\ol a,\ol b,\dots)$ we may associate \zzz
a formula $\px(d_a,d_b,\dots)$. For any $\T\subseteq M$, let \ssf{p.s.}(\T) be the collection \zzz
of those $\px(d_a,d_b,\dots)$ for which $\px(\ol a,\ol b,\dots)$ is a pseudo-standard formula whose \zzz
code is in \T.

    \refstepcounter{thisisdumb}\label{Chapter5:5.9}
    The first of the two lemmas, which does not require that \M\ be nonstandard, is a \zzz
common generalization of Scott\,\cite{Scott62} and Friedman\,\cite{Frie73} and the generalizations \zzz
of the MacDowell-Specker Theorem\,\cite{MacDSpec61} by Phillips\,\cite{Phil74} and Gaifman\,\cite{Gaif76}. \zzz
A slightly less general version was first observed by Kirby, McAloon, and \zzz
Murawski\,\cite{KirbyMcAlMura79}; their proof is an ultrapower construction which requires \zzz
that \A\ be a model of $\Pi^1_1\HY\CA\up$.

\thm{5.9 Lemma} Let $\seq{\M,\sss X}$ be a countable model of $\PA^-_\text{ex}$, $\Dx^0_0$-\Induction, \zzz
$\Dx^0_0$-\ssf{Arith-Collection}, and \WKL.
\begin{mylist}
\item Let \T\e\sss X be such that in $\seq{\M,\sss X}$ \T\ is a consistent set of \zzz
sentences of arithmetic. Then there exists an end-extension \sss N of \M\ \zzz
which is a model of \ssf{p.s.}(\T) in the natural sense with $d_a$ being interpreted \zzz
by $a$ for all $a\e M$, and is such that $\sss R(\sss N,\M)=\sss X$. If $\M=\bb N$ or \zzz
if $\seq{\M,\sss X}$ is a model of $\Pi^0_1$-overspill, we may choose \sss N to be recursively \zzz
saturated.
\item  Let $\T\subseteq M$ be such that \ssf{p.s.}(\T) is a complete theory in the \zzz
language of arithmetic and for all $k\e\w$, $(\T\cap\Sx^0_k)\e\sss X$ and $\M\modelx\Con(\T\cap\Sx^0_k)$. \zzz
Then there exists an end-extension \sss N which is a model of \ssf{p.s.}(\T) and \zzz
is such that $\sss R(\M, \sss N)=\sss X$.
\end{mylist}
Proof: For (i) we simply have to check that the usual Henkin construction \zzz
works for $\M\ne\bb N$. Since \T\ proves each true quantifier-free pseudo-standard sentence \zzz
of \M, we can suppose that these are included in \T. Let $\seq{c_i}_{i\in\w}$ \zzz
be a list of new constants, and let $\sss L(c_0\m,\dots,c_i)$ be the language of \zzz
arithmetic augmented with the constants $c_0\m,\dots,c_i$. We shall define an \zzz
increasing sequence $\seq{\T_n}_{n\in\w}$, where for each $n\e\w$, $\T_n\e\sss X$ is a \zzz
theory of $\sss L(c_0\m,\dots,c_i)$ for some $i\e\w$ which is consistent in \M. \zzz
Let $\T_0=\T\cup\myset{c_0\gtx\ol a}{a\e M}$. This is consistent, because given any \zzz
proof of a contradiction from $\T_0$ we can replace $c_0$ by \ol a for some \zzz
large $a\e M$, obtaining a contradiction from \T. Define $\T_n$ inductively \zzz
as follows.

    Suppose $n\eqx 4m+1$ and suppose our language has only $\Ex$, $\vx$ and $\nx$ \zzz
as logical symbols. Let \px\ be the $m^\text{th}$ member of some enumeration of \zzz
the formulae of $\cupx_{i\in\w}\m\sss L(c_0\m,\dots,c_i)$. Let $\T_{n+1}$ be $\T_n+\nx\px$, if this is \zzz
consistent; otherwise, if $\px=\Ex x\m\tx x$ let $\T_{n+1}$ be $\T_n+\px+\tx c$ for some \zzz
$c$ not yet considered, and for \px\ in other forms simply let $\T_{n+1} = \T_n+\px$.

    Suppose $n\eqx 4m+2$. Let $X$ be the first element of some enumeration \zzz
of \sss X not yet considered, and let $c$ be some constant not yet considered. \zzz
Let $\T_{n+1} = \T_n\cup\myset{\ol a\e c}{a\e X}\cup\myset{\ol a\notE c}{a\notE X}$.

    Suppose $n\eqx 4m+3$. Since we are already using ``$\e$'' to denote \zzz
membership between \M\ and \sss X, we shall let $\hat\e$ denote the membership \zzz
relation between \ul{integers}, as defined on page \pageref{Chapter1:DefinitionMembership}. Choose $X\e\sss X$ so that
\[
\T_{n+1} = \T_n\cup\myset{\ol a\,\hat\e\, c_m}{a\e X}\cup\myset{\ol a\,\hat\notE\, c_m}{a\notE X}
\]
is consistent. We can do so, for consider the tree of binary sequences
\[
\myset{b\e M}{\text{there is no proof of 0\m=1 with code$\m\ltx\lenm b$ from } \T_n+\myset{\ol a\,\hat\e\,c_m}{b_a\text{=1}}+\myset{\ol a\,\hat\notE\,c_m}{b_a\text{=\,0}}}\,.
\]
This is in \sss X, and it is infinite in $\seq{\M,\sss X}$, because otherwise by $\Dx^0_0$\ssf{-Arith-Collection} \zzz
we could piece together proofs to obtain a contradiction from \zzz
$\T_n$. By \WKL, we can choose an infinite branch in \sss X to obtain $X$.

    Suppose $n\eqx 4m+4$. It is at this stage that we ensure that the end-extension \zzz
is r.s. We can suppose $\M\ne\bb N$, for otherwise we simply use \zzz
the standard argument. Let $\tx x$ be the $m^\text{th}$ $\Dx^0_0$ formula, and suppose that \zzz
for some $i\e\w$ the set $\myset{x\e\w}{\tx x}$ consists of the codes of formulae \zzz
of the form $\px xc_0\dots c_i$. Suppose for each $p\e\w$, $\Ex x\aax\myset{\px x}{\tx\,\godel{\mx\px\mx}\ax\godel{\mx\px\mx}\ltx p}$ \zzz
is consistent in \M\ with $\T_n$. Then by $\Pi^0_1$-overspill, this is consistent \zzz
for some nonstandard $p\e M$. Let $c$ be some constant not yet considered, \zzz
and let $\T_{n+1}=\T_n\cup\myset{\px c}{\tx\,\godel{\mx\px\mx}\ax\godel{\mx\px\mx}\ltx p}$.

    This completes our construction. The remainder of the proof is standard and \zzz
is left to the reader.

    The proof of (ii) is obtained by a simple modification of the above, \zzz
essentially due to D.~Jensen and A.~Ehrenfeucht\,\cite{Jens76} and, independently, \zzz
Guaspari\,\cite{Guas79}. But since we shall not require (ii), the proof is omitted.\done

    The next lemma is an immediate corollary of the embedding techniques \zzz
of Friedman\,\cite{Frie73}, Wilkie\,\cite{Wilk7x}, and Wilmers\,\cite{Wilm77}.

\thm{5.10 Lemma} Let \M\ and \sss N be countable.
\begin{mylist}
\item Suppose \M\ and \sss N are models of $\PA^-_\text{ex}$ + $\Sx_{k+1}$-overspill \zzz
with the same standard systems and with $\Th_{\Sx_{k+1}}(\sss N)\subseteq\Th_{\Sx_{k+1}}(\M)$. \zzz
Also suppose $\sss N\modelx\Sx_{k+1}\HY\Collection$. Then \sss N is isomorphic to a $k$-elementary \zzz
initial segment of \M.
\item  If \sss N and \M\ are elementarily \zzz
equivalent, recursively saturated models of \PA\ with the same standard \zzz
system, then \sss N is isomorphic to a proper elementary initial segment of \M.\done
\end{mylist}
Proof of 5.8: By $\Ex X.\,X=\emptyset^{(k)}$, we have that the set of true $\Pi^0_k$ sentences \zzz
of \M\ is an element of \sss X. By Application (xv) on page \pageref{Chapter4:Applicationxv} of Chapter~IV, $\Sx^0_k$-\Induction\ \zzz
implies the consistency of these true $\Pi^0_k$ sentences plus $\Sx_k$-\Collection, and so by \zzz
5.9 we may choose a recursively saturated $k$-elementary end-extension \sss N \zzz
of \M\ which is a model of $\Sx_k$-\Collection\ and is such that $\sss R(\M,\sss N)=\sss X$. \zzz
Now for any $c\e M$, $\Th_{\Sx_k}(\sss N,c)=\Th_{\Sx_k}(\M,c)$ and so by 5.10 there is an \zzz
embedding of \sss N into a proper initial segment of \M\ fixing $c$.

    For 5.8.ii, we note that there is a recursive predicate $P(x,k,X)$ such \zzz
that $\myset{x}{P(x,k,\emptyset^{(k)})}$ is the set of true $\Pi^0_k$ sentences. So if \M\ is \zzz
a model of \PA\ and if $\Sx^1_1$-overspill and $\Ex X.\,X=\emptyset^{(k)}$ for all $k\e\w$ hold \zzz
in $\seq{\M,\sss X}$, then there is some nonstandard $k\e M$ for which
\[
\Ex X\,\big(X=\emptyset^{(k)}\ax\Con\m(\myset{x}{P(x,k,X)})\m\big)
\]
holds. Now by 5.9, we may choose a r.s.\ elementary end-extension \sss N of \zzz
\M\ such that $\sss R(\M,\,\sss N)=\sss X$, and then by 5.10, embed \sss N in \M.

    This completes the proof of 5.8.\done

    \refstepcounter{thisisdumb}\label{Chapter5:5.11}
    Our final theorem considers the order types of various classes of \zzz
elementary initial segments of r.s.\ models of \PA.

\thm{5.11 Theorem} Let \M\ be a recursively saturated model of \PA, and let
\begin{align*}
A=&\myset{I\prece^\text{end}\M}{\text{$I$ not r.s.}}\\
B=&\myset{I\prece^\text{end}\M}{\text{$I$ r.s. and \w\ codes $I$}}\\
C=&\myset{I\prece^\text{end}\M}{\text{$I$ r.s.}}
\end{align*}
Then:
\begin{mylist}
\item For all $I\prece^\text{end}\M$, $I\e A$ iff there exists $a\e I$ such that the \zzz
elements of $M$ definable from $a$ are cofinal in $I$.
\item There exists an order-preserving map from $M$ into $A$ (ordered by \zzz
inclusion), and so by (i), the cardinality of $A$ is equal to that of $M$.
\item For all $I\e A$, $I=\capx\myset{J\e B}{I\subsetneq J}$, and so $A$ is densely ordered.
\item $C=\myset{\cupx X}{\emptyset\ne X\subseteq A,\ X \text{ has no greatest element}}$.
\item There exists an order-preserving map from $M$ into $B$, and so the \zzz
cardinality of $B$ equals that of $M$.
\item For all $I\e B$, $I=\capx\myset{J\e B}{I\subsetneq J}=\cupx\myset{J\e B}{J\subsetneq I}$, and so $B$ is \zzz
densely ordered.
\item $C=\myset{\cupx X}{\emptyset\ne X\subseteq B}$, and so by (vi), if $M$ is countable then $C\setminus\{M\}$ \zzz
has the order-type of the real numbers.
\end{mylist}

    Let us first make some remarks. Many of these results are well-known \zzz
for countable $M$: see, for example, Barwise\,\cite{Barwx75} or Schlipf\,\cite{Schl78}. Part (i) \zzz
is essentially due to W.~Marek and H.~Kotlarski; see Kotlarski\,\cite{Kotl78}. For the \zzz
case when \M\ is countable, much of (ii) to (vii) has also been independently \zzz
obtained by Kotlarski\,\cite{Kotl78} and Murawski\,\cite{Mura78}.

    We remark that Paris and Kirby\,\cite{ParisKirby78} show that if $\M\prec^\text{end}\sss N$ with $\M\modelx\s$, \zzz
then \Collection\ holds in \M\ and \sss N, and so it is necessary for us to \zzz
consider models of \PA\ in this result. It will be obvious, however, that \zzz
an analogue of 5.11 could be given for \ul{rank extensions} of r.s.\ models of \zzz
\ZF\ or for \ul{$L$-extensions} of r.s.\ models of $\ZF^-+V=L$. The set-theoretic \zzz
case has in addition another natural class of initial segments symbiotic \zzz
with those above: those \ul{internal} sets $B\e A$ for which $B\prec A$.

\noindent
Proof of 5.11: We shall only consider i-fulfilment in this proof, as \zzz
given on page \pageref{Chapter1:Definitionifulfilment}. Let \ssf{True}(\sx) \zzz
be the schema
\[
\Ax\m i\ltx\lenm\sx\;\Av x\lex\sx_i\big(\px\xvect x\imp(\px\xvect x)^\sx_i\big)\,,
\]
where \px\ has the free variables indicated.

    We shall first prove (vii). Suppose $I\e C$. For each $a\e I$ choose a \zzz
$\sx\e I$ to realize the recursive type
\[
\{a\ltx\sx_0\}\cup\myset{\lenm\sx\gtx\ol n}{n\e\w}\cup\ssf{True}(\sx)\,,
\]
and let $J_a=\cupx_{i\in\m\w}\m\sx_i$. Then $J_a\e B$ and $I=\cupx_{a\in I}\m J_a$. Conversely, it is \zzz
clear that recursive saturation is preserved under unions of elementary \zzz
chains.

    This also gives the second part of (vi). For the first, let $I\e B$ \zzz
and let $s$ be a sequence coded in $M$ such that $I=\cupx_{i\in\m\w}\m s_i$. Let $a\e M\setminus I$; \zzz
by recursive saturation choose $t\e M$ to be an increasing sequence such \zzz
that $\lenm t$ is nonstandard, $t_{\lenm t}\lex\lenm s$, $s_{t_{\lenm t}}\ltx a$, and so that $\ssf{True}(\seq{s_{t_i}}_{i\le\lenm t})$ \zzz
holds; and, let $J_a=\cupx_{i\in\w}\m s_{t_i}$. Then $I=\capx_{a\in M\setminus I}\,J_a$, and (vi) is proved,

    For each $a\e M$, let $I_a$ be the initial segment determined by those \zzz
elements definable from $a$. Then $I_a\prec \M$ by Lemma 5.2, and $I_a$ is not \zzz
recursively saturated, for the type $\myset{\Ex y\m\tx ay\imp\Ex y\ltx x\,\tx ay}{\tx zy\text{ any formula}}$ \zzz
is not realized in $I_a$. $A$ and $B$ are easily seen to be symbiotic, because \zzz
for any $X\subseteq A$ with no greatest element there exists $Y\subseteq B$ such that \zzz
$\cupx X=\cupx Y$, and so $\cupx X$ is r.s. Thus each member of $A$ has the form $I_a$, \zzz
and so (i) holds.

    Let $I\eqx I_a$, let $b\e M\setminus I$, and let $c$ satisfy
\begin{align*}
\Ex y\,\tx xy &\imp\Ex y\ltx b\,\tx xy\quad\text{and}\\
\Ex y\,\tx ay &\imp\Ex y\ltx x\,\tx ay
\end{align*}
for all formulae $\tx xy$. Then $I\subsetneq I_c<b$, and as $b$ was arbitrary, this \zzz
proves (iii). (iv) is easy: if $J\e C$, $J=\cupx_{a\in J}\m I_a$ where $I_a$ is as above.

    Only (ii) and (v) remain to be proved. Imagine, for a moment, that \zzz
there exists an \M-infinite sequence $\sx=\seq{\sx_i}_{i\in M}$ satisfying $\True(\sx)$. \zzz
(If \M\ is \ul{resplendent} then such a sequence exists: we can let \sx\ be a \zzz
cofinal set of indiscernibles for which $\seq{\M,\sx}$ is a model of full induction; \zzz
see Schlipf\,\cite{Schl78}.) Then we could define our required maps easily. Define \zzz
$M\rightarrow A$ by letting $i\mapsto I_{\sx_i}$. This is monotone, for if $\tx xy$ is any \zzz
formula,
\[
\Ex y\,\tx\m\sx_i y\;\ \text{iff}\;\ \Ex y\ltx\sx_{i+1}\,(\tx\m\sx_i y)^\sx_i\;\ \text{iff}\;\ \Ex y\ltx\sx_{i+1}\,\tx\m\sx_i y\,.
\]
Define the map $M\rightarrow B$ by $i\mapsto\cupx_{j\in\w}\m\sx_{ki+j}$, where $k$ \zzz
is some fixed nonstandard element of $M$. This is clearly monotone.

    I do not know whether such a sequence exists in general, but we shall \zzz
construct a suitable alternative. Let $\Gamma(\sx,x)=\{\lenm\sx=x\}\cup\{x\ltx\sx_0\}\cup\True(\sx)$. \zzz
Let $\aaa,\bbb,\gamma$ range over the (real) ordinals, and let \aaa\ be the cofinality \zzz
of $M$, that is, \aaa\ is the least cardinal such that there exists \zzz
a cofinal subset of $M$ of cardinality \aaa. Choose a sequence of internal \zzz
sequences $\seq{\sx^\bbb}_{\bbb<\aaa}$ as follows. Let $\sx^0$ satisfy $\Gamma(\sx^0,1)$, let $\sx^{\bbb+1}$ \zzz
satisfy $\Gamma(\sx^{\bbb+1},\sx^\bbb)$, and if $\gamma\ltx\aaa$ is a limit ordinal, choose $x\e M$ \zzz
greater than $\sx_\bbb$ for all $\bbb\ltx\gamma$, and let $\sx^\gamma$ satisfy $\Gamma(\sx^\gamma,x)$. Define \zzz
an increasing sequence of functions $\seq{f_\bbb}_{\bbb<\aaa}$, each mapping an initial \zzz
segment of $M$ into $M$, as follows. Let $f_0$ be empty; let
\[
f_{\bbb+1}(x)=
\begin{cases}
f_\bbb(x)\,,       &\text{if }x\e\domain(f_\bbb)\,,\\
(\sx^{\bbb+1})_x\,,&\text{if }x\lex\lenm{\sx^{\bbb+1}}\text{ and }x\notE\domain(f_\bbb)\,,\\
\text{undefined,}  &\text{otherwise;}
\end{cases}
\]
and if $\gamma$ is a limit ordinal, let $f_\gamma=\cupx_{\bbb<\gamma}\m f_\bbb$. Now define maps $M\rightarrow A$ \zzz
and $M\rightarrow B$ as before, but replacing by \sx\ by $f_\aaa$.

    This completes the proof of 5.11.\done
\cleardoublepage
\thispagestyle{plain} 
\sectioncentred{More Model-theoretic Applications: \texorpdfstring{\w}{omega}-models}

    We shall study in this chapter \w-models of various theories, \zzz
(including theories in the infinitary language \Lw) which extend our base theory \zzz
\s\ given on page \pageref{Chapter1:DefinitionMembership}. To give an illustration of our techniques, to compare them with other \zzz
known techniques, and to help the reader understand the more general results \zzz
given later, we shall first consider the example of \w-models of the language \zzz
of second-order arithmetic, that is, analysis.

    Let \BI, Bar-Induction, be the schema:
\begin{flalign*}
&\BI\!: &&\Ax X\,\PPL \wf(\prec_{_X})\imp\ppl\Ax n\,\big((\Ax m\prec_{_X}\!n\;\tx m)\imp\tx n\big)\imp\Ax n\m\tx n\ppr\m\PPR &&&
\end{flalign*}
where $\prec_{_X}\,=\myset{(m,n)}{\seq{m,n}\e X}$. Let us consider a theory \T\ in the \zzz
language of analysis extending \s\ + $\Sx^1_1$-\BI\ which has an \w-model. To keep \zzz
things simple, suppose \T\ is r.e., although for the following results it \zzz
suffices to let \T\ be $\Pi^1_1$.

    Gandy, Kreisel, and Tait\,\cite{Gandy60} showed:

\thm{6.1 Theorem} \T\ does not have a \ul{minimum} model.\done

\noindent
Here we are considering \w-models ordered by \ul{inclusion}. Briefly, they \zzz
showed that the intersection, of all models of \T\ consists of the hyper-arithmetic \zzz
sets, \rrm{HYP}, and as $\Sx^1_1$-\BI\ is false in \rrm{HYP}, \T\ does not have a \zzz
minimum model.

    Friedman\,\cite{Frie73} proves that:

\thm{6.2 Theorem} If \T\ extends the full schema of \AC, \T\ does not have a \zzz
\ul{minimal} model.

\noindent
Simpson\,\cite{Simp73} asks whether weaker conditions on \T\ suffice for this result. \zzz
Indeed, we shall show:

\thm{6.3 Theorem} \T\ does not have a minimal model.

    We shall first give a very quick sketch of Friedman's result 6.2 (since \zzz
the existing versions in print are unnecessarily complicated) and then \zzz
briefly explain the modification which allows us to remove \AC.

    First note that as $\Ex X\,\ppl ``\myset{(X)_n}{n\e\w}\modelx\T\m\m"\ppr $ is true in any \bbb-model, \zzz
we can restrict our attention to non-\bbb-models. Fix some non-\bbb-model \A\ \zzz
of \T\ (i.e.\ $\A\subseteq\sss P\w$ and $\seq{\w,\A,\dots}\modelx\T$)\footnote{Footnote added 2019:
That is, for \w-models of analysis, the same notation is used for both the model and
its \sss P\w\ component.} and choose some non-well-founded \zzz
linear ordering $\prec\mx\e\,\A$ for which $\A\modelx\wf(\prec)$. Let \aaa\ be the \zzz
ordinal of the well-founded part of $\prec$.

    First, to give Friedman's proof, suppose \A\ is a model of \AC. We \zzz
consider $\tau\e\A$ such that $\tau$ is (the code of) a tree of finite sequences \zzz
$\sx=\seq{\sx_0,\dots,\sx_{\smash{\lenm\sx}}}$ of sets for which
\begin{mylist}
\item $\ol{\sx_i}\modelx(\T\cap i)$, for all $i\lex\lenm\sx$, and
\item $\ol{\sx_i}\prece_i\ol{\sx_{i+1}}$, for all $i\ltx\lenm\sx$,
\end{mylist}
where in general in this chapter $\ol X=\myset{(X)_n}{n\e\w}$, and $\prece_i$ means $\Sx^1_i$-elementary \zzz
substructure.

\thm{Claim~1} For all $\bbb\ltx\aaa$ there exists such a $\tau\e\A$ of rank \bbb.

\noindent
Theorem 6.2 follows from this claim. For as \T\ includes $\Sx^1_1$-\BI, we may \zzz
apply overspill to obtain a tree $\tau$ of nonstandard rank. $\tau$ has (in the \zzz
real world) an infinite branch $\seq{\sx_0,\sx_1,\sx_2,\dots}$. Let $\B=\cupx_{i\in\w}\m\ol{\sx_i}$. \zzz
Then it is clear that $\B\subseteq\A$ and $\B\modelx\T$. That $\A\ne\B$ follows by \zzz
Cantor's diagonal argument.

    The proof of the claim is by induction on \bbb, and to be able to carry \zzz
out the inductive step, we need to prove a bit more. Say a finite \zzz
sequence \sx\ is \ul{\good} if (i) and (ii) hold, and
\begin{mylist}
\setcounter{enumi}{2}
\item $\ol{\sx_i}\prece_i\A$, for all $i\lex\lenm\sx$.
\end{mylist}
\thm{Claim~2} For each $\bbb\ltx\aaa$ and each \good\ sequence $\sx\e\A$ there exists such a \zzz
tree $\tau\e\A$ of rank $\bbb+\lenm\sx$ which is such that each sequence in $\tau$ extends \sx.

\noindent
Proof of Claim~2: By induction on \bbb. Suppose the claim holds for \bbb\ \zzz
and that we wish to show that it holds for $\bbb+1$. This is easy but requires \zzz
the essential use of \AC. Suppose $\sx\e\A$ is \good, say $\sx=\seq{\sx_0,\dots,\sx_n}$. \zzz
By \AC, we can find a set $A\e\A$ such that $\ol A\prece_n\A$ and $\ol{\sx_n}\subseteq\ol A$. \zzz
The sequence
\[
\seq{\sx_0,\dots,\sx_n, A}
\]
is also \good, and so we may apply the inductive hypothesis to it to obtain \zzz
the appropriate $\tau$ for \sx.

    At limit \bbb, we in addition need to piece together the $\tau$ of \zzz
appropriate ranks using $\Sx^1_1$-\AC. This is straightforward.

    This concludes our sketch of the proof of the claim and of 6.2.\done

    Turning next to the proof of 6.3, let \A\ be a non-\bbb-model of \T, \zzz
and let $\prec$ be as before. We now consider trees $\tau\e\A$ of finite sequences \zzz
$\sx=\seq{\sx_0,\dots,\sx_{\smash{\lenm\sx}}}$ such that:
\begin{mylist}
\renewcommand{\labelenumi}{\roman{enumi}$'$.}
\item $\ol{\sx_i}$ is finite for all $i\lex\lenm\sx$, and \zzz
\item $\seq{\ol{\sx_i},\ol{\sx_{i+1}},\dots,\ol{\sx_{\smash{\lenm\sx}}}}$ fulfils the $i^\text{th}$
axiom of \T\ for all $i\ltx\lenm\sx$.
\end{mylist}
(More precisely, condition (i$'$) of course means this. \sx\ is a subset of \zzz
\w; $\sx_i$ is the set $\sx_i=\myset{m}{\seq{m,i}\e\sx}$; and $\sx_i$ codes the \zzz
collection $\myset{X_j}{j\e\w}$ where $X_j=\myset{x}{\seq{x,j}\e\sx_i}$. We require in \zzz
(i$'$) that $X_j=\emptyset$ for all but finitely many $j$.)

    As before, it suffices to show that \bfx{Claim~1} holds. For if it does, \zzz
then by overspill there exists such a $\tau\e\A$ of nonstandard rank. Choose \zzz
(externally) an infinite branch $\seq{\sx_0,\sx_1,\dots}$ of $\tau$, and let $\B=\cupx_{i\in\w}\m\ol{\sx_i}$. \zzz
Then $\B\subseteq\A$; \B\ is a model of \T\ by Lemma 1.1; and $\A\ne\B$ as before.

    To show that the new version of \bfx{Claim~1} holds, we need, once again, \zzz
to prove a bit more. Say a sequence $\sx=\seq{\sx_0,\dots,\sx_{\lenm\sx}}$ is \ul{\good} if (i$'$) and \zzz
(ii$'$) hold, and
\begin{mylist}
\renewcommand{\labelenumi}{\roman{enumi}$'$.}
\setcounter{enumi}{2}
\item for all $i\ltx\lenm\sx$, there exist satisfaction functions $f,g,\dots,h$ \zzz
for the $i^\text{th}$ axiom of \T\ such that for all $j$ with $i\lex j\ltx\lenm\sx$,
\[
\ol{\sx_{j+1}}\supseteq f''\m\ol{\sx_j}\cup g''\m\ol{\sx_j}\cup\dots\cup h''\m\ol{\sx_j}\,.
\]
\end{mylist}
Then \bfx{Claim~2} also holds for this new version of $\good$ness, and the proof is \zzz
as before. But now, given a \good\ sequence $\seq{\sx_0,\dots,\sx_{\smash{\lenm\sx}}}$, we do not need \zzz
any strong axioms to extend it to another good sequence
\[
\seq{\sx_0,\dots,\sx_{\smash{\lenm\sx}}\m, A}\,,
\]
for it is sufficient (and, by the definition (i$'$), necessary) to take $\ol A$ \zzz
to be a finite collection of sets.

    Finally, we note the well-known result that $\Sx^1_1$-\BI\ implies $\Sx^1_1$-\AC. \zzz
Thus we have strong enough axioms to carry out the inductive step at limit \zzz
ordinals.

    This concludes the sketch of the proof of \bfx{Claim~2}, and of 6.3.\done

    Next we shall give a more general development of the ideas contained \zzz
in the above proof, starting with an extension of the definition of fulfilment \zzz
to the language \Lw. It will be convenient to suppose that \Lw\ \zzz
is defined so that all formulae have only finitely many free variables. \zzz
We shall consider sequences $\sx=\seq{(A_i,B_i)}_{i\in I}$ of ordered pairs, where \zzz
each $A_i$ is, as before, a subset of the domain of individuals, and each \zzz
$B_i$ is a subset of \Lw. Often we shall write $\sx=\seq{(\sx^0_i,\sx^1_i)}_{i\in I}$. For \zzz
reasons which will become clear, we will henceforth suppose that $A_i$ \zzz
and $B_i$ \ul{are finite for all} $i\e I$. We will also suppose that $\seq{A_i}_{i\in I}$ \zzz
and $\seq{B_i}_{i\in I}$ are increasing chains, with the former closed under the functions \zzz
of \sss L. To the previous definition of fulfilment on page \pageref{Chapter1:DefinitionFulfil}, add the clauses:
\begin{align*}
\ppl\aax\Phi\,\ppr ^\sx_i  &=\Ax j\gex i\,\aax\myset{\px^\sx_j}{\px\e\Phi\cap B_j}\,,\ \text{$j$ a new variable,}\\
\ppl\vvx\Phi\,\ppr ^\sx_i  &=\vvx\myset{\px^\sx_i}{\px\e\Phi\cap B_{i+1}}\,.\\
\ppl\nx\!\aax\Phi\,\ppr^\sx_i &=\ppl\vvx\myset{\nx\px}{\px\e\Phi}\m\m\ppr ^\sx_i\,,\ \text{and}\\
\ppl\nx\!\vvx\Phi\,\ppr^\sx_i &=\ppl\aax\myset{\nx\px}{\px\e\Phi}\m\m\ppr ^\sx_i\,.
\end{align*}
This choice of definition is rather arbitrary, and other choices are possible. \zzz
This is because each conjunction or disjunction can be regarded as having \zzz
a domain unique to itself, whereas all the quantifiers \Ex\ and \Ax\ have a \zzz
common domain. (In this chapter, the words ``conjunction'' and ``disjunction'' \zzz
refer to both the finitary and infinitary kinds.)

    The analogue of Lemma 1.1 from page \pageref{Chapter1:Lemma11} holds; we shall combine parts (i) and (ii) of \zzz
1.1, as this is all we shall require. Recall that $\sx\up\,(n+1)=\seq{\sx_i}_{i\lex n}$.

\thm{6.4 Lemma} \label{Chapter6:Lemma64}Let $\sx=\seq{(A_i,B_i)}_{i\in\w}$ and let \px\ be a sentence of \Lw. \zzz
Suppose that $\cupx_i\m B_i$ includes all subformulae of \px, that $\A=\cupx_i\m A_i$, and \zzz
that $\sx\up\,(n+1)$ fulfils \px\ for all $n\e\w$. Then \px\ is true in \A.

\noindent
Proof: We shall sketch a game-theoretic argument. We can suppose that \px\ \zzz
is in negation-normal-form. The set of subformulae of \px\ may be considered \zzz
as a tree using the subformula relation. Two players, \Ax\ and \Ex, play \zzz
the following game (called the \ul{\px-game}), in the process choosing a branch \zzz
$\seq{\px=\px_0,\px_1,\px_2,\dots,\px_\Omega}$ of this tree and an increasing
sequence $s_0\subseteq s_1\subseteq\dots\subseteq s_\Omega$, \zzz
with each $s_j$ a valuation in $\A$ of the free variables of $\px_j$. \Ex\ \ul{wins} \zzz
iff $\px_\Omega(s_\Omega)$ is true in \A. Suppose $\px_j$ and $s_j$ have been chosen, and \zzz
consider the $(j+1)^\text{st}$ \ul{stage}:
\begin{mylist}
\renewcommand{\labelenumi}{(\alph{enumi})}
\item If $\px_j$ is atomic or the negation of an atomic formula, let $\Omega\eqx j$, \zzz
and the game halts.
\item If $\px_j=\ssf Qx\m\tx$, let $\px_{j+1}=\tx$ and let $s_{j+1}=s_j\cup\{(x,a)\}$, where $a\e\A$ \zzz
is chosen by \ssf Q.
\item If $\px_j=\tx\vx\yx$\, (or $\tx\ax\yx$), let $s_{j+1}=s_j$ and let $\px_{j+1}$ be \tx\ or \yx, \zzz
as chosen by \Ex\ (or \Ax, respectively).
\item If $\px_j=\smash\vvx\Phi$\, (or $\smash\aax\Phi$), let $s_{j+1}=s_j$ and let $\px_{j+1}\e\Phi$ be chosen by \Ex\ \zzz
(or \Ax, respectively).
\end{mylist}
\noindent
Then \px\ is true in \A\ iff \Ex\ has a winning strategy for the \px-game. \zzz
We shall show that \Ex\ in fact has a winning strategy such that for any \zzz
play according to this strategy,
\[
\label{Chapter6:eq:1}\tag{1}
\parbox{3.75in}{if for any $j\lex\Omega$, if \Ax\ has made all of their choices up to the
$(j+1)^\text{st}$ stage from $A_i\cup B_i$, then $(\px_js_j)^{\sx\upharpoonright(n+1)}_i$
holds of all $n>i$.}
\]
We shall simultaneously define this strategy and verify \eqref{Chapter6:eq:1}, using induction \zzz
on $j$. If $j\eqx 0$, \eqref{Chapter6:eq:1} is just the premise of the lemma. Consider the \zzz
$(j+1)^\text{st}$ stage. If it is \Ax's move, then any choice by them will \zzz
preserve \eqref{Chapter6:eq:1}. Suppose it is \Ex's move, that $\px_j=\Ex x\m\px_{j+1}x$, say, and \zzz
that \Ax\ has made all of their choices from $A_i\cup B_i$. Then by the inductive \zzz
hypothesis $(\px_js_j)^{\smash{\sx\upharpoonright(n+1)}}_i$ holds for all $n>i$. As $A_{i+1}$ is finite, there exists \zzz
an $x\e\A_{i+1}$ such that $(\px_{j+1}(x,s_j))^{\smash{\sx\upharpoonright(n+1)}}_i$ holds for infinitely many (and hence all) \zzz
$n>i$; let \Ex\ choose such an $x$. Cases (c) and (d) are similar.\done

    In this chapter, the base theory \s\ will be assumed to include \zzz
$\Sx_1$-\Collection\ and $\Dx_0$-\Separation. The reader is reminded of the simple \zzz
notions and conventions regarding non-standard models given on page \pageref{Introduction:BasicConventions}. \zzz
We shall be considering non-\bbb-models $\A=\seq{A,\dots}$ of \s. In each such \zzz
structure we shall (usually implicitly) fix some linear ordering $\prec\!\e\,A$ \zzz
which is well-founded internally but not externally, and let $\aaa,\bbb,\dots$ \zzz
range over its domain. \aaa\ is a \ul{standard ordinal} of \A\ if it is in the \zzz
well-founded part of this linear ordering. \ul{The} standard ordinal of \A\ is \zzz
the order-type of its standard ordinals; because $\A\modelx\Sx_1$-\Replacement, this \zzz
is independent of the choice of $\prec$. As usual, it will be convenient to use \zzz
the notations and conventions associated with von Neumann ordinals.

    We shall consider trees $\tau$ of finite sequences \sx. Say a tree $\tau$ \ul{fulfils} \zzz
a sentence \px\ of \Lw\ if \sx\ fulfils \px\ for each $\sx\e\tau$. Given \aaa\ \zzz
in the domain of $\prec$, say $\tau$ is of \ul{rank} \aaa\ if there is a map \ul{\rk} in $A$ \zzz
from $\tau$ to $\aaa+1=\myset{\bbb}{\bbb\prece\aaa}$ satisfying
\[
\rk(\sx)=\sup\myset{\rk(\sx')+1}{\sx'\e\tau\text{ and $\sx'$ extends \sx}}\,.
\]
with $\rk(\text{root of $\tau$})=\aaa$; we shall write $\rk(\tau)=\aaa$.

    The next result is the analogue of Lemma 1.2 on page \pageref{Chapter1:Lemma12}.

\thm{6.5 Lemma} Let \px\ be a sentence of \Lw\ and let \A\ be a model of \s\ \zzz
with $\px\e A$ and \px\ true in \A. Then for each standard ordinal \aaa\ \zzz
of \A\ there exists a tree in \A\ of rank \aaa\ which fulfils \px, and if \zzz
we have sufficient \ssf{Bar-Induction}, the construction can be formalized.

\noindent
Proof: In Lemma 1.2, we considered finite pieces of satisfaction functions \zzz
for \px. The infinitary analogue of satisfaction functions is a winning \zzz
strategy for the \px-game, and we shall consider finite pieces of the same.

    A (\ul{winning}) \ul{strategy} for \Ex\ is a partial function \bb S from $(A\cup\sub(\px))^{<\w}$ \zzz
into $A\cup\sub(\px)$ such that \Ex\ plays (and wins, respectively) the \px-game by choosing \zzz
at each turn \bb S\xvect a, where \vect a is the listing of \Ax's choices up \zzz
to that point. If \bb S is a strategy for \Ex, let $\bb S_0$ be \bb S with its \zzz
domain restricted to arguments \vect a for which \Ex\ actually needs the value \zzz
\bb S\xvect a in the course of some play; $\bb S_0$ is a \ul{minimal} strategy. Say \zzz
$\sx=\seq{(A_i,B_i)}_{i\le\lenm\sx}$ is \ul{closed} under a strategy \bb S for \Ex\ if
\[
\bb S_0''(\A_i\cup B_i)^{<\w}\subseteq A_{i+1}\cup B_{i+1}
\]
for all $i\ltx\lenm\sx$.

\thm{Claim~1} Given finite sets $A_0$, $B_0$, a function $F\mymap\w\Lw$, and a \zzz
strategy \bb S for \Ex, there exists $\sx=\seq{(A_i,B_i)}_{i<\w}$ which is closed \zzz
under \bb S and the functions of \sss L, and satisfies $Fi\e B_{i+1}$ for all $i\ltx\w$.

\noindent
This claim is not quite obvious, because $\bb S_0$ may contain arbitrarily \zzz
large \n-tuples in its domain. Suppose $A_i$ and $B_i$ are finite. Let $C$ \zzz
be the smallest set of formulae including $B_i$ and containing \yx\ and $\yx'$ if \zzz
it contains $\Ex x\m\yx$, $\Ax x\m\yx$, $\yx\vx\yx'$, or $\yx\ax\yx'$. By K\"onig's lemma, $C$ is \zzz
finite, and the lengths of the sequences \vect a from $A_i\cup B_i$ which must  \zzz
considered may be bounded by the cardinality of $C$. This establishes the \zzz
claim.

    Let $\sx=\seq{(A_i,B_i)}_{i\le\n}$. The \ul{\sx-game} (\ul{for} \px) is like the \zzz
\px-game, except that \Ax's choices are restricted to $A_{n-1}\cup B_{n-1}$, and \zzz
at each turn of \Ex, if \Ax\ has made all of their previous choices from \zzz
$A_i\cup B_i$, then \Ex\ must choose from $A_{i+1}\cup B_{i+1}$; and if either \Ax\ or \Ex\ \zzz
cannot make a choice, the game is drawn. A strategy for \Ex\ for the \zzz
\sx-game is \ul{good} if for each play according to this strategy, $\px_i s_i$ is true \zzz
for all $i\lex\Omega$, where $\px_i$, $s_i$, and $\Omega$ are as in 6.4 on page \pageref{Chapter6:Lemma64}. The following is \zzz
an easy consequence of these definitions.

\thm{Claim~2} If \bb S is a winning strategy for \Ex\ for the \px-game and if \sx\ \zzz
is closed under \bb S, then \bb S is a \good\ strategy for the \sx-game. \zzz
Conversely, any \good, minimal strategy for a \sx-game may be extended to \zzz
a winning strategy \bb S for \Ex\ for the \px-game which is such that \sx\ is \zzz
closed under \bb S.

\noindent
(We may now see the reasons for the above definitions; the notion of \zzz
\ul{a winning strategy for \Ex\ for the \px-game} cannot be expressed directly \zzz
in our language, while the notion of \ul{a \good\ strategy for \Ex\ for the} \zzz
\ul{\sx-game} can be so expressed by the use of partial truth definitions \zzz
because such strategies are essentially finite; this will be required \zzz
when we give an internal version of 6.5.)

    Given a sequence \sx\ and a tree $\tau$, let $\sx_\smile\tau$ be the tree
$\myset{\sx_\smile\sx'}{\sx'\e\tau}$, where the latter $_\smile$ is sequence concatenation.

\thm{Claim~3} For all \sx\ for which there exists a \good\ strategy for \Ex\ for \zzz
the \sx-game of \px, and for all standard ordinals \aaa\ of \A, there \zzz
exists a tree $\tau$ of rank \aaa\ such that $\sx_\smile\tau$ fulfils \px.

\noindent
The lemma follows immediately from this claim. We use induction on \aaa. \zzz
If $\aaa\eqx 0$, take $\tau$ to be empty. Suppose the claim holds for all $\aaa\ltx\bbb$. \zzz
Let \sx\ satisfy the premise. By \bfx{Claim~2}, there exists $\sx'=\sx_\smile\seq{(A',B')}$ \zzz
for which there exists a good strategy for the $\sx'$-game. By the inductive \zzz
hypothesis, for all $\aaa\ltx\bbb$ there exists a $\tau$ of rank \aaa\ such \zzz
that $\sx'_\smile\tau$ fulfils \px. By $\Sx_1$-\ssf{Strong Replacement}, there exists a \zzz
set $K$ such that for all $\tau'\e K$, $\sx_\smile\tau'$ fulfils \px, and for all $\aaa\ltx\bbb$ \zzz
there exists $\tau'\e K$ of rank $\gex\,\aaa$. (If \bbb\ is a successor, $K$ may \zzz
be taken to be a singleton.) Then $\tau=\seq{(A',B')}_\smile\cupx K$ is the required \zzz
tree.\done

    We shall want to consider an extension of the formula hierarchy into \zzz
the transfinite. The details here do not matter; all we require of, say, \zzz
$\Sx_\aaa$ is that it be closed under finitary conjunction, finitary disjunction, \zzz
and existential quantification (and, if $\aaa\eqx 0$, infinitary conjunctions and \zzz
disjunctions), and that there exists a suitable universal formula, e.g.\ a \zzz
formula $\px_\aaa(e,x)\,\in\,L_{\aaa^+}\cap\Sx_\aaa$ (where $\aaa^+=$ the next admissible) such that \zzz
in any model \A\ of \s, for all $\Sx_\aaa$ formulae $\yx x$
in \A, $\A\modelx\Ax x\,\ppl \yx x\equiv\px_\aaa(\godel\yx,x)\m\ppr $. \zzz
For any subclass $\Gamma$ of \Lw, let $\Gamma^{\bb B}=\Gamma\cap\bb B$. Let $\Sx^{\smash{(\bb B)}}_{(\aaa)}$ be the closure \zzz
of $\Sx^{\bb B}_\aaa$ under conjunction, disjunction, existential and bounded universal \zzz
quantification. We note that in the presence of the schema of
\begin{flalign*}
&\Sx_\aaa\HY\ssf{(Strong) Replacement}\!: &&\Ax x\e a\,\Ex y\,\tx\imp\Ex f\,
\Ax x\e a\,(f\mx x\ne\emptyset\ax\Ax y\e f\mx x\m\,\tx)\,,\ \tx\e\Sx_\aaa\,, &&&
\end{flalign*}
every formula of the form $\Ax x\ltx y\,\tx$ with $\tx\e\Sx_\aaa$ is equivalent to a $\Sx_\aaa$ \zzz
formula.

    \refstepcounter{thisisdumb}\label{Chapter6:6.5.ii}As a corollary of the above proof, we have the following extension.

\thm{6.5(ii) Corollary} Let \A\ be a model of \s, and let \T\ and \bb B be coded in \zzz
\A, with \T\ a set of sentences of \Lw\ and \bb B a transitive set satisfying \zzz
some weak closure condition (e.g.\ $\sss L_{\bb B}$ is a fragment of \Lw). Let $F\e A$ \zzz
be, in \A, a map from \w\ into \Lw, and suppose its range includes $\Sx^{\bb B}_\aaa$ \zzz
and all subformulae of \T. For an ordinal \bbb\ of \A, say
\[
\tau\ \ \text{\bbb-fulfils}\ \ \T + \Tr\Sx^{\bb B}_\aaa\ \ \text{(with respect to $F$)}
\]
if $\tau\e A$ is a tree of rank $\gex\,\bbb$, and for all $\sx\e\tau$, all $i\ltx\lenm\sx$,
\begin{mylist}
\itemsep0ex
\item $Fi\e\sx^1_i\,$
\item for all $\px\e\sx^1_i\cap\T$, $\,\px^\sx_i$, and
\item for all $\px\xvect x\m\e\m\sx^1_i\cap\Sx_\aaa\cap\bb B$, all $\vect x\ve\sx^0_i$, $\,\px\xvect x\imp(\px\xvect x)^\sx_i$.
\end{mylist}
\mypagebreak{[2]}
Say $\tau$ \bbb-fulfils $\dots$ if it does so w.r.t.\ some such $F$.

    Then if \A\ is a model of \T\ + $\Sx^{\bb B}_{\aaa+1}$-\Collection, then for each standard ordinal \zzz
\bbb\ of \A, there exists a tree $\tau$ in \A\ which \bbb-fulfils \T\ + $\Tr\Sx^{\bb B}_\aaa$.\done

   \refstepcounter{thisisdumb}\label{Chapter6:6.6}
    Our next theorem is the generalization of 6.3 and the analogue of 5.3 from page \pageref{Chapter5:Theorem53}.

\thm{6.6 Theorem} Let \A\ be a model of \s; let $\sss L_{\bb B}$ be a fragment of \Lw\ \zzz
coded in \A\ by a code which is countable in \A, and let \T\ be a theory \zzz
of $\sss L_{\bb B}$ coded in \A. Let \aaa\ be a standard ordinal of \A\ and suppose \zzz
that $\Sx^{\bb B}_{\aaa+1}$-overspill holds in \A. Consider the set
\[
\sss S=\Myset{\B\prece^\text{end}_{\Sx^{\bb B}_\aaa}\A}{\B\modelx\T}
\]
\begin{mylist}
\vspace{-3.5ex}
\renewcommand{\labelenumi}{(\alph{enumi})}
\item Suppose that \aaa\ is \ul{finite}, that \T\ is included in the closure of \zzz
the set of finitary formula under conjunction, disjunction, and existential \zzz
quantification, and that \T\ extends \s\ plus the first-order schema of \zzz
\Collection. Then
  \begin{mylist}
  \renewcommand{\labelenumii}{\textbullet}
  \vspace{.8ex}
  \itemsep0ex
  \item \sss S has an indicator.
  \item $\sss S\setminus\{\A\}$ is non-empty iff \A\ is a model of the $\Sx^{\smash{(\HF)}}_{(\aaa+1)}$ theory of \T, and
  \item \sss S is cofinal in \A\ iff \A\ is a model of the $\Pi^{\smash{(\HF)}}_{(\aaa+2)}$ theory of \T.
  \end{mylist}
\item Suppose \A\ is locally countable and either that \A\ is a model of \zzz
$\Sx^{\bb B}_{\aaa+1}$\ssf{-Strong-Replacement} or that \T\ extends \s\ + $\Sx^{\bb B}_{\aaa+1}$\ssf{-Strong-Replacement}. Then
  \begin{mylist}
  \renewcommand{\labelenumii}{\textbullet}
  \vspace{.8ex}
  \itemsep0ex
  \item \sss S has an indicator.
  \item $\sss S\setminus\{\A\}$ is non-empty iff \A\ is a model of the $\Sx^{\smash{(\bb B)}}_{(\aaa+1)}$ theory of \T, and
  \item \sss S is cofinal in \A\ iff \A\ is a model of the $\Pi^{\smash{(\bb B)}}_{(\aaa+2)}$ theory of \T.
  \item If \bb B is a resolvable admissible set with height equal to the ordinal of \A, then \sss S is symbiotic with
  \[
  \Myset{\B\e\sss S}{\B\,\text{ \bb B-saturated}}\,,
  \]
  that is, those \B\ which realize a \bb B-r.e.\ type (which may have parameters from \B) if every \bb B-finite subset of the type is realized.
  \end{mylist}
\end{mylist}

\noindent
Proof: We shall only consider the condition for $\sss S\setminus\{\A\}$ being nonempty, \zzz
and the symbiosis. The remainder of the theorem will then be clear from \zzz
our constructions.

    The necessity of the condition is clear, so let us first consider its \zzz
sufficiency in (a) when $0\ne\aaa$. We need one preliminary definition: for each \zzz
(real) set $a$, define the $\Dx_0$ formula $\tx_a(x)$ by
\[
\tx_a(x)=\Ax y\e x\vvx_{b\in a}\m\tx_b(y)\ax\aax_{b\in a}\m\Ex y\e x\,\tx_b(y)\,.
\]
The formal sentence $\Ex x\m\tx_a(x)$ will often be paraphrased by \ul{$a$ exists}.

    Let \T\ be as given in (a). Then by \hyperref[Chapter6:6.5.ii]{6.5.ii} for each standard ordinal \bbb\ \zzz
of \A, \T\ proves
\begin{mylist}
\renewcommand{\labelenumi}{(\arabic{enumi})}
\itemsep0ex
\item\quad if \T\ exists and is hereditarily countable and if \bbb\ exists,
\item\quad then there exists $\tau$ which \bbb-fulfils \T\ + $\Tr\Sx_\aaa$.
\end{mylist}
(If $\aaa\eqx 0$, this cannot be expressed as a $\Sx_1$ sentence, and so we shall \zzz
postpone this case until later.) If \A\ is a model of the $\Sx^{\smash{(H\!F)}}_{(\aaa+1)}$ consequences \zzz
of \T, we have that (2) holds in \A\ for all standard \bbb. By overspill \zzz
we may find a witness $\tau$ of nonstandard rank. Choose any infinite branch $\sx=\seq{\sx_i}_{i\in\w}$ \zzz
of $\tau$, and let $\sss C=\A\up\cupx_{i\in\w}\m\sx^0_i$. \sss C is a $\Sx^\HF_\aaa$-elementary \zzz
substructure of \A\ which is a model of \T. Let \sss B be the initial \zzz
segment of \A\ determined by \sss C. Then \sss B is as required, for Lemma \hyperref[Chapter6:Lemma64]{6.4} \zzz
also shows that each axiom of \T\ is absolute between \sss B and \sss C. Finally \zzz
$\sss A\ne\sss B$ as follows. We have\mypagebreak{[2]}

\[
\A\modelx\Ax n\,\Ex!\seq{b_0,\dots,b_n}\,\ppl b_0=\cupx_{\sx\in\tau}\m\m\sx^0_{\smash{\lenm\sx}}\ax\Ax i\ltx n\,(b_{i+1}=\cupx b_i)\m\ppr \,.
\]
The relation $x=\cupx y$ is $\Dx_0$ and so by $\Sx_1$-\Collection\ we \zzz
have that $b=\cupx_{n\in\w}\m b_n$ is an element of \A. Let $c=\myset{x\e b}{x\notE x}$. \zzz
Then $c\notE b$, and as $\sss B\subseteq b$, we have $c\notE\sss B$. This proves (a) for $\aaa\ne 0$.

    Before we proceed with (b), let us consider the easier problem of \zzz
finding a $\Sx^{\bb B}_\aaa$-elementary \ul{substructure} of \A\ which is a model of \T\ \zzz
when $\Sx^{\bb B}_{\aaa+1}$-\Replacement\ holds in \A\ and \T\ is arbitrary. (We are not \zzz
assuming \T\ has any coding ability.) We shall show that this exists if \zzz
(and only if) \A\ is a model of those sentences contained in the closure \zzz
of $\Sx^{\bb B}_{\aaa+1}$, under conjunction, disjunction and existential quantification \zzz
which are implied by \T. Let \bbb\ be a standard ordinal of \A\ and \zzz
define the following sets in \A:
\[
\bb B_\bbb=\myset{x\e\bb B}{\rk(x)\ltx\bbb}\,,\quad \T_\bbb=\T\cap\bb B_\bbb\,.
\]
\refstepcounter{thisisdumb}Let $G\e A$ be a surjection from \w\ onto $\sss L_{\bb B_\bbb}$. For each $n\e\w$ and each \zzz
standard ordinal $\gamma$ of \A, define the formula (where the middle terms refer to \zzz
i,ii in the definition of \bbb-fulfil on page \pageref{Chapter6:6.5.ii})
\[
F_{\gamma,n}(\sx_0,\dots,\sx_{n-1})=\ppl \aax_{i\in n}Gi\e\sx^1_i\ppr \ax\T^{\seq{\sx_0,\dots,\sx_{n-1}}}_\bbb\ax
\big(\Tr\Sx_\aaa\big)^{\seq{\sx_0,\dots,\sx_{n-1}}}\ax\Ex\sx_n\aax_{\delta<\gamma}F_{\delta,n+1}(\sx_0,\dots,\sx_n)
\]
where this is assumed to be defined without any coding as in our proof of \zzz
the Completeness Theorem on page \pageref{Chapter1:GodelsProof} and of \hyperref[Chapter5:5.3.ii]{5.3.ii} on page \pageref{Chapter5:5.3.ii}. Then for all $\gamma,\,n$, \T\ proves
\[
\Ex\sx_0,\dots,\sx_{n-1}\,F_{\gamma,n}(\sx_0,\dots,\sx_{n-1})
\]
since we can see (nonconstructively) that this is true in any model of \T. \zzz
Now suppose $F_{\gamma,n}(\sx_0,\dots,\sx_{n-1})$ holds in \A\ for some $\sx=\seq{\sx_0,\dots,\sx_{n-1}}\e A$. We \zzz
claim that there exists in \A\ a tree $\tau$ of rank $\gamma$ such that $\sx_\smile\tau$ fulfils \zzz
$\T_\bbb+\Tr\Sx^{\smash{\bb B_\bbb}}_\aaa$ (with respect to $G$). We use induction on $\gamma$. By the inductive \zzz
hypothesis, there exists $\sx_n\e A$ so that for all $\delta\ltx\gamma$ there exists $\tau^\delta$ of rank $\gex\,\delta$ so \zzz
that $\sx_\smile\seq{\sx_n}_\smile\tau^\delta$ fulfils $\T_\bbb+\Tr\Sx^{\smash{\bb B_\bbb}}_\aaa$ (w.r.t.\ $G$).
By $\Sx^{\bb B}_{\aaa+1}$-\Replacement, we \zzz
may form a collection $K$ of such $\tau^\delta$'s. Then $\tau=\seq{\sx_n}_\smile\!\cupx K$ is as \zzz
required. Hence if \A\ is a model of the part of the theory of \T\ specified \zzz
above, by setting $\gamma\eqx \bbb$ and applying overspill and 6.4, we may obtain \zzz
a proper substructure of \A\ which is a model of \T.

    We may use this construction for (a) for the case $\aaa\eqx 0$, by taking \zzz
the appropriate initial segment.

    Next consider the case where \T\ extends \s\ + $\Sx^{\bb B}_{\aaa+1}$-\Replacement\ \zzz
where \A\ is a locally countable model of the $\Sx^{\smash{(\bb B)}}_{(\aaa+1)}$ theory of \T, and we \zzz
wish to find a $\Sx^{\bb B}_\aaa$-elementary initial segment of \A\ which is a model of \T. \zzz
We need a definition analogous to that given in the non-\w-model case, \zzz
one which captures some of the properties of trees defined non-constructively \zzz
via satisfaction functions. These further properties must be local in \zzz
that they are preserved under the union of trees. Consider the following. \zzz
Let \bbb\ be a standard ordinal of \A. Fix some surjection $G\mymap\w{\sss L_{\bb B_\bbb}}$. \zzz
For two sequences, $\sx$, $\sx'$, with $\lenm\sx=\lenm{\sx'}=n$, say $\sx\prece\sx'$ if for all \zzz
$i\ltx\n$, $\sx^0_i\subseteq\sx'\m^0_i$ and if $Gi$ has free variables \vect x, then for all $\vect x\ve\sx_i$, \zzz
$(Gi)(\vect x)^\sx_i\equiv (Gi)(\vect x)^{\sx'}_i$. For any ordinal $\gamma$ say \ul{extendible$_{\bbb,\gamma}(\tau)$} (\ul{w.r.t.} $G$) \zzz
if $\tau$ fulfils $\T_\bbb+\Tr\Sx^{\bb B_\bbb}_\aaa$, \ul{and} for all $\delta\ltx\gamma$, $\sx\e\tau$, $x\e\sx^0_{\smash{\lenm\sx}}$, $y\e x$,
\begin{mylist}
\renewcommand{\labelenumi}{\textbullet}
\itemsep0ex
\item  there exists a tree $\tau'$ of sequences each extending, or a subsequence of, \sx, and
\item there exists a tree isomorphism $f$ from $\tau\up\,\sx=\myset{s\e\tau}{s\supseteq\sx\vx s\subseteq\sx}$
to $\tau'$ such that for all $s\e\tau$, $s\prece f\mx s$,\footnote{It is sufficient to suppose
$f\mx\sx=\sx$ rather than $s\prece f\mx s$ for all $s\e\tau$; this latter requirement
is needed, however, to construct a \bb B-saturated initial segment.} and
\item for all $s\e\tau'$, if $\lenm s\gtx\lenm\sx$ then $y\e s^0_{|\sx|}$, and finally
\item  extendible$_{\bbb,\delta}(\tau')$ holds.
\end{mylist}

    Suppose \T\ includes $\Sx^{\bb B}_{\aaa+1}$-\Replacement. Then by the proof of 6.5, \zzz
it is easy to see by induction on $\gamma$ that in any model of \T\ containing \zzz
\bbb, $\gamma$, $\T_\bbb$, and $\bb B_\bbb$, for any ``good'' sequence \sx\ there is a tree $\tau$ of rank \zzz
$\gamma$ such that $\sx_\smile\tau$ is a $(\bbb,\gamma)$-extendible tree, and so in particular,
\[
\text{if \bbb, $\T_\bbb$, and $\bb B_\bbb$ exist, there is a  $(\bbb,\bbb)$-extendible tree of rank \bbb}\,.\label{Chapter6:eq:3}\tag{3}
\]
If $\aaa\ne 0$, \eqref{Chapter6:eq:3} may be expressed in \T\ by a $\Sx^{\smash{(\bb B)}}_{(\aaa+1)}$ sentence. So if \A\ is a \zzz
model of the $\Sx^{\smash{(\bb B)}}_{(\aaa+1)}$ consequences of \T, then by overspill there exists such \zzz
a tree $\tau$ of nonstandard rank. We now construct our initial segment. \zzz
Let $\tau\eqx \tau^0$, and choose $\sx^0\e\tau^0$ of nonstandard rank. Pick $x\e(\sx^0)^0_{\smash{{\lenm{\sx^0}}}}$ and \zzz
$y\e x$. Then there is a tree $\tau^1$ satisfying the four points in the definition of extendible \zzz
which is $(\gamma,\gamma)$-extendible of some nonstandard $\gamma$. Choose $\sx^1\e\tau^1$ of \zzz
nonstandard rank and $x\e(\sx^1)^0_{\smash{\lenm{\sx^1}}}$ and $y\e x$. Continuing in this way we have an \zzz
increasing sequence of finite sequences
\[
\sx^0\subseteq\sx^1\subseteq\sx^2\subseteq\dots\,.
\]
Let $\sx=\seq{\sx_i}$ be the limit. Then \sx\ fulfils \T\ + $\Tr\Sx^{\bb B}_\aaa$, and so $\B=\A\,\up\cupx_{i\in\w}\m\sx^0_i$ \zzz
is a $\Sx^{\bb B}_\aaa$-elementary substructure of \A\ which is a model of \T. By \zzz
judicious choices of the $x$'s and $y$'s, and using the local countability \zzz
of \A, we may ensure that \B\ is an initial segment. Finally, because \zzz
\T\ includes $\Sx^{\bb B}_{\aaa+1}$-\Collection, we can specify that all the $\tau$'s are \zzz
contained in some $b\e A$, and in this way ensure that $\B\ne\A$.

    If $\aaa\eqx 0$, or if \T\ does not necessarily include $\Sx^{\bb B}_{\aaa+1}$-\Replacement\ but \zzz
\A\ is a model of $\Sx^{\bb B}_{\aaa+1}$-\Replacement, then we can combine the above two \zzz
constructions. That is, we can express \eqref{Chapter6:eq:3} directly without coding, and \zzz
so expressed, it is provable in \T. So if \A\ is a model of the $\Sx^{\smash{(\bb B)}}_{(\aaa+1)}$-consequences  \zzz
of \T, then \eqref{Chapter6:eq:3} holds in \A, and moreover we may express \eqref{Chapter6:eq:3} \zzz
uniformly in \A\ using $\Sx^{\bb B}_{\aaa+1}$-\Replacement\ in \A\ and codes of \T\ and \bb B. \zzz
Now use overspill and construct the initial segment as before.

    Next let us consider the following problem. Suppose \A\ is a model \zzz
of \T\ + $\Sx^{\bb B}_{\aaa+1}$-\Replacement\ + $\Sx^{\bb B}_{\aaa+1}$-overspill which codes \T\ and a resolvable countable \zzz
admissible set \bb B, and that we wish to find a proper $\Sx^{\bb B}_\aaa$-elementary \zzz
\ul{substructure} of \A\ which is a \bb B-saturated model of \T. We may suppose \zzz
that \T\ contains $\Av x\,(\tx\vx\nx\tx)$ for all $\tx\xvect x\ve\sss L_{\bb B}$. Let $R\mymap{\mathrm{ord}(\bb B)}{\bb B}$ \zzz
be a resolution of \bb B. Let $\tau\e A$ be any tree which fulfils \T\ + $\Tr\Sx^{\bb B}_\aaa$ of nonstandard rank \zzz
and let $C=\cupx_{\sx\in\tau}\m\sx^0_{\smash{\lenm\sx}}$. Fix some (external) listing of all pairs $\seq{\m\Phi,\xvect a}$, \zzz
where $\Phi=\Phi(x,\vect a)$ is a \bb B-recursive subset of $\sss L_{\bb B}$ with only the free variables \zzz
listed and $\vect a\ve C$. Choose an infinite branch of $\tau$ in the following manner. \zzz
Suppose we have chosen $\sx\e\tau$ nonstandard rank. Let $\seq{\m\Phi,\xvect a}$ be the \zzz
first pair in our listing which we have not yet considered such that $\vect a\ve\sx^0_{\smash{\lenm\sx}}$. \zzz
Suppose
\[
\parbox{4.125in}{for all standard \bbb, there exists $i\e\w$ and an extension $\sx'\e\tau$ of \sx\
of rank $\gtx\bbb$ such that for all extensions $\sx''\e\tau$ of $\sx'$, $\ppl \Ex x\aax(\Phi\cap R\bbb)\m\ppr ^{\sx''}_i$.}
\label{Chapter6:eq:4}\tag{4}
\]
Then by overspill, we may choose $i\e\w$ and a $\sx'$ of nonstandard rank, \zzz
and let this $\sx'$ be in the branch. If not \eqref{Chapter6:eq:4}, then choose $\sx'\e\tau$ to be \zzz
any extension of \sx\ of nonstandard rank.

    Let $\seq{\sx_i}_{i\in\m\w}$ be any infinite branch obtained in this way, and let \zzz
$\B=\cupx_{i\in\w}\m\sx^0_i$. Then \B\ is \bb B-saturated, for let $\Phi(x,\xvect a)$ be any \zzz
\bb B-recursive type. If \eqref{Chapter6:eq:4} held when we came to consider the pair $\seq{\m\Phi,\xvect x}$, \zzz
then $\Phi$ is clearly realizable in \B. Suppose \eqref{Chapter6:eq:4} did not hold, and \zzz
at that point we had already chosen $\sx=\seq{\sx_i}_{i\le\n}$. Then there exists a \zzz
standard \bbb\ such that for all $i\e\w$, all extensions $\sx'\e\tau$ of \sx\ of \zzz
rank $\gtx\bbb$, there is an extension $\sx''$ of $\sx'$ that $\nx\tx^{\sx''}_i$ holds, where \zzz
$\tx=\Ex x\aax(\Phi\cap R\bbb)$. Choose $i$ so that the sentence $\Av a\,(\tx\vx\nx\tx)$ is the \zzz
$i^\text{th}$ element of our implicit listing of \T. Then $(\nx\tx)^{\sx''}_i$, and so $(\nx\tx)^{\sx'}_i$. \zzz
Thus for all extensions $\sx'\e\tau$ of nonstandard rank, $(\nx\tx)^{\sx'}_i$, and so \tx\ \zzz
is false in \B, that is, $\Phi$ is not \bb B-finitely realizable.

    To obtain a \bb B-saturated \ul{initial segment}, we in addition ensure that \zzz
the tree $\tau$ be $(\bbb,\bbb)$-extendible for some nonstandard \bbb, and then in \zzz
choosing an infinite branch, we alternate the two constructions given \zzz
above.

    If \T\ does not necessarily extend \s\ + $\Sx^{\bb B}_{\aaa+1}$-\Replacement\ but \A\ is \zzz
a model of $\Sx^{\bb B}_{\aaa+1}$-\Replacement, we combine all three of the above techniques \zzz
to obtain the required result.

    This completes our proof of 6.6.\done

    As a corollary we have:

\thm{6.7 Corollary} Let \T\ be any r.e.\ extension of $\ZF^-$. \zzz
For any $k\e\w$, each nonstandard model of \T\ has a proper \zzz
$k$-elementary initial segment which is a model of \T.\done

\noindent
A weaker version of 6.7 (requiring the power-set axiom) appeared in the \zzz
unpublished Friedman\,\cite{Frie7y} (but 6.7 was obtained independently).

   \refstepcounter{thisisdumb}\label{Chapter6:6.8}
    Another corollary concerns \w-models of analysis.

\thm{6.8 Corollary} \label{Chapter6:Corollary68}Let \A\ be an \w-model of $\Sx^1_1$-\AC. Then the following \zzz
are equivalent. (Also, see addendum on page \pageref{Chapter7:Chapter6Addendum}.)
\begin{mylist}
\item Let $\T=\myset{\px x}{\px\e\T}$ be a $\Pi^1_1$ collection of formulae of analysis, where \zzz
the $\Pi^1_1$ definition may have parameters from \A. Let $Z\e\A$ and suppose \zzz
$\A\modelx\px Z$ for all $\px\e\T$. Then there exists $B\e\A$ and $\B\subseteq\myset{(B)_n}{n\e\w}$ \zzz
such that $Z\e\B$ and $\B\modelx\px Z$ for all \px\e\T.
\item The set $W=\myset{X\e\A}{\text{$\prec_X$ well-founded}}$ is not $\Sx^1_1/\A$.
\end{mylist}
Proof: Suppose (ii). If \A\ is a \bbb-model,then (i) holds, so suppose \zzz
\A\ is not a \bbb-model. Then there is some pseudo-well-ordering $\prec$ in \A, \zzz
and, moreover, we have $\prec$-overspill holding on this ordering, for otherwise \zzz
we would have that $W$ is $\Sx^1_1/\A$.

    Let $\T=\myset{\px}{\Ax f\,\Ex n\,\tx\m(\m\ol{fn},\godel\px)}$ be any $\Pi^1_1$ set of sentences true in \zzz
\A\ (where we are suppressing the parameter $Z$ and the parameters in \tx). \zzz
Let $\prec_\px$ be the Kleene-Brouwer ordering on the non-past-secured sequences of \zzz
of $\myset{s}{\tx\m(\m s,\m\godel\px)}$. Now by 6.5, for each `ordinal' \aaa\ in the standard \zzz
part of $\prec_\px$, we have
\[
\Ex\tau\,\ppl \rk(\tau)\gex\aaa\ax\Ax\px\,(\rk(\prec_\px)\ltx\aaa\imp\tau\text{ fulfils }\px)\m\ppr \,.
\]
By overspill, this holds for some non-standard \aaa, and we may obtain the \zzz
required \B\ from any infinite branch of any witness $\tau$.

    Now suppose (i) but not (ii). Then \A\ cannot be a \bbb-model. Fix \zzz
some pseudo-well-ordering $\prec\,\in\A$ and suppose
\[
\text{the well-founded part of}\m\prec\,\,\,=\myset{n}{\A\modelx\yx\ol nZ}
\]
for some $\Sx^1_1$ formula \yx\ with parameter $Z\e\A$. Now
\[
\T=\myset{\yx\ol nx}{\n\e\text{well-founded part of}\prec}
\]
is a $\Pi^1_1$ set with parameter $\prec\,\in\A$, and $\T(Z)$ is true in \A. So by (i) \zzz
there exists $X\e\A$ such that for all $\yx\ol nx\e\T$, $(\yx nZ)^X$. But then \zzz
the well-founded part of $\prec$ is $\Dx^1_1$ in \A, and as $\Dx^1_1$-\CA\ holds in \A, \zzz
this contradicts that $\prec$ is a well-ordering in \A.\done

    We were careful enough in our proof of 6.5 so that it may be readily \zzz
formalized to yield:

\thm{6.9 Corollary} For all $k>1$, and for each $\px\e\Sx_{k+2}$,
\[
\s+\Pi_k\HY\BI\provex\Ax\prec\,
\ppl \wf(\prec)\imp\Ax x\,\big(\px x\imp\Ex\tau\,(\rk(\tau)=\rk(\prec)\ax\tau\text{ fulfils } \px x)\m\big)\m\ppr .
\]
For $k\eqx 1$, we require, say, $\Pi_2$-\BI\ or $\Sx_1$-\DC. (In analysis we note that $\Pi^1_1\HY\BI\up\equiv\Sx^1_1\HY\DC\up$
by S.D. Friedman\,\cite{FrieSD79}.)\done

    \refstepcounter{thisisdumb}\label{Chapter6:6.10}Consider the schema
\begin{flalign*}
&\w\HY\RFN\!: &&\Ax X\subseteq\w\,\ppl \px X\imp\Ex\,\text{\w-model of $\px X$ containing $X$}\ppr \,. &&&
\end{flalign*}
Our next result is due to H. Friedman\,\cite{Frie75}. His argument, using the \zzz
completeness of the cut-free sequent calculus for \w-logic, is probably simpler \zzz
than ours. We note the next result immediately gives the \bfx{Fact} on page \pageref{Chapter3:FactDeferredProof}.

\thm{6.10 Corollary} $\s\provex\Sx_{k+2}\HY\w\RFN\equiv\Pi_k\HY\cBI$, for all $k\ge 1$.

\noindent
Proof; We shall first show the right to left implication. Consider the \zzz
version of 6.9 for $\frac12\ast$fulfilment, namely for $\px\e\Sx_{k+2}$
\[
\s+\Pi_k\HY\cBI\provex\Ax\prec\,\subseteq\w^2\,
\ppl \px\imp\Ex\tau\,\big(\rk(\tau)\eqx\rk(\prec)\m\ax\Ax\sx\e\tau\,\px^{\frac12\ast\sx}\big)\m\ppr .
\label{Chapter6:eq:5}\tag{5}
\]
By $\Pi^0_1\HY\CA\up$ we can consider the set:
\[
Y=\myset{\sx\e\w}{\text{\sx\ codes a finite sequence which $\textstyle\frac12\ast$fulfils \px}}\,.
\]
$Y$ is \ul{not} well-founded, for otherwise by \eqref{Chapter6:eq:5} we may find a tree $\tau\subseteq Y$ of \zzz
rank greater than that of $Y$. Since $Y\subseteq\w$, we may choose an infinite \zzz
branch internally; the union of this branch is the required \w-model. \zzz
This argument relativizes to an arbitrary $X\subseteq\w$ by adding a constant \ol X \zzz
and axioms $\myset{\ol n\e\ol X}{n\e X}$ and so we have $\Sx_{k+2}\HY\w\RFN$.

    For the converse, note that for $\Pi_k\HY\cBI$ it suffices to consider \zzz
orderings on \w. Let $\prec\,=\myset{(m,n)}{(m,n)\in\m\prec}$ be such an ordering, \zzz
and suppose
\[
\Ex\text{ parameters }\ppl \Ax n\,\big((\Ax m\prec n\,\px m)\imp\px n\big)\ax\nx\Ax n\,\px n\ppr
\label{Chapter6:eq:6}\tag{6}
\]
where \px\ is in $\Pi_k$. Then by $\Sx_{k+2}\HY\w\RFN$ there is an \w-model of \eqref{Chapter6:eq:6}, \zzz
and so $\wf(\prec)$ is false.\done

    With regard to 6.6, in certain circumstances we need only look \zzz
at first-order sentences. Call an admissible ordinal \aaa\ \ul{self-definable} if \zzz
for no $\bbb\ltx\aaa$, $L_\bbb\prec_1 L_\aaa$.

\thm{6.11 Corollary} Let \s\ be a model of \s\ + $\Sx_1$-overspill of standard \zzz
ordinal \aaa, and suppose \aaa\ is self-definable. Let $\T\e L_\aaa$ be any \zzz
theory of $\sss L_{L_\aaa}$ extending \s. Then there exists a proper $0$-elementary \zzz
substructure of \A\ which includes \aaa\ and is a model of \T\ iff \A\ is a \zzz
model of those first-order $\Sx_1$ sentences \tx\ for which $\T\modelx\!\!\!_\aaa\,\tx$, i.e., those \zzz
\tx\ which are true in all models of \T\ with standard ordinal $\gex\,\aaa$.\done

\noindent
This follows easily from the following well-known result of Kripke and \zzz
Platek:
\[
\parbox{4in}{\aaa\ is self-definable iff for each $a\e L_\aaa$, there exists a $\Sx_1$
formula $\tx x$ such that for any structure \B\ end-extending $L_\aaa$, $\B\modelx\Ex!x\,\tx x\ax\tx a$.}
\]
For a proof, see V.7.8 of Barwise\,\cite{Barwx75}. For readers who, like myself, \zzz
know very little generalized recursion theory, I shall give a brief \zzz
indication of the extent of the self-definable ordinals. The least admissible \zzz
ordinal which is \ul{not} self-definable is greater than the first recursively \zzz
inaccessible, the first recursively Mahlo, the first recursively hyper-Mahlo, \zzz
etc.; see Cenzer\,\cite{Cenz74}. The largest self-definable ordinal is the \zzz
least stable ordinal, that is, the least \aaa\ such that $L_\aaa\prec_1 L$.

    As our final result, we shall consider the analogue of 2.8 on page \pageref{Chapter2:Corollary28}. Let \T\ be \zzz
a consistent first-order r.e.\ theory extending \s\ + \Infinity, and let
\[
(\Gamma_1,\,\Gamma_2)_\w=\myset{\godel\px\e\Gamma_1}{\text{\px\ is $\Gamma_2$-conservative over \T\ with the \w-rule}}\,.
\]
Say \T\ is $\Pi^1_1$-sound if for all $\Pi^1_1$ sentences \px, if $\T\provex\!\!\!_\w\,\px$, then \px\ \zzz
is true.

    \refstepcounter{thisisdumb}\label{Chapter6:6.12}
    The following result is not, as 2.8 was, the best possible, and \zzz
further work remains to be done. (But see page \pageref{Chapter7:Chapter6Addendum}.)

\thm{6.12 Theorem} Let $k\e\w$ and let \T\ be a consistent r.e.\ extension of \s.
\begin{mylist}
\item If \T\ extends $\Sx_{k+1}$-\Collection, a sentence \px\ is $\Sx_{k+1}$-conservative \zzz
($\Pi_{k+2}$-conservative, respectively) over \T\ with the \w-rule iff all $L_{\smash{\w^{\text{CK}}_1}}$-saturated \zzz
models of \T\ have (for any element $x$) a $k$-elementary substructure (containing $x$) \zzz
which is a model of \T\ + \nx\px.
\item If $k\ge 0$ and \T\ extends $\Sx_{k+1}$-\Collection, then
\[
(\Pi_{k+1},\,\Sx_{k+1})_\w\text{ (if $k\geX 1$), and }\ (\Sx_{k+2},\,\Pi_{k+2})_\w
\]
are $\Pi^0_1$ in Kleene's \sss O, and are complete for this class of sets.
\item If either \T\ is not $\Pi^1_1$-sound and includes $\Sx^1_1$-\BI\ or if \T\ is strong \zzz
enough to prove that every $\Pi^1_1$ formula is equivalent to a $\Sx_1$ formula, then
\[
(\Pi_1,\,\Sx_1)_\w\ \text{ and }\ (\Pi_1,\,\Dx_1(\T))_\w
\]
are also complete for this class.
\end{mylist}
Proof: The proof of (i) is clear if we can show that for any sentence \tx\ \zzz
consistent with \T\ in \w-logic there exists a $L_{\smash{\w^{\text{CK}}_1}}$-saturated model of \zzz
\T\ + \tx. But this is a well-known result, closely related to the Gandy \zzz
Basis Theorem, and is due to Ressayre\,\cite{Ress77} (and perhaps others).

    For (ii), we first need a result concerning semi-representability.

\thm{6.13 Lemma} Fix some recursive set $W$, and suppose that for all $e\e\w$,
\[
\prec_e\,=\myset{\seq{x,\m y}}{\seq{x,\m y,\m e}\e W}
\]
is a linear ordering and that $\myset{\prec_e}{e\e\w}$ contains all primitive recursive \zzz
well-orderings.
\begin{mylist}
\item There is a $\Sx^1_1$ formula $\yx\mx x$ such that
\[
\prec_n\text{ is well-founded }\ \text{iff }\ \T\provex\!\!_\w\,\yx\ol n\,.
\]
\item If \T\ is $\Pi^1_1$-sound, we may choose \yx\ above to be $\Pi^1_1$.
\item If \T\ is not $\Pi^1_1$-sound and if \T\ extends $\Sx^1_1$-\BI, we may choose $\yx\e\Dx^1_1(\T)$.
\item The set
\[
\Myset{e}{\Ax n\,\big( \wf(\prec_n)\imp\wf(\prec_{(e,n)})\,\big) }
\]
is $\Pi^0_1$ in \sss O, and is complete for this class of sets, and hence so is
\[
\Myset{\yx}{\Ax n\,\big( \wf(\prec_n)\imp(\T\provex\!\!_\w\,\yx\ol n)\,\big) }
\]
where \yx\ ranges over $\Sx^1_1$ formulae, or over $\Dx_1(\T)$ formulae if \T\ is \zzz
not $\Pi^1_1$-sound and extends $\Sx^1_1$-\BI.
\end{mylist}
Proof of 6.13: (i)\, Choose \yx\ so that \T\ proves for all $n\e\w$,
\[
\yx\ol n\equiv\Ex\tau\,\ppl\m \rk(\tau)\mx=\m\prec_{\ol n}\m\ax\m\Ax\sx\e\tau\,\big(\T^{\frac12\ast\sx}\ax(\nx\yx\ol n)^{\frac12\ast\sx}\big)\m\ppr \,.
\]
Then $\prec_n$ well-founded implies $\T\provex\!\!_\w\,\nx\yx\ol n\imp\yx\ol n$, which implies $\T\provex\!\!_\w\,\yx\ol n$. \zzz
And $\T\provex\!\!_\w\,\yx\ol n$ implies $\prec_n$ well-founded.

\noindent
(ii)\, If \T\ is $\Pi^1_1$-sound, we merely consider the $\Pi^1_1$ predicate
\[
``\m\T\provex\!\!_\w\,\yx\m\ol n\m\,"
\]
where \yx\ is as above. Alternatively, consider the predicate ``$\m\T\provex\!\!_\w\,\wf(\prec_n)$\,''.

\noindent
(iii)\, Let $\prec$ be some recursive non-well-ordering which \T\ proves is a \zzz
well-ordering. Now $\Sx^1_1$-\BI\ implies that there is a hyperarithmetic hierarchy \zzz
$H$ along $\prec$. Let \yx\ be such that \T\ proves for all $n\e\w$,
\[
\yx\ol n\equiv\Ex H\,\bigl(\text{$H$ is a hyperarithmetic hierarchy along $\prec$ and }
\Ex i\e\w\,\Ex\tau\e H_i\,\ppl \rk(\tau)\mx=\prec_{\ol n}\ax\Ax\sx\e\tau\,(\T^{\frac12\ast\sx}\!\ax(\nx\yx\ol n)^{\frac12\ast\sx})\m\bigr)\m\ppr \,.
\]
We claim that $\T\provex\!\!_\w\,\yx\ol n\,$ iff $\prec_n$ is a well-ordering. The only novel \zzz
point here is to see that if \T\ + $\yx\ol n$ is consistent in \w-logic, then \zzz
for any ordinal $\aaa\ltx\w^\text{CK}_1$, there exists a tree $\tau\e\rrm{HYP}$ of rank \aaa\ which \zzz
$\frac12\ast$fulfils \T\ + $\yx\ol n$; this is clear from the proof of 6.5.

\noindent
(iv)\, Finally, (iv) is obtained immediately from the above by looking at \zzz
the definition of the dual class, the collection of sets r.e.\ in \sss O.

    This concludes the proof of Lemma 6.13.\done

\noindent
Continuation of proof 6.12: Fix $k\ge 0$, Let $\yx\mx x$ be $\Sx_1$, and consider the \zzz
fixed point
\[
\T\provex\!\!_\w\,\Phi\equiv\Ax n\,
\PPL\Ex\tau\,\ppl \rk(\tau)\mx=\m\prec_n\m\ax\m\Ax\sx\e\tau\,(\T^\sx\ax\Tr\Sx_k\m^\sx\ax\Phi^\sx)\imp\yx\mx n\ppr \m\PPR
\label{Chapter6:eq:7}\tag{7}
\]
or, alternatively, the fixed point
\[
\T\provex\!\!_\w\,\Phi\equiv\nx\Ax n\,
\PPL\m\Ax m\,\ppl \text{``$\prec_m$ embeddable into $\prec_n$''}\imp\yx \mx m\ppr \imp
\Ex\tau\,\ppl \rk(\tau)\mx=\prec_n\m\ax\m\Ax\sx\e\tau\,(\T^\sx\ax\Tr\Sx_k\m^\sx\ax\Phi^\sx)\m\ppr \m\PPR\,.
\]
Consider:
\begin{mylist}
\renewcommand{\labelenumi}{(\alph{enumi})}
\itemsep0ex
\item $\Ax n\e\w\,\ppl \wf(\prec_n)\imp(\T\provex\!\!_\w\,\yx\ol n)\m\ppr $
\item $\Phi$ is $\Sx_{k+1}  $-conservative over \T\ + \w-rule
\item $\Phi$ is $\Sx^1_1    $-conservative over \T\ + \w-rule
\item $\Phi$ is $\Dx^1_1(\T)$-conservative over \T\ + \w-rule
\end{mylist}
\mypagebreak{[2]}
We claim that (a), (b) and (c) are equivalent.

    First suppose (a), and let \A\ be any model of \T. We shall show \zzz
that there is a $k$-elementary substructure of \A\ which is a model of \T\ + $\Phi$. \zzz
If $\A\modelx\Phi$, we are done. If not, choose some witnesses $n,\tau\e A$ for $\nx\Phi$. \zzz
Now $\prec_n$ is necessarily non-well-founded, and so $\tau$ has (in the real \zzz
world) an infinite branch, from which we may obtain the required structure. \zzz
Hence (b).

    The implication (b)$\imp$(c) is trivial, and (c)$\imp$(a) follows \zzz
from 6.5.ii. By the lemma, this gives the first part of (ii). For the \zzz
second, we merely alter \eqref{Chapter6:eq:7} to ensure that there exists $\tau$ containing \zzz
any arbitrary element.

    For (iii), we simply choose \yx\ to be $\Dx_1(\T)$ in \eqref{Chapter6:eq:7}. Then we have \zzz
that (a), (b), and (d) are equivalent.

    This completes the proof of 6.12.\done
   \cleardoublepage
\thispagestyle{plain} \sectioncentred{The Paris-Harrington Statement}

    Since the original definition of fulfilment was motivated by the \zzz
Paris-Harrington statement, PHS, we thought that it would be appropriate \zzz
to conclude this thesis with an exposition of this. For newcomers, we \zzz
mention that PHS (defined below) is a natural combinatorial sentence which \zzz
is true but not provable in \PA: in fact, the instance
\[
\Ax e\,\Ex n\;n\phs(e+1)^e_3
\]
is not provable in \PA, even with the set of true $\Pi_1$ sentences as additional \zzz
axioms. Furthermore, if we let
\[
\sx(e,c)=\mu n\,.\,n\phs(e+1)^e_c
\]
then for each function $f$ provably recursive in \PA\ there exists an $e$ \zzz
such that
\[
f(x)\ltX \sx(e,x)
\]
for all $x$ (although for each $e$ the function $\lambda x\m.\m\sx(e,x)$ \ul{is} provably \zzz
recursive), and
\[
f(x)\ltX \sx(e,3)
\]
for all large $x$.

    Let $k,e,c$ be (non-negative) integers and let $X$ be a (finite) set \zzz
of integers. Say the partition relation
\[
X\phs(k)^e_c
\]
holds if $|X|\geX e$ and for each function $f\mymap{[X]^e}c$, where $[X]^e$ is \zzz
the set of increasing $e$-tuples from $X$ and where $c=\{0,1,\dots,c-1\}$, \zzz
there exists a \ul{large} subset $Y$ of $X$ of cardinality at least $k$ which \zzz
is \ul{homogeneous} for $f$, where $Y$ is \ul{large} if $|Y|\geX\min Y$ and where $Y$ \zzz
is \ul{homogeneous} for $f$ if $f$ is constant on $[Y]^e$. This is a primitive recursive \zzz
relation (in fact elementary), and so it is expressible by a formula in \zzz
the language of arithmetic which is $\Dx_1$ in \PRA. The Paris-Harrington \zzz
statement is
\begin{flalign*}
&\text{PHS:} &&\Ax k,e,c\,\Ex n\;n\phs(k)^e_c\,.\quad&&
\end{flalign*}
\indent The plan of this chapter is as follows. First the PHS is shown \zzz
to be true, and then various independence results are established, with the \zzz
PHS being shown to be equivalent to the $\Sx_1$ Reflection Principle for \PA. \zzz
(Only this latter fact makes any use of the notion of fulfilment.) Then \zzz
we consider the rate of growth of the function \sx.

    The exposition given below is an amalgam from many sources, but all \zzz
the main ideas are due to J.~Paris and L.~Harrington. The work of \zzz
L.~Kirby, G.~Mills and J.~Paris\,\cite{KirbyMillsParis79} and of J.~Ketonen and R.~Solovay\,\cite{Keto79} has \zzz
provided sharp results, namely that
\[
\PA^-_\text{ex}+\Sx_{e+1}\HY\Collection\provex\Ax k\,\Ex n\;[k,n]\phs(e+2)^{e+1}_{\ol c}
\label{Chapter7:eq:1}\tag{1}
\]
for each integer $c$, and that
\[
\text{true $\Pi_1$ sentences}+\Sx_{e+1}\HY\Collection\nprovex\Ax c\,\Ex n\;n\phs(e+2)^{e+1}_c.
\refstepcounter{thisisdumb}\label{Chapter7:eq:2}\tag{2}
\]
The proofs below give \eqref{Chapter7:eq:2} but not \eqref{Chapter7:eq:1}. On the other hand, they are much \zzz
simpler than those of the above-mentioned works, replacing complex model-theoretic \zzz
and combinatorial arguments by na\"ive ones.

\thm{7.1 Theorem} PHS is true.

\noindent
Proof: Suppose not, and let $k,e,c$ be chosen to violate the theorem. Call \zzz
$f$ \ul{good} if $f\mymap{[n]^e}c$ for some $n$ and $f$ has no large homogeneous subset \zzz
of cardinality $\geX\,k$. Let
\[
T=\Myset{g}{g\mymap{[n]^e}c \text{ for some $n$ and $g\subseteq f$ for some good $f$}}\,.
\]
Then $T$ ordered by inclusion is a finitely branching tree, and by our \zzz
initial assumption, it is infinite. By K\"onig's Lemma, it has an infinite \zzz
branch $B$. Let $F=\cupx B$. Then $F\mymap{[\w]^e}c$, and by the infinite version \zzz
of Ramsey's Theorem, there exists an infinite $Y\subseteq\w$ which is homogeneous \zzz
for $F$. Let $X\subseteq Y$ be a finite large set of cardinality $\geX\,k$. Then \zzz
$F\up X\subseteq f$ for some good function $f$. But then $X$ is homogeneous for $f$, \zzz
which contradicts the goodness of $f$.\done

    In fact, one can show that for each $e$, where $\PA^-_\text{ex}$ is as given on \pageref{Chapter1:DefinitionPAex},
\[
\PA^-_\text{ex}+\Sx_{2e+4}\HY\Collection\provex\Ax k,c\,\Ex n\;[k,n]\phs(e+2)^{e+1}_c.
\label{Chapter7:eq:3}\tag{3}
\]
For in the above proof, we can choose the branch $B$ so that the graph of $F$ is \zzz
definable: $\Dx_3$ in fact. By a direct analysis of a suitable proof of the \zzz
infinite Ramsey's Theorem (see Jockusch\,\cite{Jock72}, page 275), it can be seen \zzz
in \PA\ that $F$ has an infinite definable (it can be shown to be $\Dx_{2e+2}$) \zzz
homogeneous subset. The contradiction follows as above. Thus for each \zzz
$e$ the function $\lambda x\,.\,\sx(e,x)$ is provably recursive in \PA.

    The above argument will not work if $e$ is a free variable, for \zzz
Jockusch\,\cite{Jock72} proves, for instance, that for each $e\geX 2$, there exists a \zzz
\ul{recursive} partition of $[\w]^e$ into two classes with no infinite $\Dx_e$ \zzz
homogeneous set.

\thm{7.2 Lemma} For all $e\geX 1$
\[
\PA^-_\text{ex}+\Sx_e\HY\Collection\nprovex\Ax c\,\Ex n\;n\phs(e+2)^{e+1}_c.
\]
Proof: (We shall in fact prove more.) Let \M\ be a nonstandard model of \zzz
$\PA^-_\text{ex}$. Suppose that for some nonstandard $c\e M$,
\[
\Ex n\;n\phs(e+2)^{e+1}_c.
\]
Fix $n$ to be the least witness. Let $2^{b+2}$ be any nonstandard power of $2$ \zzz
less than  $c$, and let $\seq{\px_i}$ be any natural listing of all $\Dx_0$ formulae \zzz
whose free variables are among $v_0\m,\m v_1,\dots,v_e$. Define $F\mymap{[n]^{e+1}}c$ by
{
\setlength{\belowdisplayskip}{1ex}
\setlength{\abovedisplayskip}{1ex}
\begin{align*}
F(\xvect v)&=
\begin{cases}
2v_0+1\,,                            &\text{if }v_0\ltX 2^b,\\
\sum_{i\lt b}2^{i+1}F_i(\xvect v)\,, &\text{otherwise,}
\end{cases}
\intertext{where \vect v is (the code of) an $e+1$-tuple from \n, and}
F_i(\xvect v)&=
\begin{cases}
1\,, &\text{ if }\Sat_{\Dx_0}(\m\godel{\px_i\!},\xvect v)\,,\\
0\,, &\text{ otherwise,}
\end{cases}
\end{align*}
}
where $\Sat_{\Dx_0}$ is the usual satisfaction predicate for $\Dx_0$ formulae.

    Since $F$ can be defined internally, there exists in \M\ a large \zzz
subset $\{c_0\ltx c_1\ltx\dots\}$ of \n\ of cardinality greater than $e+1$ \zzz
which is homogeneous for $F$. Since $F(c_0\m,\dots,c_e)=F(c_1,\dots,c_{e+1})$, it \zzz
must be the case that $c_0\geX 2^b$. Thus the set is of nonstandard cardinality, \zzz
and so we may consider the initial segment $I=\cupx_{i\in\w}\m c_i$. If we can show \zzz
that $I$ is a model of $\PA^-_\text{ex}$ + $\Sx_e$-\Collection, we are done, for the partition \zzz
relation is then absolute between $I$ and $M$, and obviously $n\notE I$. We \zzz
proceed as follows. For a $t$-tuple $\vect i=\seq{i_0\m,\dots,i_{t-1}}$, let $c_{\subvect i}=\seq{c_{i_0}\m,\dots,c_{i_{t-1}}}$.

\noindent
(i)\, First note that for any $\Dx_0$ formula $\px\vect v$ with $t\lex e+1$ free \zzz
variables and for any two increasing $t$-tuples $0\leX\vect i,\,\vect j\leX c_0+t-e-2$, \,$\px c_{\subvect i}\equiv\px c_{\subvect j}$ \zzz
holds in \M. If $t\ne 0$, this is because $F(c_{\subvect i}\m_\smile\vect z)=F(c_{\subvect j}\m_\smile\vect z)$ for \zzz
$\xvect z=\seq{c_{c_0+t-e-1},\m c_{c_0+t-e},\dots,\m c_{c_0-1}}$ (where $\xvect z=\seq{}$ when $t\eqx e+1$).

\noindent
(ii)\, $I$ is closed under addition and multiplication, as follows. \zzz
$c_0+i\leX c_i$ for all $i\leX c_0-e$, and so $2c_o-e\leX c_{c_0-e}$. As $c_0\gtX 2e$, \zzz
we have $\frac32c_0\ltX c_{c_0-e}$, and hence by (i), $\frac32c_0\ltX c_1$. So \zzz
$(\frac32)^ic_0\ltX c_i$ for all $i\leX c_0-e$ and so $(\frac32)^{c_0/2}c_0\ltX c_{c_0-e}$. Thus $(c_0)^2\ltX c_{c_0-e}$, and \zzz
by (i) again, $(c_i)^2\ltX c_{i+1}$.

    \refstepcounter{thisisdumb}\label{Chapter7:DefinitionNabla}
    Define the class $\nabla_k$ of formula inductively as follows. Let $\nabla_0=\Dx_0$, \zzz
and let $\nabla_{k+1}$ be the closure of $\nabla_k\cup\myset{\Ev x\m\tx}{\tx\e\nabla_k}\cup\myset{\Av x\m\tx}{\tx\e\nabla_k}$ \zzz
under conjunction and disjunction. (NB: we do not close under bounded \zzz
quantification.) For each $\nabla_k$ formula $\px(\xvect x)$, define a $\Dx_0$ \zzz
formula $\px^\Delta(\xvect x,y_0\m,\dots,y_{k-1})$ inductively as follows. If \px\ is \zzz
$\Dx_0$, let $\px^\Delta=\px$;\, let $(\px\vx\tx)^\Delta=\px^\Delta\vx\tx^\Delta$, and similarly for conjunction; \zzz
and let $(\ssf Q \xvect v\m\px(\xvect x,\xvect v)\m)^\Delta=\ssf Q \xvect v\!\ltx y_0\;\px^\Delta(\xvect x,\xvect v,\m y_1,\dots,y_k)$ for
$\px\e\nabla_k$ and for $\ssf Q$ $\Ex$ or $\Ax$.

    Let $\vecth{d_k}$ be the $k$-tuple $\seq{c_{c_0-e-k},\dots,c_{c_0-e-1}}$.

\noindent
(iii)\, We claim that if $\px\vect x\ve\nabla_k$ with $k\ltX e$, then for each $l$, each \zzz
increasing $k$-tuple $\vect i$ with $l\ltX \vect i\!\ltX c_0-e$, and for all parameters $\vect x\!\ltX c_l$,
\[
\px^\Delta(\xvect x,\m c_{\subvect i})\equiv\px^\Delta(\xvecth x,\xvecth{d_k})\,.
\]
The proof is by induction on the complexity of \px. The only non-trivial \zzz
case is when $\px\xvect x=\ssf Q\xvect y\,\yx(\vect x,\xvect y)$ with $\yx\e\nabla_{k-1}$. Suppose the claim \zzz
is false: for some $l$ and $k$-tuple $l\ltX\xvect i\ltX c_0-e$,
\begin{flalign*}
&                          & &\Ev x\!\ltx c_l\,\big(\px^\Delta(\xvect x,c_{\subvect i})\nequiv\px^\Delta(\xvecth x,\vecth{d_k})\m\big)\,,&&\\
&\text{that is,}\hidewidth &  \Ev x\!\ltx c_l\,\big(
                       \ssf Q\xvect y&\!\ltx c_{i_0}\,\yx^\Delta(\xvect  x,\vect  y,c_{\subvect j})\nequiv
                       \ssf Q\xvect y \!\ltx c_{d'} \,\yx^\Delta(\xvecth x,\vecth y,\vecth{d_{k-1}})\m\big)\,,
\end{flalign*}
where $\vect i=\seq{i_0}_\smile\vect j$ and $d'=c_0-e-k$. The inductive hypothesis gives \zzz
that for all $\vect x,\xvect y\ltX c_{i_0}$,
\begin{flalign*}
&                        & \yx^\Delta(\xvect x,\vect y,c_{\subvect j})\equiv\yx^\Delta(&\xvecth x,\vecth  y,\vecth{d_{k-1}})&&\\
&\text{so if }\ssf Q\eqx\Ex,\hidewidth & \Ev x\!\ltx c_l\,\big(
                             \nx\Ev y\!\ltx c_{i_0} \,\yx^\Delta(\xvecth x,\vecth y,\vecth{d_{k-1}})\;&\ax
                                \Ev y\!\ltx c_{d'}  \,\yx^\Delta(\xvecth x,\vecth y,\vecth{d_{k-1}})\m\big)\,,\\
&\text{or if }\ssf Q\eqx\Ax,\hidewidth & \Ev x\!\ltx c_l\,\big(\,
                           \Av y\!\ltx c_{i_0}\,\yx^\Delta(\xvecth x,\vecth y,\vecth{d_{k-1}})\,\ax
                        &\;\nx\Av y\!\ltx c_{d'}\,\yx^\Delta(\xvecth x,\vecth y,\vecth{d_{k-1}})\m\big)\,.
\end{flalign*}
Say $\ssf Q\eqx\Ex$. Since the number of $c$'s is $3+(k-1)\leX e+1$, (i) gives
\[
\Ev x\!\ltx c_0\,\big(\nx\Ev y\!\ltx c_2\,\yx^\Delta(\xvecth x,\vecth y,\vecth{d_{k-1}})\ax
                         \Ev y\!\ltx c_3\,\yx^\Delta(\xvecth x,\vecth y,\vecth{d_{k-1}})\m\big)\,.
\]
Let $f(z)=\big\lfloor\!\sqrt[\leftroot{0}\uproot{3}s]{z-e-k-2}\,\big\rfloor$, where $s$ is the length of \vect x. By (ii), $c_0\ltX f(c_1)$, and so we may change \zzz
the latest bound on \vect x from $c_0$ to $f(c_1)$. So by (i), for all $i$ \zzz
with $0\ltX i\ltX d'$,
\[
\Ev x\!\ltx f(c_0)\,\big(\nx\Ev y\!\ltx c_i    \,\yx^\Delta(\xvecth x,\vecth y,\vecth{d_{k-1}})\ax
                            \Ev y\!\ltx c_{i+1}\,\yx^\Delta(\xvecth x,\vecth y,\vecth{d_{k-1}})\m\big)\,.\label{Chapter7:eq:*}\tag{$\ast$}
\]
Consider in \M\ a map from $\myset{i}{0\ltx i\ltx d'}$ to $\myset{\vect x}{\vect x\ltx f(c_0)}$, where for each $i$ we choose a \zzz
witness $\vect x$ to the body of \eqref{Chapter7:eq:*}. The size of the codomain is \zzz
$f(c_0)^s\leX\,c_0-e-k-2\eqx d'-2$, less than that of the domain. So by the pigeon-hole \zzz
principle in \M, there exist $i,j$ with $0\ltX i\ltX j\ltX d'$ with the same witness $\vect x$. That is,
\[
\nx\Ev y\!\ltx c_i    \,\yx^\Delta(\xvecth x,\vecth y,\vecth{d_{k-1}})\ax
   \Ev y\!\ltx c_{i+1}\,\yx^\Delta(\xvecth x,\vecth y,\vecth{d_{k-1}})\ax
\nx\Ev y\!\ltx c_j    \,\yx^\Delta(\xvecth x,\vecth y,\vecth{d_{k-1}})\ax
   \Ev y\!\ltx c_{j+1}\,\yx^\Delta(\xvecth x,\vecth y,\vecth{d_{k-1}})\,.
\]
But the middle two conjuncts are contradictory. The $\ssf Q\eqx\Ex$ case is identical, except that the positions of the negations are changed. This proves the claim.

\noindent
(iv)\, Using (iii), by an easy induction on complexity, we see that for $k\ltX e$, for any $\nabla_k$ formula $\px\vect x$ \zzz
and any $\vect x\ve I$,
\[
I\modelx\px\vect x\text{\ \ \ iff \ \ \ }\px^\Delta(\xvecth x,\vecth{d_k})\,.\label{Chapter7:eq:4}\tag{4}
\]
Also for any $\nabla_e$ formula $\px\vect x$ and any $\vect x\ve I$,
\[
I\modelx\px\vect x\text{\ \ \ iff \ \ \ }\Ax\text{ large }j\e\w\;\px^\Delta(\xvecth x,c_j,\vecth{d_{e-1}})\,.\label{Chapter7:eq:4x}\tag{$4'$}
\]
(v)\, That $\Sx_e$-\Collection\ holds in $I$ now follows by a simple application of \zzz
underspill: suppose
\[
I\modelx\Av x\!\ltx a\,\Ev y\m\tx\vect x\xvect y
\]
where $a\e I$ and where \tx\ is a $\Pi_{e-1}$ formula which may have parameters \zzz
from $I$. Then for any nonstandard $i\ltX c_0+1-2e$,
\[
\Av x\!\ltx a\,\Ev y\!\ltx c_i\,\tx^\Delta(\xvecth x,\xvecth y,\vecth{d_{e-1}})\,.
\]
By underspill, there exists a standard $i$ for which this holds, and by \eqref{Chapter7:eq:4} \zzz
we are done.

\noindent
(vi)\, We have completed the proof of the lemma, but we shall also show \zzz
that parameter-free $\nabla_e$-\Foundation\ holds in $I$, as follows. Let $\px\mx x$ be $\nabla_e$, \zzz
with no other parameters, and suppose that $I\modelx\Ex x\m\px\mx x$, and that we wish to find a least witness. Then
\begin{flalign*}
&& I\modelx\Ex x\m\px\mx x\;\;\text{iff\ \ } & \Ex i\;\Ex x\ltx c_i\,\big(I\modelx\px\m x\big)\,\text{, $i$ standard}&&&\\
&& \text{iff\ \ } & \Ex i\;\Ex x\ltx c_i\,\Ax\text{ large }j\;\px^\Delta(x,c_j,\vecth{d_{e-1}})\text{, $i,j$ standard, by \eqref{Chapter7:eq:4x}}\\
&& \text{iff\ \ } & \Ex x\ltx c_i\,\px^\Delta(x,c_j,\vecth{d_{e-1}})\text{, for a standard $i$ and non-standard $j$ (or for $\Rightarrow$, take $j\eqx i+1$)}\\
&& \text{iff\ \ } & \Ex x\ltx c_i\;\px^\Delta(x,c_j,\vecth{d_{e-1}})\text{ for \ul{all} $i,j\e M$, with $0\lex i\ltx j\ltx d'\eqx c_0-2e+1$, by (i)}\\
&& \text{iff\ \ } & \Ex x\ltx c_0  \;\px^\Delta(x,c_2,\vecth{d_{e-1}})\text{, by (i)}\\
&& \text{which implies\ \ } & \Ex x\ltx c_1/4\;\px^\Delta(x,c_2,\vecth{d_{e-1}})\text{, by (ii)}\\
&& \text{iff\ \ } & \Ex x\ltx c_0/4\;\px^\Delta(x,c_i,\vecth{d_{e-1}})\text{, by (i) again, for all $i$ with $0\ltX i\ltX d'$.}
\end{flalign*}
Consider in \M\ the map $i\mapsto\mu x\ltx c_0/4\,.\,\px^\Delta(x,c_i,\vecth{d_{e-1}})\m$, $0\ltX i\ltX d'$. Since $d'-1\gtX c_0/4$, by \zzz
the pigeon-hole principle there exist $x_0\ltX c_0/4$ and $i,j$ with $0\ltX i\ltX j\ltX d'$ \zzz
such that
\[
                  \px^\Delta(x_0,c_i,\vecth{d_{e-1}})\ax
\Ax y\ltx x_0\,\nx\px^\Delta(y,  c_i,\vecth{d_{e-1}})\ax
                  \px^\Delta(x_0,c_j,\vecth{d_{e-1}})\ax
\Ax y\ltx x_0\,\nx\px^\Delta(y,  c_j,\vecth{d_{e-1}})\,.
\]
And so
\[
\Ex x\ltx c_i/4\,\big(
                 \px^\Delta(x,c_i,\vecth{d_{e-1}})\ax
\Ax y\ltx x \,\nx\px^\Delta(y,c_i,\vecth{d_{e-1}})\ax
                 \px^\Delta(x,c_j,\vecth{d_{e-1}})\ax
\Ax y\ltx x \,\nx\px^\Delta(y,c_j,\vecth{d_{e-1}})\m\big)\,.
\]
By (i), this must hold for all $i\ltX j\ltX d'$. Setting $i\eqx 0$ shows the given map is constant. Thus \zzz
$\px^\Delta(x_0\m,c_i,\vecth{d_{e-1}})$ holds for all $i\ltX\w$, and so $I\modelx\px(x_0)$. If $I\modelx\px(y)$ \zzz
for some $y\ltX x$, then for some $i\ltX c_0$ (in fact, for all sufficiently \zzz
large $i\ltX\w$), $\px^\Delta(y,c_i,\vecth{d_{e-1}})$ holds, a contradiction.\done

\thm{7.3 Lemma} For all $e\geX 0$,
\begin{mylist}
\item $\Sx_{e+1}$-\Collection\ is a $\Pi_{e+2}$-conservative extension of $\Sx_e$-\Induction, and
\item $\Sx_e$-\Induction\ is a $\Sx_{e+2}$-conservative extension of $\Sx_e$-\Collection\ \zzz
plus parameter-free $\nabla_e$-\Foundation\ (with $\nabla_e$ defined on page \pageref{Chapter7:DefinitionNabla}), all over $\PA^-$.
\end{mylist}
Proof: The first part of the lemma is just Application (xiv) on page \pageref{Chapter4:Applicationxiv}, so let us \zzz
consider the second. Let \M\ be any model of $\PA^-$ + $\Sx_e$-\Collection\ and parameter-free \zzz
$\nabla_e$-\Foundation, and let $N\subseteq M$ consist of those elements $a$ of $M$ for \zzz
which there is a $\Sx_{e+1}$ formula $\yx\mx x$ such that
\[
\M\modelx\Ex!x\,\yx\mx x\ax\yx\mx a\,.
\]
Then $N\prece_{e+1}\M$, as follows. $N$ is clearly closed under addition and \zzz
multiplication. Suppose $\M\modelx\Ex x\m\yx\mx$, where \yx\ is $\Pi_e$. Then by \zzz
parameter-free $\nabla_e$-\Foundation,
\[
\Ex!x\e M\;\big(\M\modelx\yx\mx x\ax\Ax y\ltx x\,\nx\yx\mx y\,\big)\,.
\]
By $\Sx_e$-\Collection, this latter formula is equivalent to a $\Sx_{e+1}$ formula, and \zzz
so \yx\ has a witness in $N$. Thus $N\prece_{e+1}\M$.

    To complete the lemma, it suffices to show that $N$ is a model of \zzz
$\Sx_e$-\Foundation. Suppose
\[
N\modelx\Ex x\m\px xa\,,
\]
where \px\ is $\Sx_e$ and for notational simplicity we only allow a single \zzz
parameter $a\e N$. Then
\[
\M\modelx\Ex x,u,z\,(\px xu\ax\yx_a zu)\,,\label{Chapter7:eq:5}\tag{5}
\]
where $\yx_a$ is a $\Pi_e$ formula such that
\[
\M\modelx\Ex! u\,\Ex z\,\yx_a zu\ax\Ex z\,\yx_a za\,.
\]
Now by parameter-free $\nabla_e$-\Foundation, consider the least triple $\seq{x_0\m,u_0\m,z_0}$ which is \zzz
a witness to \eqref{Chapter7:eq:5}. Because $u_0$ is unique, $x_0$ and $z_0$ are independent of \zzz
each other, and so $x_0$ is the unique element satisfying
\[
\M\modelx\Ex u,z\,(\px xu\ax\yx_a zu\ax\Ax y\ltx x\,\nx\px yu)\,.
\]
Now by $\Sx_e$-\Collection, this is equivalent to a $\Sx_{e+1}$ formula, and so \zzz
$x_0\e N$. Since $N\prece_{e+1}\M$, $x_0$ is the least witness to \px\ in $N$ also.\done

    From 7.2 with parameter-free $\nabla_e$-\Foundation\ and from 7.3 as the PHS is $\Pi_2$, we immediately have:

\thm{7.4.i Theorem} True $\Pi_1$ sentences + $\Sx_{e+1}\HY\Collection\nprovex\Ax c\,\Ex n\;n\phs(e+2)^{e+1}_c$.\done

    From the proof of 7.2 we see that in $\PA^-_\text{ex}$ the sentence $\Ax c\,\Ex n\;n\phs(e+2)^{e+1}_c$ \zzz
implies that for each $m$ there exists a sequence of length $m$ which \zzz
fulfils the first $m$ axioms of $\Sx_{e-1}$-\Induction, where the most convenient \zzz
notion of fulfilment to use here is the combination of the second notion \zzz
motivated by Skolem functions as given on page \pageref{Chapter1:DefinitionFulfilIII}
with the version of i-fulfilment as given on \zzz
page \pageref{Chapter1:Definitionifulfilment}. The proofs are readily formalized to yield:
\[
\PA^-_\text{ex}+\Ax c\,\Ex n\;n\phs(e+2)^{e+1}_c\provex\RFN_{\Sx_1}(\Sx_{e-1}\HY\Induction)\,,
\]
and so by the remarks following 4.3, we also have
\[
\PA^-_\text{ex}+\Ax c\,\Ex n\;n\phs(e+2)^{e+1}_c\provex\RFN_{\Sx_1}(\Sx_e\HY\Collection)\,.
\]
Moreover, the proofs are uniform in $e$, and so we may obtain
\[
\PA^-_\text{ex}+\Ax e,c\,\Ex n\;n\phs(e+1)^e_c\provex\RFN_{\Sx_1}(\PA)\,.
\]
Furthermore, by formalizing the proof of \eqref{Chapter7:eq:3} we have
(where the dot notation is defined on page \pageref{Chapter2:DefinitionDot})
\[
\PA^-_\text{ex}\provex\Ax e\,\Prx{\PA}(\godel{\,\Ax c\,\Ex n\;n\phs(\dot e+1)^{\dot e}_c\m})\,,
\]
and so we have

\thm{7.4.ii Corollary} $\PA^-_\text{ex}\provex\Ax e,c\,\Ex n\;n\phs(e+1)^e_c\equiv\RFN_{\Sx_1}(\PA)$. \done

    We may take a weaker instance of the left-hand-side of the above \zzz
equivalence as follows. First we need an lemma of Paris and Harrington\,\cite{ParisHarr77}.

\thm{7.5 Lemma} A set $Y\subseteq X$ is homogeneous for $F\mymap{[X]^e}{c}$ iff every \zzz
subset of $Y$ of cardinality $e+1$ is homogeneous for $F$.

\noindent
Proof: Let $\vect x=\seq{x_0\,,\dots,x_{e-1}}$ be the first $e$ elements of $Y$. Pick \zzz
$\vect y=\seq{y_0\m,\dots,y_{e-1}}$ from $Y$ so that $F(\xvect x)\ne F(\xvect y)$ and $y_0+y_1+\dots +y_{e-1}$ \zzz
is minimal. If $i$ is the least index such that $x_i\ne y_i$, then \zzz
$\{x_0\m,\dots,x_i,y_i,\dots,y_{e-1}\}$ is of cardinality $e+1$ but not homogeneous.\done

    Recall the finite version of Ramsey's Theorem, which is provable in \zzz
\PRA:
\[
\Ax k,e,c\,\Ex n\;n\rightarrow (k)^e_c\,.
\]
\thm{7.6 Lemma} If $m\geX 3$, if $m\mx+\mx e\mx+\mx 7\rightarrow (e+1)^e_c$, and if one of the following holds:
\begin{mylist}
\parskip0pt
\item $N\phs(m+1)^m_3$
\item $N\phs(m+2)^m_2$
\item $(m,N)\phs(m+1)^m_2$, where $(m,N)=\myset{x}{m\ltx x\ltx N}$;
\end{mylist}
then $(m,N)\phs(m\mx+\mx e\mx+\mx 7)^e_c$.

\noindent
Proof: We can suppose that $e,c\geX 2$. Let $F\mymap{[N]^2}c$ be given. Let \zzz
$s=m\mx-\mx e\mx-\mx 1$ and $t\eqx 2m+6$. Since $t\!\nrightarrow\!(m\mx+\mx 1)^m_2$ (see comment at end \zzz
of chapter), we can choose $g\mymap{[t]^m}2$ with no homogeneous set of cardinality \zzz
$m\mx+\mx 1$. Define $f\mymap{[N]^{e+1}}2$ by
\[
f(v_0\m,\dots,v_e)=
\begin{cases}
1\, &\text{ if $\{v_0\m,\dots,v_e\}$ is homogeneous for $F$}\\
0\, &\text{ otherwise.}
\end{cases}
\]
Define $G\mymap{[N]^m}3$ as follows. Let $\vect v=\seq{v_0\m,\dots,v_{m-1}}\ltX N$ be given and \zzz
let $i$ be the least $i$ such that $v_i\geX t$, if this exists, and $m$ otherwise. \zzz
If (i) holds let
\vspace{.75ex}
{
\setlength{\belowdisplayskip}{.25ex}
\setlength{\abovedisplayskip}{.25ex}
\begin{align*}
G(\vect v)&=
\begin{cases}
\makebox[2in][l]{$g(\xvect v)$}           &\text{if $i\eqx m$}\\
2                                         &\text{if $i\eqx 1$ or $m-1$}\\
\rrm{parity}(i)                           &\text{if $1\ltX i\ltX m-1$}\\
f(v_0-s,\dots,v_e-s)                      &\text{if $i\eqx 0$.}
\end{cases}
\intertext{If (ii) holds, let}
G(\vect v)&=
\begin{cases}
\makebox[2in][l]{$g(\xvect v)$}           &\text{if $i\eqx m$}\\
\rrm{parity}(i)                           &\text{if $m-e\leX i\ltX m$}\\
\rrm{parity}(i+f(v_i\mx-\mx s,\dots,v_{i+e}\mx-\mx s)) &\text{if $i+e\ltX m$.}
\end{cases}
\intertext{And if (iii) holds, let}
G(\vect v)&=
\begin{cases}
\rrm{parity}(i+f(v_i\mx-\mx s,\dots,v_{i+e}\mx-\mx s)) &\text{if $i+e\ltX m$.}\\
\makebox[2in][l]{$0\text{ or } 1$}        &\text{otherwise.}
\end{cases}
\end{align*}
}

    Let $X\subseteq N$  be homogeneous for $G$ and such that $X$ satisfies \zzz
the appropriate conditions of (i), (ii) or (iii). It is straightforward \zzz
to check that in each case, $|X|\geX\min X\geX t$, and so we can suppose \zzz
$X=\{c_0\ltx c_1\ltx\dots\ltx c_{c_0-1}\}$. Let $Y$ be the first $m\mx+\mx e\mx+\mx 7$ elements of \zzz
$X$, and define $H\mymap{[Y]^e}c$ by
\[
H(v_0\m,\dots,v_{e-1})=F(v_0\mx-\mx s,\dots,v_{e-1}\mx-\mx s)\,.
\]
By our hypothesis, we can choose $Z\subseteq Y$ homogeneous of cardinality $e+1$. \zzz
Now the cardinality of $Z\cup\{c_{m+e+7},\dots,c_{t-1}\}$ is $(e+1)+(t-(m+7-e))=m$, \zzz
and the value of $G$ on this set is $0$, and so $G''[X]^m=\{0\}$. Let \zzz
$X'=\{c_0\mx-\mx s,\dots,c_{c_0-s-1}\mx-\mx s\}$. Obviously $m\mx+\mx e\mx+\mx 7\leX\min X'\leX|X'|$. We claim \zzz
that $X'$ is homogeneous for $F$. For if $ \{x_0\mx-\mx s,\dots,x_e\mx-\mx s\}$ is any subset \zzz
of $X'$ of cardinality $e+1$, then the cardinality of $\{x_0\m,\dots,x_e\}\cup\{c_{c_0-s},\dots,c_{c_0-1}\}$ \zzz
is $e+1+s\eqx m$, and $G(x_0\m,\dots,x_e,c_{c_0-s},\dots,c_{c_0-1})=0$. By Lemma 7.5, \zzz
we are done.\done

\thm{7.7 Theorem} (i)\, If $f\mymap\w\w$ is provably recursive in $\PA^-$ + \zzz
$\Sx_{e+1}$-\Collection, then
\[
f(x)\ltx\sx(e+1,x)
\]
for all large $x$, for $e\geX 1$.

\noindent
(ii)\, If $f\mymap\w\w$ is provably recursive in \PA, then
\[
f(x)\ltx\sx(x,3)
\]
for all large $x$.

\noindent
Proof: (i)\, Suppose $f$ is provably recursive in $\PA^-$ + $\Sx_{e+1}$-\Collection. \zzz
Let \M\ be any proper elementary extension of \bb N, and fix $c$ nonstandard. \zzz
Let $I$ be as in Lemma 7.2, but with a constant for $c$ added to the language. \zzz
Now by the proof of Lemma 7.3, relativized to $c$, we have that \zzz
$f(c)\e I$ whereas $\sx(e,c)\notE I$. That is,
\[
\M\modelx\text{for all large $x$, }f(x)\ltx\sx(e+1,x)\,.
\]
Now the same must be true in \bb N.

\noindent
(ii)\, Let \M\ be any proper elementary extension of \bb N, and fix $e$ \zzz
nonstandard. By 7.6 there is a primitive recursive function $g$ such \zzz
that
\[
\sx(g(e,c),3)\gtx\sx(e,c)\,.
\]
If $c$ is nonstandard, we may construct an initial segment which is a \zzz
model of \PA\ containing $c$ and $e$, and so $g(e,c)$ and $f(g(e,c))$, but \zzz
not $\sx(e,c)$, and so not $\sx(g(e,c),3)$. The theorem follows as \zzz
before.\done

    By Lemma 7.6, we see that \PA\ does not prove $\Ax e\,\Ex n\;n\phs(e+1)^e_3$, \zzz
$\Ax e\,\Ex n\;n\phs(e+2)^e_2$, and $\Ax e\,\Ex n\;(e,n)\phs(e+1)^e_2$. But we have left \zzz
unanswered the question of whether or not \PA\ proves.
\[
\Ax e\,\Ex n\;n\phs(e+1)^e_2\,.
\]
This question does not appear to be answerable using our simple-minded \zzz
techniques. There is an analogous (apparently very difficult) problem \zzz
concerning the ordinary Ramsey partition relation: reasonable bounds \zzz
are not known for the function
\[
\lambda k\,.\,\mu n\,.\,n\rightarrow(k+1)^k_2\,.
\]
One can easily show that
\[
(2k+1)\nrightarrow(k+1)^k_2\,\quad(k\geX 2)\,.
\]
Isbell\,\cite{Isbe69} has proved that $12\nrightarrow(4)^3_2$, and using this we may obtain\footnote{
Proof of $(2k+6)\nrightarrow(k+1)^k_2$ for $k\ge 3$. We can suppose $k> 3$. Let:
\begin{align*}
\text{green} &= \myset{n}{n\ltX 12}\text{, the first twelve numbers}\\
\text{red}   &= \myset{n}{12\leX n\ltX k-3+12}\text{, the next $k-3$ numbers}\\
\text{blue}  &= \myset{n}{k-3+12\leX n\ltX 2(k-3)+12}\text{, and the next $k-3$ numbers}
\end{align*}
Choose f\,\mymap{[12]^3}{2} to witness $12\nrightarrow(4)^3_2$.
Define F\,\mymap{[2k+6]^k}{2} by
\begin{flalign*}
&             &F(X) &=\rrm{parity}(D+E)&&\\
&\text{where} &D(X) &=
\begin{cases}
\text{index of first blue number in natural order, if it exists}\\
k\text{, otherwise}
\end{cases}\\
&\text{where} &E(X) &=
\begin{cases}
f(Y) \text{, if $|X\cap 12|\gex 3$, where $Y$ consists of the three least elements of $X$}\\
0\text{, otherwise}
\end{cases}
\end{flalign*}
Suppose $F$ has a homogeneous set $Z\subseteq 2k+6$ of cardinality $k+1$. Then $|Z\cap 12|\leX 3$.
Hence $Z$ must contain both red and blue numbers. But $F(Z\setminus\{\text{any red number}\})$ cannot
equal $F(Z\setminus\{\text{any blue number}\})$.\done
}
\[
(2k+6)\nrightarrow(k+1)^k_2\,\quad(k\geX 3)\,.
\]
These lower bounds appear to me to be very low, but I cannot improve on \zzz
them. Known upper bounds are also very poor.

\vspace*{.5in}
\thm{Addendum} \label{Chapter7:Chapter6Addendum}

\noindent
Corollary 6.8 on page \pageref{Chapter6:Corollary68} may be strengthened by requiring that \A\ be \zzz
only an \w-model of $\Dx^1_1$-\CA.

\noindent
Using the same technique, we may also greatly \zzz
improve part (ii) of Theorem 6.12: e.g.\ for any \w-consistent r.e.\ theory \zzz
\T\ extending $\Dx^1_1$-\CA\ and for any $k\geX 3$, the sets
\[
(\Ex m\,\Pi_k,\,\Sx_k)_\w\text{ and }(\Sx_k,\,\Pi_k)_\w
\]
are $\Pi^0_1$ in Kleene's \sss O, and are complete for this class of sets. (Here a \zzz
sentence is in ``$\m\Ex m\,\Pi_k$'' if it is a $\Pi_k$ formula prefaced by existential \zzz
numerical quantifiers.)
                 \cleardoublepage
\thispagestyle{plain} \fancyhead[RE,LO]{Bibliography}
\sectioncentredunnumbered{Bibliography}
\renewcommand{\section}[2]{}%

\makeatletter
\renewcommand\@biblabel[1]{19#1}
\makeatother
\newcommand{\myauthor} {\vspace{.5em}\item[]\hspace{-3.0em}}

                                         \cleardoublepage

\ifcourier
\end{small}
} 
\fi

\end{document}